\documentclass[10pt, a4paper]{amsart}

  \usepackage{ifpdf}
  \ifpdf
    \usepackage[pdftex,pagebackref,colorlinks,linkcolor=black,citecolor=blue,urlcolor=blue,a4paper,hypertexnames=true,plainpages=false]{hyperref}
    \hypersetup{pdfauthor={Roman Sauer},
                pdftitle ={Amenable covers, volume, and L2-Betti numbers of aspherical manifolds},
                }
     
  \else
    \usepackage[pagebackref,colorlinks,linkcolor=black,citecolor=blue,urlcolor=blue,a4paper,hypertexnames=true,plainpages=false]{hyperref}
  \fi

\usepackage{amsrefs}
\usepackage[all]{xy}
\usepackage{enumerate}
\usepackage{amssymb}
\usepackage{mathrsfs}

\SelectTips{cm}{}

\newcommand{\cala}{\mathcal{A}}

\newcommand{\bbC}{{\mathbb C}}
\newcommand{\bbN}{{\mathbb N}}
\newcommand{\bbR}{{\mathbb R}} 
\newcommand{\bbZ}{{\mathbb Z}} 

\newcommand{\scrF}{\mathscr{F}}

\newcommand{\scrC}{\mathscr{C}}
\newcommand{\scrH}{\mathscr{H}}
\newcommand{\scrS}{\mathscr{S}}
\newcommand{\scrG}{\mathscr{G}}

\newcommand{\scrR}{\mathscr{R}}
\newcommand{\scrN}{\mathscr{N}}
\newcommand{\scrU}{\mathscr{U}}
\newcommand{\scrV}{\mathscr{V}}
\newcommand{\scrY}{\mathscr{Y}}

\newcommand{\abs}[1]{{\left\lvert #1\right\rvert}}
\newcommand{\absfix}[1]{{\bigl\lvert #1\bigr\rvert}}
\newcommand{\betti}{b^{(2)}}

\newcommand{\id}{\operatorname{id}}
\newcommand{\im}{\operatorname{im}}

\newcommand{\norm}[1]{{\left\lVert #1\right\rVert}}
\newcommand{\pr}{\operatorname{pr}}

\newcommand{\tr}{\operatorname{tr}}
\newcommand{\vol}{\operatorname{vol}}

\newcommand{\defq}{\mathrel{\mathop:}=}

\newcommand{\noqed}{\renewcommand{\qedsymbol}{}}

\newcommand{\symm}{\operatorname{\mathbf{S}}}

\newcommand{\Cn}{\operatorname{C}}
\newcommand{\Dn}{\operatorname{D}}
\newcommand{\Hn}{\operatorname{H}}
\newcommand{\Clf}{\operatorname{C}^{\text{lf}}}

\newcommand{\Flf}{F^{\text{lf}}}
\newcommand{\Glf}{G^{\text{lf}}}

\newcommand{\lz}{\mathrm{L}^\infty(X;\bbZ)}

\newcommand{\lc}{\mathrm{L}^\infty(X;\bbC)}
\newcommand{\lzo}{\mathrm{L}^\infty}
\newcommand{\lzog}{\mathrm{L}^\infty\rtimes\Gamma}
\newcommand{\zg}{{\bbZ\Gamma}}

\newcommand{\sing}{\operatorname{Sing}}

\newcommand{\const}{\operatorname{const}}
\newcommand{\incl}{\operatorname{incl}}
\newcommand{\openstar}{\operatorname{star}}
\newcommand{\diam}{\operatorname{diam}}
\newcommand{\carr}{\operatorname{carr}}

\newcommand{\sign}{\operatorname{sign}}
\newcommand{\sd}{\operatorname{sd}}
\newcommand{\minvol}{\operatorname{minvol}}

\newcommand{\leaf}{\operatorname{\mathscr{L}}}
\newcommand{\mass}{\operatorname{mass}}

\newcommand{\supp}{\operatorname{supp}}
\newcommand{\ricci}{\operatorname{Ricci}}
\newcommand{\ltwo}{{\mathrm{L}^2}}

\newcommand{\bDelta}{\boldsymbol{\Delta}}

\newtheorem{mthm}{Theorem}

\newtheorem{theorem}{Theorem}[section]
\newtheorem{lemma}[theorem]{Lemma}
\newtheorem{setup}[theorem]{Assumption}

\newtheorem*{setup_subsection}{Assumption for Subsection \thesubsection}

\newtheorem*{conjecture_o}{Conjecture}

\newtheorem*{corollary_o}{Corollary}

\theoremstyle{definition}
\newtheorem{definition}[theorem]{Definition}
\newtheorem*{definition_o}{Definition}

\newtheorem{example}[theorem]{Example}
\newtheorem{remark}[theorem]{Remark}
\newtheorem*{remark_o}{Remark}

\numberwithin{equation}{section}

\newcommand{\comment}[1]                      
{
{{\bf Comment: } {\ttfamily #1}}
}

\title[Amenable Covers, Volume and L$^2$-Betti Numbers]{Amenable
  Covers, Volume 
and L$^2$-Betti Numbers of Aspherical Manifolds} 

\author{Roman Sauer}
\address{Department of Mathematics, University of Chicago, 5734 S. University Avenue, Chicago, IL 60637} 
\thanks{The author thanks Clara L\"oh for corrections and 
suggestions for improving the manuscript. He 
acknowledges support of the German science foundation (DFG), 
made through grant SA 1661/1-1}
\email{romansauer@member.ams.org}
\urladdr{www.romansauer.de}

\begin{document}
\frenchspacing
\maketitle

\begin{abstract}
We provide a proof for an inequality between volume and $L^2$-Betti numbers 
of aspherical manifolds for which Gromov outlined a strategy based on 
general ideas of Connes. The implementation of that strategy 
involves measured equivalence relations, 
Gaboriau's theory of $L^2$-Betti numbers of $\scrR$-simplicial complexes, and 
other themes of measurable group theory. Further, we prove new 
vanishing theorems 
for $L^2$-Betti numbers that generalize a classical result 
of Cheeger and Gromov. As one of the corollaries, 
we obtain a gap theorem which 
implies vanishing of $L^2$-Betti numbers of an aspherical manifold when its 
minimal volume is sufficiently small. 
\end{abstract}

\enlargethispage{2\baselineskip}
\tableofcontents

\section{Introduction}\label{sec:introduction} 

\subsection{Discussion of results}\label{subsec:statement of results}
The \textit{minimal volume} $\minvol(M)$
of a smooth manifold~$M$ is the infimum of volumes 
of complete Riemannian metrics on $M$ whose sectional curvature is pinched
between $-1$ and $1$. Questions of significant geometric interest are: 
Is $\minvol(M)>0$? 
If yes, is the infimum attained by a metric of $M$? If yes,
is this metric unique, or does it satisfy certain regularity properties? 
In the seminal paper~\cite{gromov(1982)} Gromov introduced the
notion of \textit{simplicial volume} and initiated the use of
simplicial volume and \textit{bounded cohomology} as basic
tools in the study of the minimal volume. The central fact is that 
the simplicial volume provides a homotopy invariant that bounds the
minimal volume from below. 

In the present paper we transfer classical vanishing theorems for simplicial
volume and bounded cohomology and theorems relating volume and simplicial volume 
of manifolds 
to \textit{$L^2$-Betti numbers} of closed, aspherical manifolds or spaces. 

A manifold~$M$ is called
\textit{aspherical} if its 
universal covering is contractible or, equivalently, $M$ 
is a model of the classifying space~$B\pi_1(M)$ of its 
fundamental group. $L^2$-Betti numbers for regular coverings of closed
Riemannian manifolds were introduced by Atiyah~\cite{atiyah}, 
and their range of definition and application 
was widened over the 
years~\cites{connes,dodziuk,cheeger+gromov,farber,lueck(1998)}. 
By now there is a definition of $L^2$-Betti numbers
$\betti_i(Y;\Gamma)\in [0,\infty]$, $i\ge 0$, of an arbitrary space
$Y$ with the action of a discrete group 
$\Gamma$~\cite{lueck(2002)}*{chapter~6}. 
The most important case for
us is the universal covering $\widetilde{M}$ of a space $M$ with the
natural action of its 
fundamental group $\pi_1(M)$. Here we omit $\pi_1(M)$ in the notation and 
simply write $\betti_i(\widetilde{M})$. 

We provide a proof (Section~\ref{sec:proof of main inequality}) 
of the following inequality 
stated by Gromov in~\cite{gromov(1999)}*{Section 5.33 on p.~297}, 
where he outlines a strategy based on general ideas of Connes. 
A major part of this paper deals with a rigorous implementation of this 
strategy, which takes a considerable effort. See 
Subsection~\ref{subsec:on the approach} for an overview.  

\begin{mthm}\label{thm:folvol bound by packing inequalities}
Let $C>0$ and $n\in\bbN$. Then there is a constant $\const_{C,n}>0$ with
the following property: If $M$ is an $n$-dimensional, 
closed, aspherical Riemannian manifold such that its universal covering 
$\widetilde{M}$ with the induced metric has the property (packing
inequality) that 
each ball of radius $1$ 
contains at most $Cr^{-n}$ disjoint balls of radius $r$ for every
$0<r<1$, then 
\begin{equation*}
\betti_i(\widetilde{M})\le\const_{C,n}\vol(M)\text{ for all $i\ge 0$.}
\end{equation*}
\end{mthm}

Let us briefly recall the well known 
relation between a \textit{lower Ricci curvature
bound} and a \textit{packing inequality} like in the hypothesis of the
preceding theorem. If $M$ satisfies the lower Ricci curvature bound 
$\ricci(M)\ge -(n-1)$ then also 
$\ricci(\widetilde{M})\ge -(n-1)$. Write $B(m,r)\subset\widetilde{M}$ 
for the ball of radius $r$ around $m\in\widetilde{M}$, and $B_{hyp}(r)$ 
for a ball of radius $r$ in hyperbolic $n$-space. 
According to the Bishop-Gromov inequality~\cite{gallot}*{Theorem~4.19}, 
if $R>r$, then 
\begin{equation*}
\frac{\vol(B(m,R))}{\vol(B(m,r))}\le\frac{\vol(B_{hyp}(R))}{\vol(B_{hyp}(r))}. 
\end{equation*}
Since the right hand side can bounded by $\const_n(R/r)^n$ for 
$r,R\le 1$, a mere volume estimate implies that at most
$\const_nr^{-n}$ disjoint balls of radius $r$ fit into a ball of
radius $1$ for a constant $\const_n$ only depending on the dimension
$n$. 
Thus the following corollaries  
are direct consequences of Theorem~\ref{thm:folvol
  bound by packing inequalities}. 

\begin{corollary_o}
For each $n\in\bbN$ there is a constant $\const_n>0$ with the
following property: If $M$ is an $n$-dimensional, 
closed, aspherical Riemannian manifold with lower Ricci curvature bound 
$\ricci(M)\ge -(n-1)$, then 
\begin{equation*}
\betti_i(\widetilde{M})\le\const_n\vol(M)\text{ for all
  $i\ge 0$.} 
\end{equation*}
\end{corollary_o}

\begin{corollary_o}[Main inequality for $L^2$-Betti numbers] 
For each $n\in\bbN$ there is a constant $\const_n>0$ with the
following property: If $M$ is an $n$-dimensional, closed,
aspherical manifold $M$, then   
\begin{equation*}
\betti_i(\widetilde{M})\le\const_n\minvol(M)\text{ for all
  $i\ge 0$.}
\end{equation*}
\end{corollary_o}

These two corollaries are the analog of Gromov's \textit{main
  inequality}~\cite{gromov(1982)}*{0.5} with \textit{simplicial
  volume} replaced by \textit{$L^2$-Betti numbers} and the additional
hypothesis of asphericity. However, Gromov's main inequality 
also applies to non-compact complete manifolds, and the constant is 
explicitly known. 

Notice that for aspherical $M$ we have 
$\betti_i(\widetilde{M})=\betti_i(\pi_1(M))$, and 
$\betti_i(\pi_1(M))$ is an orbit equivalence invariant of 
the fundamental group 
$\pi_1(M)$ by Gaboriau's work~\cite{gaboriau(2002b)}. It is a 
particularly interesting aspect of the previous corollary that 
it provides a non-trivial orbit equivalence invariant that bounds 
the minimal volume from below. 

According to the Hopf-Singer 
conjecture~\cite{lueck(2002)}*{Conjecture~11.1}),  
$\betti_i(\widetilde{M})=0$ if $2i\ne\dim(M)$ for every closed, aspherical 
manifold $M$. Examples of closed aspherical even-dimensional manifolds 
where the $L^2$-Betti number in the middle dimension is positive and 
the Hopf-Singer conjecture holds true include K\"ahler hyperbolic
manifolds~\cite{gromov(1991b)} 
and closed locally symmetric spaces of fundamental
rank zero~\cite{borel}. 

We prove the following new vanishing Theorem 
in Section~\ref{sec:amenable covers}. 
Recall that a subset $U\subset M$ of a topological space is 
\textit{amenable} if for every $x\in U$ 
the image $\im(\pi_1(U;x)\rightarrow\pi_1(M;x))$ is amenable. 

\begin{mthm}\label{thm:vanishing theorem}
Let $M$ be an $n$-dimensional, closed, triangulated (\textit{e.g.}~smooth),  
aspherical manifold. Assume 
that $M$ is covered by open, amenable sets such that every point
belongs to at most $n$ sets. Then 
\begin{equation*}
\betti_i(\widetilde{M})=0\text{ for all $i\ge 0$.}
\end{equation*}
\end{mthm}

Note that the original definitions of $L^2$-Betti numbers of Atiyah and 
Dodziuk do not apply to topological manifolds without 
triangulations, but L\"uck's theory does. Theorem~\ref{thm:vanishing
  theorem} actually stays true without the 
assumption~\textit{triangulated} but we omit the necessary
modifications in this case due to the length of the paper. 

We remark that a naive spectral sequence approach to the theorem 
above fails since the amenable subsets are in general not aspherical,
thus may have non-vanishing $L^2$-cohomology. 

The following result for simplicial complexes 
is also obtained in Section~\ref{sec:amenable
  covers}. 

\begin{mthm}\label{vanishing theorem for simplicial complexes}
Let $M$ be a finite aspherical simplicial complex. 
Assume that $M$ is covered 
by open amenable sets such that every point belongs to at most $n$
subsets. Then 
\begin{equation*}
\betti_i(\widetilde{M})=0\text{ for all $i\ge n$.}
\end{equation*}
\end{mthm}
Theorems~\ref{thm:vanishing theorem} and~\ref{vanishing theorem for
  simplicial complexes} together are 
the analog of Gromov's 
\textit{vanishing theorem}~\cite{gromov(1982)}*{Section~3.1} for
$L^2$-Betti numbers. 
Furthermore, Theorem~\ref{vanishing theorem for simplicial complexes}
generalizes a classical theorem of Cheeger and 
Gromov~\cite{cheeger+gromov}, which says that 
the $L^2$-Betti numbers of aspherical simplicial complexes 
with amenable fundamental group vanish. 

Gromov shows~\cite{gromov(1982)}*{Section~3.4} 
that for every dimension $n\in\bbN$ there is a constant 
$\epsilon_n>0$ with the property: 
If $M$ is an $n$-dimensional, closed Riemannian 
manifold such that $\ricci(M)\ge -(n-1)$ and $\vol(B(1))\le\epsilon_n$
for each unit ball $B(1)\subset M$ then there is an amenable cover as
in the hypothesis of Theorem~\ref{thm:vanishing theorem}; his proof
can be simplified a bit using the more recent Margulis lemma for Ricci
curvature by Cheeger and 
Colding~\citelist{\cite{cheeger+colding}*{Theorem
    8.7}\cite{gromov(1982)}*{Section~3.4}}. 
Hence we obtain the following interesting gap theorems. 

\begin{corollary_o}
For every $n\in\bbN$ there is a constant $\epsilon_n>0$ with the following
property: If $M$ is an $n$-dimensional, closed, aspherical Riemannian manifold
such that $\ricci(M)\ge -(n-1)$ and $\vol(B(1))\le\epsilon_n$ for each unit
ball $B(1)\subset M$ then 
\begin{equation*}
\betti_i(\widetilde{M})=0\text{ for all $i\ge 0$.}
\end{equation*}
\end{corollary_o}

\begin{corollary_o}[Isolation theorem for $L^2$-Betti numbers]
For every $n\in\bbN$ there is a constant $\epsilon_n>0$ with the
following property: 
If $M$ is an $n$-dimensional, closed, aspherical, smooth
manifold such that $\minvol(M)\le\epsilon_n$, then
\begin{equation*}
\betti_i(\widetilde{M})=0\text{ for $i\ge 0$.}
\end{equation*} 
\end{corollary_o}

This is the analog of Gromov's 
\textit{isolation theorem}~\cite{gromov(1982)}*{Section~0.5} for
$L^2$-Betti numbers.  

The Atiyah conjecture for a 
group~$\Gamma$~\cite{lueck(2002)}*{Chapter~10} 
predicts that the $L^2$-Betti numbers of any closed aspherical
manifold with fundamental group~$\Gamma$ 
are integers. Provided the Atiyah conjecture holds true,  
the previous corollaries would follow from 
Theorem~\ref{thm:folvol bound by packing inequalities} and its
corollaries. So far the Atiyah conjecture has been verified for
certain inductively defined classes of groups but not for fundamental groups 
of manifolds satisfying the hypothesis of the corollaries above nor
for any geometrically defined class of groups. 

As mentioned above, all stated results are analogs of theorems about 
simplicial volume or bounded cohomology in Gromov's 
paper~\cite{gromov(1982)}. Unlike there, 
the hypothesis of asphericity is essential in the present context. For example, 
Theorems~\ref{thm:folvol bound by packing inequalities}
and~\ref{thm:vanishing theorem}  
both fail for the $3$-sphere, which has vanishing minimal volume. 

Theorems~\ref{thm:folvol bound by packing inequalities}
and~\ref{thm:vanishing theorem} 
would follow from Gromov's 
results~\cite{gromov(1982)} if the following conjecture, 
which Gromov formulated as a question, would be 
true~\citelist{\cite{gromov(1993)}*{Section 8A
    on~p.~232}\cite{gromov(1999)}*{Remark e)
    on~p.~304}\cite{lueck(2002)}*{Chapter~14}}. 

\begin{conjecture_o}
For every dimension $n\in\bbN$ there is a constant $\const_n>0$ such
that for every  
$n$-dimensional, closed, aspherical, orientable manifold $M$ we have 
\begin{equation*}
\betti_i(\widetilde{M})\le\const_n\norm{M}\text{ for all $i\ge 0$.}
\end{equation*}
Here $\norm{M}$ denotes the simplicial volume of $M$. 
\end{conjecture_o}

\subsection{Conventions}\label{subsec:conventions}

The following framework is used throughout the paper. 
\begin{setup}\label{setup:orbit equivalence relation}
Let $(X,\mu)$ be a standard Borel space equipped with an atom-free 
probability Borel measure $\mu$. Let $\Gamma$ be a countable group acting 
(essentially) freely and $\mu$-preservingly on $X$. The \textit{orbit
  equivalence relation} $\scrR$ is the equivalence relation on $X$
given by 
\begin{equation*}
\scrR=\bigl\{(\gamma x,x)\in X\times X;~\gamma\in\Gamma,x\in X\bigr\}.
\end{equation*}
\end{setup}

Here the 
$\Gamma$-action is \textit{essentially free} if $X_0=\{x\in
X;~\Gamma_x\ne\{1\}\}$ is a $\mu$-null set. Upon replacing $X$ with 
$X\backslash X_0$, we can always achieve that the action is strictly
free. Therefore  
the word \textit{essentially} will be frequently omitted. 
Recall that $(X,\mu)$ as a measure space 
is isomorphic to $([0,1],\lambda)$ where $\lambda$ denotes the
Lebesgue measure. 

An orbit equivalence relation is an example of a 
\textit{countable measured equivalence 
relation}. For more information about that notion we refer 
to~\citelist{\cite{feldman+moore(a)}\cite{gaboriau(2002b)}*{Section~$0$}
\cite{furman(1999b)}*{Section~$2$}}. Every countable measured
equivalence relation arises as an orbit equivalence relation of 
a (not necessarily free) action of a countable group. 

For later reference, we record the following assumptions. 

\begin{setup}\label{setup:framework for simplicial complexes}
Let $M$ be a connected, finite simplicial complex with fundamental
group $\Gamma$. Let $(X,\mu)$ and $\scrR$ be as in
Assumption~\ref{setup:orbit equivalence relation}. All metric notions
about $M$ refer to the unique length metric that restricts to the 
standard Euclidean metric on simplices. 
\end{setup}

\begin{setup}\label{setup:framework for manifolds}
Let $M$ be an $n$-dimensional, connected, closed, oriented,
triangulated  
manifold with fundamental group $\Gamma$. 
Let $(X,\mu)$ and $\scrR$ be as in
Assumption~\ref{setup:orbit equivalence relation}. 
\end{setup}

Without loss of generality, we assume for the proofs of
Theorems~\ref{thm:folvol bound by packing inequalities}
and~\ref{thm:vanishing theorem} that the manifold $M$ is connected and 
orientable, that is, satisfies Assumption~\ref{setup:framework
  for manifolds}. If $M$ is non-orientable then there is a two-fold
orientation cover $p:\bar{M}\rightarrow M$. If $M$ satisfies the
hypothesis of Theorem~\ref{thm:folvol bound by packing inequalities}, 
then $\bar{M}$ with the induced Riemannian metric does so. If $M$
has an amenable cover $\{U_i\}_{i\in I}$ as in the 
hypothesis of Theorem~\ref{thm:vanishing theorem}, then
$\{p^{-1}(\bar{M})\}_{i\in I}$ is one for $\bar{M}$. Because of 
multiplicativity~\cite{lueck(2002)}*{Theorem~1.35 (9)} 
\begin{equation*}
\betti_i(\widetilde{\bar{M}})=2\betti_i(\widetilde{M})\text{ for all
  $i\ge 0$}
\end{equation*}
it is sufficient to prove the theorem for $\bar{M}$. 

A countable family $(X_i)_{i\in I}$ of Borel subsets 
of a measure space $(X,\mu)$ is called a \textit{(countable) Borel
  partition} if $\bigcup_{i\in I}X_i$ is a $\mu$-conull set and 
$\mu(X_i\cap X_j)=0$ for $i\ne j$. By abuse of notation, we just 
write $X=\bigcup_{i\in I}X_i$. 
The abbreviation 
\textit{a.e.} means either \textit{almost every} or \textit{almost
  everywhere}. 

\subsection{On the approach}\label{subsec:on the approach}

We present an elaborate version of 
Gromov's strategy~\cite{gromov(1999)}*{Section~5.33} 
to attack Theorem~\ref{thm:folvol bound by packing inequalities}, 
which itself is motivated by general ideas of Connes. 
The appropriate framework involves techniques from 
\textit{measured equivalence
  relations}, Connes' and Gaboriau's theory of 
\textit{$\scrR$-simplicial complexes} and other themes 
of \textit{measurable group theory}~\cite{shalom(2005)}. 

The way we set up the general framework is flexible enough 
to run the proofs of  
Theorems~\ref{thm:vanishing theorem} 
and~\ref{vanishing theorem for simplicial complexes}. However, 
a crucial difference is the geometric construction of 
$\scrR$-covers which involves the Ornstein-Weiss Lemma. 

\subsubsection{General Remark on bounding ($L^2$)-Betti
  numbers}\label{subsubsec:naive approach}
One method to bound the $i$-th Betti number and the 
$i$-th $L^2$-Betti number 
of a topological space~$M$ from above by a constant~$C$ 
is to realize $M$ as a homotopy retract 
in a simplicial complex~$S$ (that is, there are maps 
$f:M\rightarrow S$ and $g:S\rightarrow M$ with $g\circ f\simeq\id_M$)  
such that the number of $i$-simplices of~$S$ is at most $C$. 

In the sequel let  
$M$ be as in Assumption~\ref{setup:framework for manifolds}, and
assume that $M$ is 
aspherical, that is, $\widetilde{M}$ is model of $E\Gamma$, 
Then 
one could alternatively try to find a free $\Gamma$-simplicial
complex~$S$ 
with at most $C$~equivariant $i$-simplices and an 
equivariant map~$f:\widetilde{M}\rightarrow S$. The advantage 
of asphericity and 
working equivariantly is that one automatically gets an equivariant 
map~$g:S\rightarrow\widetilde{M}$ and an equivariant 
homotopy~$g\circ f\simeq\id_{\widetilde{M}}$ from the universal property of
$E\Gamma$, which then leads to the same estimate. 
\subsubsection{The category of $\scrR$-spaces}
We use a similar method (which does not work anymore to bound 
Betti numbers) in the category of $\scrR$-spaces
(Section~\ref{sec:Topological Properties of R-spaces}) 
instead of $\Gamma$-spaces, where $\scrR$ is as in
Assumption~\ref{setup:orbit equivalence relation}. 
An $\scrR$-space is the realization 
of a fiberwise locally finite 
$\scrR$-simplicial complex (Definition~\ref{def:R-simplicial
  complex}) in the 
sense of~\cite{gaboriau(2002b)}*{Section~2}. One example is 
$X\times\widetilde{M}$ with the $\scrR$-action $(\gamma x,x).(x,m)=(\gamma
x,\gamma m)$. The morphism in this category, called \textit{geometric
  $\scrR$-maps} (Definition~\ref{def: topological R-maps}), are
fiberwise continuous and proper. The category of $\scrR$-spaces is an
extension of the combinatorial framework~\cite{gaboriau(2002b)} to 
a topological one. 

\subsubsection{Homotopy retracts of $\scrR$-spaces by $\scrR$-covers}
An $\scrR$-cover $\scrU$ of $X\times\widetilde{M}$ is an equivariant family of 
sets of the form $A\times U$ with $A\subset X$ Borel and
$U\subset\widetilde{M}$ open (Definition~\ref{def:measurable covering
  and packing}). Such an $\scrR$-cover gives rise to a nerve 
construction $\scrN(\scrU)$, which is naturally an $\scrR$-simplicial
complex, and a geometric $\scrR$-map 
\begin{equation}\label{eq:to nerve}
\phi: X\times\widetilde{M}\rightarrow\scrN(\scrU). 
\end{equation}
Also in this context asphericity yields a homotopy retract, that
is, a geometric $\scrR$-map 
$\psi:\scrN(\scrU)\rightarrow X\times\widetilde{M}$ with
$\psi\circ\phi\simeq\id_{X\times\widetilde{M}}$. 

\subsubsection{Bounding the $L^2$-Betti number in the top dimension
  $n$ under the hypothesis of Theorem~\ref{thm:folvol bound by packing
    inequalities}}
\label{subsubsec:bounds in the top dimension}
Now let $M$ be as in Theorem~\ref{thm:folvol bound by packing
  inequalities}. It is a standard trick that one obtains from 
the packing type hypothesis on $\widetilde{M}$ 
a cover whose multiplicity is bounded in terms
of the constant~$C$~\cite{berger}*{Lemma 125 on p.~333}. 
But in general there is no
way to make such a cover $\Gamma$-equivariant without loosing control
over its multiplicity. If there would be, we could 
prove Theorem~\ref{thm:folvol bound by packing inequalities} by
proceeding as in~\ref{subsubsec:naive approach}. 

However, it is possible to construct an $\scrR$-cover $\scrU$ such that 
the induced cover on $\{x\}\times\widetilde{M}$ 
has the correct multiplicity for a.e. $x\in X$ 
(Theorem~\ref{thm:suitable-coverings}). 
Then one can modify $\phi$ from~(\ref{eq:to nerve}) by a geometric
$\scrR$-homotopy such that 
the \textit{equivariant number of $n$-simplices} of
$\im(\phi)\subset\scrN(\scrU)$ 
in the sense of Definition~\ref{def:weighted number of simplices} 
(equivalently: the measure~$\nu^n(\scrR\backslash\Sigma^{(n)})$
defined in~\cite{gaboriau(2002b)}*{Section 2.2.3}) is dominated 
by $\const_{C,n}\vol(M)$ where $\const_{C,n}$ is a constant only
depending on the dimension~$n$ and the constant~$C$ from
Theorem~\ref{thm:folvol bound by packing inequalities}.

Now an application of Gaboriau's theory~\cite{gaboriau(2002b)} would 
yield the stated bound 
\begin{equation}\label{eq:bound in top dimension}
\betti_n(\widetilde{M})\le\const_{C,n}\vol(M)
\end{equation}
in the top dimension, which alone is useless if one assumes 
the Hopf-Singer conjecture for $M$ (mentioned in
Subsection~\ref{subsec:statement of results}). 
We introduce another 
tool to obtain the bound in all degrees: \textit{Foliated
  singular homology and the support mass.}
\subsubsection{Foliated singular homology}\label{subsubsec:foliated
  singular homology}
We define a homology theory~$\scrH_\ast$ 
(Section~\ref{sec: homology theories of R-spaces}) 
for the category of 
$\scrR$-spaces, called \textit{foliated singular homology} that 
is a singular version of the sheaf-theoretic 
\textit{tangential homology} of measured foliations~\cite{moore+schochet}.  

Actually, $\scrH_\ast(\Sigma)$ for an
$\scrR$-space~$\Sigma$ is really defined 
in terms of the laminated quotient space~$\scrR\backslash\Sigma$, but 
for conceptual and technical reasons it is better to work with the
$\scrR$-space instead of its quotient. For example, the universal
property that gives us the map $\psi$ is not transparent on the
quotient level. 

The definition of $\scrH_n(\Sigma)$ is modelled on Gromov's
description of Connes's \textit{foliated simplicial
volume}~\cite{gromov(1991)}*{2.4.B}. Indeed, the foliated simplicial
volume for the measurable foliation $\Gamma\backslash (X\times\widetilde{M})$ 
is most naturally defined in terms of cycles representing 
the fundamental class in 
$\scrH_n(X\times\widetilde{M})$ (see Remark~\ref{rem:fundamental class in
    tensored complex} for the notion of \textit{fundamental
    class}). For our purposes, we deal with another, yet related
  numerical invariant on
$\scrH_n(\Sigma)$, called the \textit{support mass} and defined 
in Subsection~\ref{subsec:norms}. 

It turns out that in the situation 
of~\ref{subsubsec:bounds in the top dimension} we 
not only get the bound~(\ref{eq:bound in top dimension}), 
but, using the properties of $\scrH_\ast$ developed in 
Section~\ref{sec: homology theories of R-spaces}, 
we can also bound the 
support mass of the fundamental class in 
$\scrH_n(X\times\widetilde{M})$ by $\const_{C,n}\vol(M)$
(Theorem~\ref{thm:homotopy further}). This implies 
(Theorem~\ref{thm:inclusion and support norm}) that 
the image of the fundamental class in
$\Hn_n(M;\bbZ)$ under the homomorphism 
induced by inclusion of coefficients 
\begin{equation*}
\Hn_n(M;\bbZ)=\Hn_n\bigl(\bbZ\otimes_{\bbZ\Gamma}\Cn_\ast(\widetilde{M})\bigr)\rightarrow 
\Hn_n\bigl(\lz\otimes_{\bbZ\Gamma}\Cn_\ast(\widetilde{M})\bigr)=
\Hn_n^\Gamma\bigl(\widetilde{M};\lz\bigr) 
\end{equation*} 
can be represented by a cycle $\sum_{i=1}^mf_i\otimes\sigma_i$ with 
$f_i\in\lz$ and $\sigma_i\in\sing_n(\widetilde{M})$ such that 
\begin{equation*}
\sum_{i=1}^m\mu\bigl(\supp(f_i)\bigr)<\const_{C,n}\vol(M).
\end{equation*}
Such a 
representing cycle of the fundamental class leads to an estimate for 
$\betti_i(\widetilde{M})$ for all $i\ge 0$ (in our case:
$\betti_i(\widetilde{M})\le\const_{C,n}\vol(M)$) by a Poincar\'{e} duality 
argument. See Theorem~\ref{thm: l2 betti and mass} in the Appendix. 

\subsubsection{The $\scrR$-covers used for Theorems~\ref{thm:vanishing
    theorem} and~\ref{vanishing theorem for simplicial complexes}}
The proofs of Theorems~\ref{thm:vanishing theorem} and~\ref{vanishing
  theorem for simplicial complexes} differ 
from that of Theorem~\ref{thm:folvol bound by packing inequalities} 
in the construction of 
the $\scrR$-covers $\scrU$. 
Whereas in the case of 
Theorem~\ref{thm:folvol bound by packing inequalities} 
the input from ergodic theory is very modest, we now
employ  
the \textit{generalized Rokhlin lemma} of Ornstein and 
Weiss~\cite{ornstein+weiss(1987)} to construct suitable $\scrR$-covers 
$\scrU_\delta$ for any $\delta>0$ 
such that the number of equivariant $k$-simplices of 
$\scrN(\scrU)$ is less than $\const\cdot\delta$ for $k\ge n$. 
The characteristic feature of elements $A\times U$ of $\scrU_\delta$,
where $A\subset X$ and $U\subset\widetilde{M}$, 
is that $\mu(A)$ gets small and $U$ gets large when~$\delta$ tends to
zero. 

The proof of Theorem~\ref{vanishing theorem for simplicial complexes} 
is easier than that of Theorem~\ref{thm:vanishing theorem} since 
we do not need the technique in~\ref{subsubsec:foliated
  singular homology}. Instead we apply Gaboriau's 
techniques~\cite{gaboriau(2002b)}.

\section{$\scrR$-Spaces and $\scrR$-simplicial
  Complexes}\label{sec:Topological Properties of R-spaces}    

In Subsections~\ref{subsec:R-spaces and
  R-simplicial complexes} to~\ref{subsec:Simplicial approximation},   
$\scrR$ denotes an arbitrary countable measured equivalence relation
on a standard Borel space~$(X,\mu)$ equipped with a probability Borel 
measure~$\mu$. In Subsection~\ref{subsec:Measurable coverings and
  packings} we refer to the more special situation of 
Assumption~\ref{setup:orbit equivalence relation}. 

\subsection{$\scrR$-simplicial complexes}\label{subsec:R-spaces and
  R-simplicial complexes}   
To fix the terminology and for the convenience of the reader, we
recall in 
Subsection~\ref{subsec:R-spaces and R-simplicial complexes} 
basic notions about $\scrR$-simplicial complexes as presented 
in~\cite{gaboriau(2002b)}. 
\begin{definition}\label{def: X-space}
An \emph{$(X,\mu)$-space} or \emph{$X$-space} 
is a standard Borel space~$\Sigma$ together with a Borel 
map~$p_\Sigma:\Sigma\rightarrow X$ such that the 
fibers~$\Sigma_x=p_\Sigma^{-1}(x)$, $x\in X$, are countable. 
A \emph{map of $X$-spaces} $f:\Sigma\rightarrow\Phi$ 
is a Borel map such that $\pr_\Phi\circ f=\pr_\Sigma$. 
\end{definition}

We denote the \textit{fiber product} 
of 
$X$-spaces~$\Sigma\times_X\Sigma'=\{(u,v)\in\Sigma\times\Sigma';~p_{\Sigma}(u)=
p_{\Sigma'}(v)\}$ by $\Sigma\ast\Sigma'$. 
\begin{definition}\label{def:R action}
An \textit{$\scrR$-action} on an $X$-space $\Sigma\rightarrow X$ is 
map of $X$-spaces
$\scrR\ast\Sigma\rightarrow\Sigma$, $((y,x),u)\mapsto (y,x).u$, 
such that 
\begin{enumerate}[a)]
\item $(y,x).\Sigma_x\subset\Sigma_y$,
\item $(x,x).u=u$ for every $x\in X$ and every $u\in\Sigma_x$, and 
\item $(z,y).\bigl((y,x).u\bigr)=(z,x).u$ for all $x,y,z\in X$ and 
every $u\in\Sigma_x$. 
\end{enumerate}
\end{definition}

\begin{remark}\label{rem: R action from Gamma action}
If $\scrR$ is as in Assumption~\ref{setup:orbit equivalence relation},
then an \textbf{$\scrR$-action on $\Sigma$ is nothing else but a measurable
$\Gamma$-action on $\Sigma$ such that
$\gamma\cdot\Sigma_x\subset\Sigma_{\gamma x}$.} 
\end{remark}

Fibered products of $X$-spaces with $\scrR$-actions
carry by definition the diagonal $\scrR$-action. 
An $\scrR$-action that has a \textit{Borel fundamental domain $\scrF$},
that is a Borel subset of $\Sigma$ whose intersection with every
$\scrR$-orbit consists of 
exactly one element, is called \textit{discrete}. 
\begin{definition}\label{def:measure on X spaces}
The \textit{natural measure} on the $(X,\mu)$-space $\Sigma$ is
defined as 
\begin{equation*}
\nu^\mu(U)=\int_X\#\left(p_\Sigma^{-1}(x)\cap
  U\right)d\mu(x)\text{ for all $U\subset\Sigma$ Borel}.
\end{equation*}
If the choice of $\mu$ is clear from the context, the superscript~$\mu$ 
in $\nu^\mu(U)$ is omitted. 
\end{definition}
\begin{definition}\label{def:transversal measure}
Let $\Sigma$ be an $(X,\mu)$-space with a discrete $\scrR$-action. 
The measure~$\nu_t^\mu$ on $\scrR\backslash\Sigma$ obtained from the 
restriction~$\nu^\mu\vert_\scrF$ of the natural measure 
to a Borel fundamental domain~$\scrF\subset\Sigma$ after 
identifying 
$\scrR\backslash\Sigma$ with $\scrF$ is called the \textit{transversal 
measure} on $\scrR\backslash\Sigma$. As before, the superscript~$\mu$
is omitted in the notation if the choice of~$\mu$ is clear from the
context. 
\end{definition}

\begin{remark_o}
It is easily verified that the preceding definition does not depend 
on the choice of the Borel fundamental domain. 
\end{remark_o}

\begin{definition}\label{def:R-simplicial complex}
An \textit{$\scrR$-simplicial complex} $\Sigma$ consists of the following
data: 
\begin{enumerate}[--]
\item an $X$-space $\Sigma^{(0)}\rightarrow X$ with a discrete $\scrR$-action. 
\item for each $n\in\bbN$ a Borel subset $\Sigma^{(n)}\subset
 \underbrace{\Sigma^{(0)}\ast\Sigma^{(0)}\ast\ldots\ast\Sigma^{(0)}}_{\text{
    $n+1$ times}}$\\ (ordered $n$-simplices) 
\end{enumerate}
subject to four conditions: 
\begin{enumerate}[a)]
\item (permutations) $\Sigma^{(n)}$ is invariant under permutations of
  the coordinates. 
\item (non-degeneracy) $(v_0,\ldots, v_n)\in\Sigma^{(n)}$ implies
  $v_0\ne v_n$. 
\item (boundary condition) $(v_0,\ldots, v_n)\in\Sigma^{(n)}$ implies
  $(v_1,\ldots, v_n)\in\Sigma^{(n-1)}$. 
\item (invariance) $\scrR.\Sigma^{(n)}=\Sigma^{(n)}$. 
\end{enumerate}
\end{definition}

\begin{definition}\label{def:R simplicial map}
An \textit{$\scrR$-simplicial map} $\phi:\Sigma\rightarrow\Phi$ 
between $\scrR$-simplicial
complexes is a Borel map over $X$ such that 
$\phi_x:\Sigma_x\rightarrow\Phi_x$ is simplicial for a.e. $x\in X$ and 
$\phi$ is \textit{$\scrR$-equivariant} in the sense that 
for a.e. $x\in X$, every $y\in X$ with $(y,x)\in\scrR$, 
and every $m\in\Sigma_x$ we have $\phi((y,x).m)=(y,x).\phi(m)$. 
\end{definition}

The reader may notice that while the notion of an $\scrR$-action
(Definition~\ref{def:R action}) is 
defined in a strict sense the conditions of the previous definition 
are only required up to null-sets. 

\begin{definition}\label{def:subdivision}
The \textit{barycentric subdivision} $\sd(\Sigma)$ of $\Sigma$ 
is the $\scrR$-simplicial complex obtained by taking 
the barycentric subdivision on each fiber $\Sigma_x$. 
The $n$-fold barycentric subdivision is denoted by $\sd^{(n)}(\Sigma)$. 
Abstractly, $\sd(\Sigma)$ is defined as follows. 
There is a natural action of the symmetric group $\symm(n+1)$ on 
$\Sigma^{(n)}$ which commutes with the $\scrR$-action. The quotient 
$\symm(n+1)\backslash\Sigma^{(n)}$ is still an $X$-space with a
discrete $\scrR$-action. 
The set of $0$-simplices is 
\[\sd(\Sigma)^{(0)}=\coprod_{n\ge 0}\symm(n+1)\backslash\Sigma^{(n)}.\]
Thus $\sd(\Sigma)^{(0)}$ is partially
ordered by inclusion. The set of $n$-simplices $\sd(\Sigma)^{(n)}$  
consists of $(n+1)$-tuples $(v_0,\ldots,v_n)$ with the property that 
the $v_i$ lie in some common fiber of $\sd(\Sigma)^{(0)}$ and 
$\{v_0,\ldots,v_n\}$ can be totally ordered. 
\end{definition}

\begin{definition}\label{def:weighted number of simplices}
The \textit{(weighted) number of (non-oriented) equivariant $n$-simplices} of an 
$\scrR$-simplicial complex $\Sigma$ is defined as 
$\nu_t\bigl((\symm(n+1)\times\scrR)\backslash\Sigma^{(n)}\bigr)$ and denoted by 
$\sigma_n(\Sigma)$. 
\end{definition}

\begin{example}\label{exa:easy R simplicial complex}
Retain the notation of Assumption~\ref{setup:orbit equivalence
  relation}. 
Let $S$ be a free $\Gamma$-equivariant simplicial complex 
Then $X\times S$ is an $\scrR$-simplicial complex via 
\begin{equation*}
(\gamma x,x)(x, s)=(\gamma x, \gamma s). 
\end{equation*}
The set of $n$-simplices is $\Sigma^{(n)}=X\times S^{(n)}$. 

Since an $\scrR$-action is defined in a strict sense, 
we may and will assume here that the $\Gamma$-action on~$X$ 
is strictly free (cf.~Subsection~\ref{subsec:conventions}). 

At any rate, in the proofs of
Theorems~\ref{thm:folvol bound by packing
  inequalities},~\ref{thm:vanishing theorem} and~\ref{vanishing
  theorem for simplicial complexes} we fix the group $\Gamma$ and thus
could also speak of $\Gamma$- instead of $\scrR$-actions
(cf.~Remark~\ref{rem: R action from Gamma action}), but the framework
of 
$\scrR$-actions is more conceptual, especially 
in the definition and the functorial properties of singular
foliated homology (Section~\ref{sec: homology theories of R-spaces}). 
\end{example}

\subsection{Standard
  embeddings}\label{subsec:standard embeddings}
Let $\Sigma$ be an $\scrR$-simplicial complex. Pick a countable set
$I$ and an isomorphism $\scrR\rightarrow X\times I$ of 
$X$-spaces~\cite{feldman+moore(a)}.  

\begin{definition}\label{def:delta complex}
Let $\bDelta(I)$ be the simplicial complex 
whose vertex set is $I$ and whose simplices consist of all finite subsets
of $I$. The realization of $\bDelta(I)$ is denoted by the same symbol;
to which meaning we refer should be clear from the context. 
\end{definition}

Let $\scrF$ be a Borel fundamental domain of $\Sigma^{(0)}$. Then there is
a countable Borel partition $\scrF=\coprod_{n\in\bbN} \scrF_n$ such
that $p_\Sigma$ 
is injective on each $\scrF_j$ (theorem of
selection;~see~\textit{e.g.}~\cite{sauer(2005)}*{Lemma~3.1} 
for the version needed here). The injective map 
\begin{equation*}
\Sigma^{(0)}=\scrR.\scrF\rightarrow \scrR\times\bbN\xrightarrow{\cong}
(X\times I)\times\bbN=X\times\Delta(I\times\bbN)^{(0)}, 
\end{equation*}
where the first arrow maps $(y,x).u$ with 
$u\in\scrF_n$ and $x=p_\Sigma(u)$ to
$(y,x,n)$ and the second arrow comes from the isomorphism above, 
extends to an embedding 
\begin{equation}\label{eq:standard embedding}
\Sigma\rightarrow
X\times\bDelta(I\times\bbN). 
\end{equation}
We say that (\ref{eq:standard embedding}) is a 
\textit{standard embedding of $\Sigma$}. 
Note that under Assumption~\ref{setup:orbit equivalence relation} 
we could take $I=\Gamma$. 

\subsection{$\scrR$-spaces and geometric
  $\scrR$-maps}\label{subsec:maps}

The \textit{geometric realization} $\abs{\Sigma}$ of an 
$\scrR$-simplicial complex $\Sigma$ is the disjoint union 
of the geometric realizations $\abs{\Sigma_x}$ of the fibers. 

\begin{definition}\label{def: R-space}
By definition, an \textit{$\scrR$-space} is the geometric realization of some 
$\scrR$-simplicial complex $\Sigma$ such that $\Sigma_x$ is locally
finite for a.e. $x\in X$. 
Here $\Sigma$ is part of the data of the $\scrR$-space. 

A standard embedding $\phi$ as in~(\ref{eq:standard embedding}) 
yields an injection 
$\abs{\phi}:\abs{\Sigma}\rightarrow X\times\abs{\bDelta(I\times\bbN)}$. 
The standard Borel structure of $\abs{\Sigma}$ is defined as the restriction 
of the product Borel structure on
$X\times\abs{\bDelta(I\times\bbN)}$. 
\end{definition}

\begin{lemma}
The Borel structure in the previous definition does not depend on the
choice of a standard embedding, and $\im(\abs{\phi})\subset
X\times\abs{\bDelta(I\times\bbN)}$ is a Borel subset. 
\end{lemma}

We skip the proof of this easy lemma. Its first part follows 
from Lemma~\ref{lem:compacta and countable variance} below (whose 
proof is also left to the reader) and the fact that, 
if $I$ is countable, the subsets $A\times K$, where $A\subset X$ is Borel and
$K\subset\abs{\bDelta(I)}$ compact, form a basis of the Borel algebra of 
$X\times\abs{\bDelta(I)}$.

\begin{lemma}\label{lem:compacta and countable variance}\hfill
\begin{enumerate}[a)]
\item Let $\phi:X\times I\rightarrow X\times J$ be an isomorphism of
  $X$-spaces. The induced bijection
  $\abs{\phi}:X\times\abs{\bDelta(I)}\rightarrow
  X\times\abs{\bDelta(J)}$ has the following property: 
For every compact subset $K\subset
\abs{\bDelta(I)}$ there is a countable Borel partition of 
$X$ into sets $X_n$ such that $\abs{\phi}\vert_{ X_n\times K}$ is a product
$\id_{X_n}\times f_n$. In particular, for every Borel subset
$A\subset X$, $\abs{\phi}(A\times K)$ is a union of sets  
$A_n\times K_n$ with $K_n$ compact and $A_n=X_n\cap A$.  
\item For an $\scrR$-simplicial complex $\Sigma$ 
let $\phi$ and $\theta$ be standard embeddings 
$\Sigma\rightarrow X\times\bDelta(I\times\bbN)$. 
Let $K\subset\abs{\bDelta(I\times\bbN)}$ be compact and 
$A\subset X$ be a Borel subset such that $A\times K\subset\im(\abs{\phi})$. 
Then there is a 
countable Borel partition of $A$ by sets $A_n$ such that 
each restriction 
$\abs{\theta}\circ\abs{\phi}^{-1}\vert_{ A_n\times K}$ is a product 
map $\id_{A_n}\times f_n$ with continuous 
$f_n:K\rightarrow\abs{\bDelta(I\times\bbN)}$. 
\end{enumerate}
\end{lemma}

\begin{definition}\label{def:geometric simplices}
The $n$-simplices of 
an $\scrR$-simplicial complex correspond to subspaces of its realization 
homeomorphic to $\Delta^n$, which we call \textit{geometric $n$-simplices}. 
\end{definition}

\begin{remark}\label{rem:simplicial path on R spaces by convention}
We equip the fibers of an $\scrR$-space with the unique \textit{length metric}
that restricts to the standard Euclidean metric on geometric
simplices. All metric notions refer to this metric. 
Since a.e. fiber is locally compact, the weak 
topology coincides with the topology obtained from the metric on
a.e. fiber. 
\end{remark}

\begin{remark}\label{rem:on terminology}
Note that we use Greek capital letters 
($\Sigma,\Phi,\ldots$) to denote both $\scrR$-spaces 
and $\scrR$-simplicial complexes. The $\abs{\Sigma}$-notation is only 
used if we want to refer explicitly to the 
underlying $\scrR$-simplicial complex. 
\end{remark}

\begin{remark}\label{rem:induced action}
We leave it to the reader to verify that 
the $\scrR$-action on $\Sigma$ induces an $\scrR$-action on
$\abs{\Sigma}$ in the same sense as in Definition~\ref{def:R action} except 
\textit{map of $X$-spaces} should be replaced by \textit{Borel map
  over $X$}. This action is fiberwise continuous. The union of 
geometric simplices corresponding to fundamental domains 
of $\Sigma^{(n)}$ over all $n\ge 0$ is a Borel fundamental domain of
$\abs{\Sigma}$. 
\end{remark}

\begin{definition}\label{def:leaf}
Since $\Sigma$ has an $\scrR$-fundamental domain, the restriction of 
the projection $\pr:\Sigma_x\rightarrow\scrR\backslash\Sigma$ is injective
for all $x\in X$. 
The image $\pr(\Sigma_x)\subset\scrR\backslash\Sigma$, which can be
identified with $\Sigma_x$, is called \textit{the leaf at $x\in X$} and 
denoted by $\leaf_\Sigma(x)$. 
\end{definition}

The quotient $\scrR\backslash\Sigma$ can be seen as a foliated space, 
foliated by the leaves $\leaf_\Sigma(x)$. 

\begin{definition}\label{def:countable variance}
Fix a choice of standard embeddings $\Sigma,\Phi\rightarrow 
X\times\bDelta(I\times\bbN)$. 
A map $\phi:\abs{\Sigma}\rightarrow\abs{\Phi}$ over $X$ is said to be 
of \textit{countable variance}, if for any set 
$A\times K$, $A\subset X$ Borel, $K\subset\abs{\bDelta(I\times\bbN)}$
compact, there is a countable Borel partition $A=\bigcup_{n\in\bbN} A_n$ such
that $A_n\times K\subset\abs{\Sigma}$ and 
each restriction $\phi\vert_{ A_n\times K}$ is a product. 
\end{definition}

Independence of the chosen standard embeddings is implied by 
Lemma~\ref{lem:compacta and countable variance}. 
Further, it is clear that countable variance implies measurability. 

\begin{definition}\label{def: topological R-maps}
A \textit{geometric $\scrR$-map} between $\scrR$-spaces is
a map over $X$ of countable variance 
that is $\scrR$-equivariant (in the sense of 
Definition~\ref{def:R simplicial map}), continuous on a.e.~fiber 
and proper on
a.e.~fiber. We say that a geometric $\scrR$-map between $\scrR$-spaces is 
\textit{$\scrR$-simplicial} if it is the realization of an
$\scrR$-simplicial map of the underlying $\scrR$-simplicial
complexes. 
\end{definition}

\begin{remark}\label{rem:induced map on the leaves}
A geometric $\scrR$-map $\phi:\Sigma\rightarrow\Phi$ 
induces proper maps $\leaf_\Sigma(x)\rightarrow\leaf_\Phi(x)$ for
a.e. $x\in X$. 
\end{remark}

It would be more natural and less cumbersome in some places 
to replace 
the condition of \textit{countable variance} by
\textit{measurability}. However, for technical reasons we need countable
variance at some places, notably in the proofs of
Theorems~\ref{thm:functoriality} and~\ref{thm:inclusion and support
  norm}. 

The easy proof (similar to Lemma~\ref{lem:compacta and countable
  variance}) of the following lemma is left to the reader. 

\begin{lemma}\label{lem:realization of simplicial map is geometric}
Let $\Sigma, \Phi$ be $\scrR$-simplicial complexes with locally finite
fibers a.e and 
$\phi:\Sigma\rightarrow\Phi$ an $\scrR$-simplicial map. If 
$\phi_x$ is proper for a.e. $x\in X$, then the
realization~$\abs{\phi}$, 
defined fiberwise as
$\abs{\phi_x}:\abs{\Sigma_x}\rightarrow\abs{\Phi_x}$, is a geometric
$\scrR$-map. 
\end{lemma}

\subsection{Simplicial approximation of geometric
  $\scrR$-maps}\label{subsec:Simplicial approximation}  
In this Subsection we introduce simplicial approximation 
theorems in the context of $\scrR$-space. 
We start by recalling some terminology concerning simplicial approximation. 
The smallest simplex in (the realization) of a simplicial complex that
contains the point $m$ is denoted by $\carr(m)$. The \emph{open star}
of a vertex $v$ in a simplicial complex is denoted by $\openstar(v)$. 
Notice that $x\in\openstar(v)\Leftrightarrow v\in\carr(x)$. 
Let $f,g$ be two maps 
from a topological space $M$ to a simplicial complex $K$. Then
$g$ is called an \emph{approximation of $f$} if $g(x)\in\carr(f(x))$
for all $x\in M$. If $M$ and $g$ are simplicial, then $g$ is called a
\emph{simplicial approximation of $f$}. If $\phi,\psi$ are geometric
$\scrR$-maps between $\scrR$-spaces, we call $\phi$ an 
($\scrR$-simplicial) approximation of $\psi$ if $\phi_x$ is a
(simplicial) approximation of $\psi_x$ for a.e. $x\in X$. 
The \textit{Lebesgue number} of an open
cover of a metric space is the supremum of all $r\ge 0$ such that 
every set of diameter less than $r$ is contained in an element of the
cover.

\begin{lemma}\label{lem:homotopy for approximation}
Let $\phi:\abs{\Sigma}\rightarrow\abs{\Phi}$ be a geometric
$\scrR$-map. Let $\psi:\abs{\Sigma}\rightarrow\abs{\Phi}$ be an
$\scrR$-equivariant map that is of countable variance and continuous
in a.e. fiber. Suppose that $\psi_x$ is an approximation of $\phi_x$
for a.e. $x\in X$. Then $\psi_x$ is proper for a.e. $x\in X$, thus 
$\psi$ is a
geometric $\scrR$-map. Furthermore, there is a geometric $\scrR$-homotopy
between $\phi$ and $\psi$. 
\end{lemma}

\begin{proof}
Define $H:\abs{\Sigma}\times [0,1]\rightarrow\abs{\Phi}$ to be the map
such that $H\vert_{\{m\}\times [0,1]}$ parametrizes the straight line segment 
connecting $\psi(m)$ and $\phi(m)$ (within the simplex $\carr(\phi(m))$). 
Equivariance of $H$ is clear. 
First we verify that $H_x$ is proper (in particular, $\psi_x$ is). 
Let $K\subset\abs{\Phi_x}$ be a compact subcomplex and 
$\{k_1,\ldots,k_p\}$ be the finitely many vertices of $K$. 
We show that 
\begin{equation}\label{eq:preimage compact}
H_x^{-1}(K)\subset
\bigcup_{i=1}^p\phi_x^{-1}\bigl(\openstar(k_i)\bigr)\times [0,1]. 
\end{equation} 
Let $(m,t)\in H_x^{-1}(K)$. Then 
\begin{equation*}
H_x(m,t)\in\carr\bigl(\phi_x(m)\bigr)\cap K\ne\emptyset. 
\end{equation*}
There is 
a vertex $k_i$ in $K$ such that
$k_i\in \carr(\phi_x(m))$, thus 
$m\in\phi_x^{-1}(\openstar(k_i))$ 
showing~(\ref{eq:preimage compact}). Since $\Phi_x$ is locally finite 
and $\phi_x$ is proper,~(\ref{eq:preimage compact}) implies that 
$H_x^{-1}(K)$ is compact. 
Next we verify that $H$ is of countable
variance. 
Choose standard embeddings $\Sigma,\Phi\rightarrow
X\times\bDelta(I\times\bbN)$. Let $K\subset\bDelta(I\times\bbN)$ be 
a compact subcomplex. Let $\{A_{n}\}_{n\in\bbN}$ and $\{B_{n}\}_{n\in\bbN}$ be 
Borel partitions of~$X$ such that $\phi\vert_{A_n\times K}$ and
$\psi\vert_{B_n\times K}$ are product maps. If $\{C_n\}_{n\in\bbN}$ is a 
refinement of~$\{A_{n}\}_{n\in\bbN}$ and $\{B_{n}\}_{n\in\bbN}$, then 
$H\vert_{C_n\times (K\times [0,1])}$ is a product map for all~$n\in\bbN$. 
So $H$ is the geometric $\scrR$-homotopy between $\phi$ and $\psi$. 
\end{proof}

\begin{definition}\label{def:lebesgue number}
Let $\phi:\abs{\Sigma}\rightarrow\abs{\Phi}$ 
be a geometric $\scrR$-map. 
Let $L(x)$ be the Lebesgue number of the pullback under $\phi_x$ 
of the open star cover of $\abs{\Phi_x}$. The \textit{Lebesgue
  number of $\phi$} is defined as the essential infimum of 
$\{L(x);~x\in X\}$. 
\end{definition}

\begin{theorem}\label{thm: simplicial approximation}
Let $\phi:\abs{\Sigma}\rightarrow\abs{\Phi}$ be a geometric $\scrR$-map. 
If the Lebesgue number of~$\phi$ is positive, then there is an 
$n\in\bbN$ and an 
$\scrR$-simplicial approximation
$\psi:\sd^{(n)}(\Sigma)\rightarrow\Phi$ of $\phi$. 
Further, $\abs{\psi}$ is a geometric $\scrR$-map, 
and $\phi$ and $\abs{\psi}$ are geometrically $\scrR$-homotopic. 
\end{theorem}

\begin{proof}
We examine the classical proof to see that 
the simplicial approximations on the fibers assemble to an
$\scrR$-simplicial map. 

Let $\delta>0$ be the Lebesgue number of $\phi$. Take $n\in\bbN$ large
enough so that the diameter of a (geometric) simplex in
the $n$-fold barycentric subdivision of $\abs{\Sigma_x}$ with respect to the
length metric of $\abs{\Sigma_x}$ is less than $\delta/2$
for a.e. $x\in X$. Then the 
diameter of the open star of a vertex $v$ in 
$\bigl\lvert\sd^{(n)}(\Sigma)_x\bigr\rvert=\abs{\Sigma_x}$, is
less than $\delta$, 
thus its image under $\phi$ is contained in the open star of some vertex $w$
of $\abs{\Phi_x}$. For purposes of the 
proof we can now forget $n$ and assume that the image of an open 
star in $\abs{\Sigma_x}$ is contained in some open star of $\abs{\Phi_x}$. 

Pick standard embeddings $\Sigma,\Phi\rightarrow
X\times\bDelta(I\times\bbN)$. 
Let $\scrF$ be a Borel fundamental domain for
$\Sigma^{(0)}\subset X\times (I\times\bbN)$; $\scrF$ is the
disjoint union of sets $A_k\times\{(i_k,n_k)\}$, $k\in\bbN$, with
$A_k\subset X$ Borel, $i_k\in I$ and $n_k\in\bbN$. For a.e. $x\in X$
and every $k\in\bbN$ there are $j\in I$ and $m\in\bbN$ such that 
we have for the open stars
\begin{equation}\label{eq:open stars inclusion}
\phi_x\bigl(\openstar(x,i_k,n_k)\bigr)\subset\openstar(x,j,m).
\end{equation}
Since for fixed~$k\in\bbN$ and $x$ running through~$A_k$ there are
only countable many~$(j,m)$ appearing in~(\ref{eq:open stars
  inclusion}), we can assume after further refining the Borel 
partition~$(A_k)_{k\in\bbN}$ that for every~$k\in\bbN$ there are 
$j=j_k\in I$ and $m=m_k\in\bbN$ such that~(\ref{eq:open stars inclusion})
holds for a.e.~$x\in A_k$. 

Define $\psi^{(0)}$ on 
$\scrF$ by sending $(x,i_k,n_k)$ to $(x, j_k, m_k)$ for $x\in A_k$, and 
extend it to $\Sigma^{(0)}$ by equivariance. 
It follows (fiberwise) from 
the classical proof~\cite{spanier}*{Theorem~3 on~p.~127} that 
the map induced by $\psi$ between the $n$-fold fibered products 
of $\Sigma^{(0)}, \Phi^{(0)}$ restricts to
$\Sigma^{(n)}\rightarrow\Phi^{(n)}$ for any $n\ge 0$. 
The rest follows from Lemma~\ref{lem:homotopy for approximation}. 
\end{proof}

\subsection{$\scrR$-covers}\label{subsec:Measurable coverings and packings}
Throughout Subsection~\ref{subsec:Measurable coverings and packings}, 
we retain the notation of Assumption~\ref{setup:framework for
  simplicial complexes}. 
We introduce certain equivariant covers, so-called 
$\scrR$-covers, on the 
$\scrR$-space $X\times\widetilde{M}$ and a nerve 
construction. 

\begin{definition}\label{def:measurable covering and packing} 
Let $I$ be a free $\Gamma$-set. For $i\in I$ let 
$A_i\subset X$ be a Borel subset and 
and $U_i\subset\widetilde{M}$ be an open subset. \\
The family $\scrU=\{A_i\times U_i\}_{i\in I}$ is called an  
\textit{$\scrR$-cover of $X\times\widetilde{M}$} if 
\begin{enumerate}[a)]
\item $A_{\gamma i}=\gamma A_i$, $U_{\gamma i}=\gamma U_i$ for all 
  $\gamma\in\Gamma$ and $i\in I$, 
\item $\scrU_x\defq\{U_i;~x\in A_i\}_{i\in I}$ is locally finite in 
$\widetilde{M}$ for a.e. $x\in X$, 
\item for fixed $m\in\widetilde{M}$ and for a.e. $x\in X$ it is 
$(x,m)\in\bigcup_{i\in I}A_i\times U_i\subset X\times\widetilde{M}$.  
\end{enumerate} 
The family $\scrU=\{A_i\times U_i\}_{i\in I}$ 
is called an \textit{$\scrR$-packing} if 
\begin{enumerate}[a')]
\item $A_{\gamma i}=\gamma A_i$, $U_{\gamma i}=\gamma U_i$ for all 
  $\gamma\in\Gamma$ and $i\in I$, 
\item $U_i\cap U_j\ne\emptyset\Rightarrow\mu(A_i\cap A_j)=0$ for all 
$i\ne j$ in $I$. 
\end{enumerate}
\end{definition}

\begin{lemma}\label{lem:properties of R-coverings} 
Every $\scrR$-cover $\scrU=\{A_i\times U_i\}_{i\in I}$ of
$X\times\widetilde{M}$ such that $\mu(A_i)>0$ for every $i\in I$ 
has the following properties. 
\begin{enumerate}[a)]
\item The index set $I$ is countable. 
\item For every compact $K\subset\widetilde{M}$ 
there is a countable Borel partition 
$X=\bigcup_{j\in J} X_j$ such that for almost all  
$x,y\in X_j$ and every $k\in K$ we have: 
$(x,k)\in A_i\times U_i\Leftrightarrow(y,k)\in A_i\times U_i$. 
\item $\scrU_x$ is a cover of $\widetilde{M}=\{x\}\times\widetilde{M}$ 
for a.e. $x\in X$. 
\end{enumerate} 
\end{lemma}

\begin{proof} 
Let $T\subset\widetilde{M}$ be a countable dense subset, and 
set $I_m=\{i\in I;~m\in U_i\}$ for $m\in T$. 
Then $I=\bigcup_{m\in T}I_m$. 
For any $m\in T$ we define a Borel partition $X=\bigcup_{n\ge 1} X_m(n)$ by 
\begin{equation*}
X_m(n)=\Bigl\{x\in X;\#\{i\in I;~(x,m)\in A_i\times U_i\}=n\Bigr\}.
\end{equation*}
Note here that since $\scrU_x$ is locally finite, the sets 
$\{i\in I;~(x,m)\in A_i\times U_i\}$ are finite for a.e. $x\in X$. 
Each set 
$I_m(r, s)=\{i\in I_m;~\mu(X_m(r)\cap A_i)>\frac{1}{s}\}$ 
is finite since $\sum_{i\in I_m(r,s)}\mu(X_m(r)\cap A_i)\le
r\mu(X_m(r))<\infty$. That $I$ is countable follows from 
\begin{equation*}
I=\bigcup_{m\in T}\bigcup_{r,s\ge 1}I_m(r,s).
\end{equation*}

Set $I_K(x)=\{i\in I;~x\in A_i, K\cap U_i\ne\emptyset\}$. 
Since $\scrU_x$ is locally finite $I_K(x)$ is finite for a.e. $x\in
X$. In particular, $I_K$ ranges as a function of
$x$ in the countable 
set of finite subsets of $I$. 
Let $X=\bigcup_{j\in J}X_j$ be a countable Borel partition such that 
$I_K(x)$ is a constant set on each $X_j$. This proves the second
assertion. 

To prove the third assertion, consider a compact subset 
$K\subset\widetilde{M}$ such that $\Gamma K=\widetilde{M}$. 
Let $X=\bigcup_{j\in J} X_j$ of $X$ be as in the second assertion. 
Suppose $\scrU_x$ is not a cover,
\textit{i.e.} $\{x\}\times\widetilde{M}\not\subset\bigcup_{i\in I}A_i\times U_i$ 
on a subset of positive measure. This implies that 
$\{x\}\times K\not\subset\bigcup_{i\in I}A_i\times U_i$ on a
subset $Y\subset X$ of 
positive measure that has a non-trivial intersection with 
some $X_{j_0}$. Pick $y_0\in Y\cap X_{j_0}$. Let $m\in K$ be such that 
$(y_0,m)\not\in\bigcup_{i\in I}A_i\times U_i$. Then 
$(y,m)\not\in\bigcup_{i\in I}A_i\times U_i$ for a.e.~$y\in Y\cap X_{j_0}$ 
contradicting c) in Definition~\ref{def:measurable covering and
  packing}. 
\end{proof}

\begin{lemma}\label{lem:properties of R packings}
Let $\scrU=\{A_i\times U_i\}_{i\in I}$ be an $\scrR$-packing such that 
$\mu(A_i)>0$ for every~$i\in I$. Then 
\begin{enumerate}[a)]
\item $I$ is countable, and 
\item for a.e. $x\in X$ the sets in $\{U_i;~x\in A_i\}$ are pairwise
  disjoint. That is, $\scrU_x$ is a packing for a.e. $x\in X$. 
\end{enumerate}
\end{lemma}

\begin{proof}
Let $T$ and $I_m$ for $m\in T$ be like in the previous proof. 
Then $I(m)$ has to be countable since otherwise there would exist 
$i\ne j$ in $I$ with $\mu(A_i\cap A_j)>0$. Thus, $I$ is countable. 

Suppose there is a Borel subset $A\subset X$ with $\mu(A)>0$ such that 
for every $x\in A$ there are $i\ne j$ with $x\in A_i\cap A_j$ and 
$U_i\cap U_j\ne\emptyset$. We may assume that $I=\bbN$. 
For $x\in A$, let $i(x)<j(x)$ be minimal in $I$ with this property. 
Since $I\times I$ is countable, there is a Borel subset 
$B\subset A$ with $\mu(B)>0$ such that $(i(x),j(x))$ is constant 
for $x\in B$, which contradicts c') in Definition~\ref{def:measurable
  covering and packing}. 
\end{proof}

\begin{definition}\label{def:R-cover}
The \textit{nerve} $\scrN(\scrU)$ of an $\scrR$-cover 
$\scrU=\{A_i\times U_i\}_{i\in I}$ of $X\times\widetilde{M}$
is the $\scrR$-simplicial complex 
whose $0$-simplices are 
\begin{equation*}
\scrN(\scrU)^{(0)}=\{(x,i); x\in A_i,i\in I\}\subset X\times I.
\end{equation*}
and whose set of 
$n$-simplices $\scrN(\scrU)^{(n)}\subset\scrN(\scrU)^{(0)}
    \ast\ldots\ast\scrN(\scrU)^{(0)}$ is given by 
\begin{equation*}
\scrN(\scrU)^{(n)}=\Bigl\{(x,i_0,\ldots,i_n);~x\in\bigcap_{s=0}^n
A_{i_s},~\bigcap_{s=0}^nU_{i_s}\ne\emptyset, i_k\ne i_l\text{ for $k\ne l$}\Bigr\}.
\end{equation*}
The nerve $\scrN(\scrU)$ carries the $\scrR$-action that comes from 
the $\Gamma$-action $\gamma(x,i_0,\ldots,i_n)=(\gamma x,\gamma
i_0,\ldots,\gamma i_n)$ (cf.~Remark~\ref{rem: R action from Gamma
  action}). 
\end{definition}

\begin{remark}\label{rem:standard embedding of the nerve}
The map $\scrN(\scrU)\rightarrow X\times\bDelta(I),
(x,i_0,\ldots,i_n)\mapsto (x,i_0,\ldots,i_n)$, is a standard embedding 
in the sense of Section~\ref{subsec:standard embeddings}. 
\end{remark}

\begin{remark}
By Lemma~\ref{lem:properties of R-coverings}, $\scrU_x$ is a cover of
$\widetilde{M}$ for a.e. $x\in X$, 
and $\scrN(\scrU)_x$ is the nerve of the cover $\scrU_x$, which is
locally finite for a.e. $x\in X$ since $\scrU_x$ is so. 
\end{remark}

\begin{remark}\label{rem:nerve as R-space}
The $\scrR$-space obtained from the 
realization of $\scrN(\scrU)$ is also called the \textit{nerve of $\scrU$} 
and denoted by the same notation. It will be clear from the context
to which meaning we refer. 
\end{remark}

\begin{lemma}\label{lem:number of cells in the nerve}
Let $\scrU=\{A_i\times U_i\}_{i\in I}$ be an $\scrR$-cover. Let
$I'\subset I$ be a complete set of $\Gamma$-representatives. Set 
\begin{equation*}
f(i_0,\ldots,i_n)=\begin{cases}
                      \mu(A_{i_0}\cap\ldots\cap A_{i_n}) &\text{ if
                        $i_l\ne i_k$ for $l\ne k$ and 
                        $\bigcap_{k=0}^nU_{i_k}\ne\emptyset$,}\\
                      0 &\text{ otherwise.}
                  \end{cases}
\end{equation*}
Then we have 
(cf.~Definition~\ref{def:weighted number of simplices}) 
\begin{equation}\label{eq:number of cells in the nerve}
\sigma_n\bigl(\scrN(\scrU)\bigr)=\frac{1}{(n+1)!}\cdot\nu_t\bigl(\scrR\backslash\scrN(\scrU)^{(n)}\bigr)=
\sum_{i_0\in I'} \sum_{(i_1,\ldots,i_n)\in I^n}f(i_0,\ldots,i_n). 
\end{equation}
\end{lemma}

\begin{proof}
It is clear that 
\begin{equation}\label{eq:fundamental domain for nerve}
\scrF=\bigl\{(x,i_0,i_1,\ldots,i_n)\in\scrN(\scrU)^{(n)};~i_0\in I'\bigr\} 
\end{equation}
is an $\scrR$-fundamental domain of $\scrN(\scrU)^{(n)}$. Let 
$\scrF_x=\scrF\cap\scrN(\scrU)_x$. From
Definition~\ref{def:transversal measure} we see that 
\begin{equation}\label{eq:before fubini}
\nu_t\bigl(\scrR\backslash\scrN(\scrU)^{(n)}\bigr)=\int_X\#\scrF_x~d\mu(x).
\end{equation}
The right hand side in~(\ref{eq:before fubini}) coincides
with~(\ref{eq:number of cells in the nerve}) by Fubini's theorem 
applied to the product measure space $X\times I^{(n+1)}$ with the counting
measure on $I^{(n+1)}$. 
\end{proof}

For the following, recall that all metric notions on simplicial
complexes refer to the 
length metric that restricts to the standard Euclidean metric on
simplices.

\begin{lemma}\label{lem:map to nerve;finite index set}
Let $\scrU=\{A_i\times U_i\}_{i\in I}$ be an $\scrR$-cover of
$X\times\widetilde{M}$ such that $\Gamma\backslash I$ is finite. 
Then there is a geometric $\scrR$-map
$\phi:X\times\widetilde{M}\rightarrow\scrN(\scrU)$ such that 
\begin{enumerate}[a)]
\item $\phi$ is $\scrR$-simplicial (after a multiple barycentric
  subdivision of the domain). 
\item For a.e. $x\in X$ the preimage under $\phi_x$ 
of the open star of the vertex $i\in I$ in $\scrN(\scrU)_x$ 
is contained in $U_i$. 
\end{enumerate}
\end{lemma}

\begin{proof}
Let $K\subset\widetilde{M}$ be a compact set that contains the 
open $1$-neighborhood of a $\Gamma$-fundamental domain of
$\widetilde{M}$. 
Then the cover $\{\gamma
K;~\gamma\in\Gamma\}$ of $\widetilde{M}$ has Lebesgue number~$1$. 
Since $\Gamma\backslash I$ is finite and $\Gamma$ acts properly on
$\widetilde{M}$, there can be only finitely 
many $i\in I$ with $K\cap U_i\ne\emptyset$. In particular, there is 
a subset $X'\subset X$ of full measure such that 
only finitely many covers appear as the restriction of 
some $\scrU_x$ , $x\in X'$, to $K$. Each of them has positive Lebesgue 
number with respect to restricted metric on $K$. 
If $\epsilon'>0$ is the minimal such number, then 
$\scrU_x$ has Lebesgue number $\epsilon\defq\min\{\epsilon',1\}$ 
for a.e. $x\in X$. 

For each $U_i$, let $\bar{B}_{\epsilon/4}(\partial U_i)$ be the
closed $\epsilon/4$-neighborhood of the boundary $\partial U_i$. 
Set $V_i\defq U_i-\bar{B}_{\epsilon/4}(\partial U_i)$. Then $\scrV=\{A_i\times
V_i\}_{i\in I}$ is an $\scrR$-cover of $X\times\widetilde{M}$ such that 
$\scrV_x$ has Lebesgue number $\epsilon/4$ and $\bar{V}_i\subset
U_i$. Let $I'\subset I$ be a system of $\Gamma$-representatives. 
By Uryson's lemma, for each $i\in I'$  
there is a function $\tau_i: \widetilde{M}\rightarrow
[0,1]$ such that $\tau_i\vert_{\bar{V_i}}\equiv 1$ and
$\supp(\tau_i)\subset U_i$. Extend the definition of 
$\tau_i$ to all $i\in I$ by 
$\tau_{\gamma i}(m)=\tau_i(\gamma^{-1}m)$. Now define 
$\tau:X\times\widetilde{M}\rightarrow\scrN(\scrU)$ by 
\begin{equation*}
\tau(x,m)=\Bigl(x,\frac{1}{\sum_{i\in I}\chi_{A_i}(x)\tau_i(m)}\sum_{i\in
    I}\chi_{A_i}(x)\tau_i(m)i\Bigr). 
\end{equation*}
Here $\chi_{A_i}$ denotes the characteristic function of $A_i$. 
Then $\tau$ is proper as $\tau$ 
clearly satisfies the second assertion and has 
Lebesgue number~$\epsilon/4$. 
Obviously, $\tau$ is equivariant. Countable variance follows 
from Lemma~\ref{lem:properties of R-coverings}~b).  

By Theorem~\ref{thm: simplicial
  approximation}, $\tau$ possesses a simplicial approximation $\phi$,
which still satisfies the second assertion. 
\end{proof}

\begin{definition}\label{def:uniformly bounded}
An $\scrR$-cover 
$\scrU=\{A_i\times U_i\}_{i\in I}$ 
of $X\times\widetilde{M}$ 
is called \textit{uniformly bounded} if there is an $R>0$ such that 
$\diam(U_i)<R$ for every $i\in I$. 
\end{definition}

\begin{definition}
A geometric $\scrR$-map $\phi:\Sigma\rightarrow\Phi$ 
between $\scrR$-spaces $\Sigma$ and $\Phi$ 
is called \textit{metrically coarse} if for all 
$R>0$ there is an $S>0$ such that for a.e. $x\in X$ and for all
$m,m'\in\Sigma_x$ we have 
\[d(m,m')<R\Rightarrow d(\phi(m),\phi(m'))<S.\]
\end{definition}

\begin{theorem}\label{thm:extension theorem for coarse maps}
Assume that $M$ is aspherical. Let $\Sigma$ be a finite-dimensional 
$\scrR$-simplicial complex and 
$\Phi\subset\Sigma$ 
an $\scrR$-simplicial subcomplex that contains $\Sigma^{(0)}$. 
Let $\phi:\abs{\Phi}\rightarrow X\times\widetilde{M}$ 
be a metrically 
coarse geometric $\scrR$-map. Then $\phi$ can be extended to 
a metrically coarse geometric $\scrR$-map $\phi:\abs{\Sigma}\rightarrow
X\times\widetilde{M}$. 
\end{theorem}

\begin{proof}
Let $n\ge 1$. We show how to extend~$\phi$ from
$\abs{\Sigma^{(n-1)}}\cup\abs{\Phi}$ to
$\abs{\Sigma^{n}}\cup\abs{\Phi}\subset\abs{\Sigma}$.  
Let $\scrF$ be a fundamental domain for the $\scrR\times\symm(n+1)$-action
on $\Sigma^{(n)}\backslash\Phi^{(n)}$. 
There is a countable partition~$\scrF=\coprod_{j\in J}\scrF_j$ such
that the 
projection $\scrF\rightarrow X$ is injective on 
each $\scrF_j$~\cite{sauer(2005)}*{Lemma~3.1}. 
Let $X_j\subset X$ be the image of 
$\scrF_j$. For every $j\in J$ there is an  
embedding $X_j\times e_j\rightarrow\abs{\Sigma^{(n)}}$ with
$e_j=\abs{\bDelta^n}$ whose restriction to $\{x\}\times e_j$ is an 
affine isomorphism onto the (geometric) simplex given by the unique element 
in $\scrF_j\cap\Phi_x$ for $x\in X_j$. 
Identifying $X_j\times e_j$ with its image, we can write: 
\begin{equation*}
\absfix{\Sigma^{(n)}}\cup\absfix{\Phi}=\bigcup_{j\in J}\scrR.\bigl(X_i\times
e_j\bigr)\cup\absfix{\Sigma^{(n-1)}}\cup\absfix{\Phi}
\end{equation*}
Since $\phi$ (as defined so far) is of countable variance, 
we can assume, after possibly refining $(X_j)_{i\in J}$, that 
$\phi\vert_{X_j\times\partial e_j}$ is a product $\id\times f_j$ with 
continuous $f_j:\partial e_j\rightarrow\widetilde{M}$. 
Since $\phi$ is metrically coarse, the diameters of 
$f_j(\partial e_j)$ are uniformly bounded by a constant $R$. 
Note that $\widetilde{M}$ is 
uniformly contractible because $M$ is compact and aspherical. 
That is, every $R$-ball of $\widetilde{M}$ can
be contracted within an $S$-ball for some $S>0$. So we can find
an extension $F_j:e_j\rightarrow\widetilde{M}$ of $f_j$ for every
$j\in J$ such that $F_j(e_j)$ has diameter at most $S$, thus 
obtaining 
a metrically coarse extension of $\phi$ to $\abs{\Sigma^{(n)}}\cup\abs{\Phi}$ 
by equivariance. Countable variance of that extension is easy to
verify. 

It remains to show that
$\phi:\abs{\Sigma^{(n)}}\cup\abs{\Phi}\rightarrow X\times\widetilde{M}$ 
is fiberwise proper. This follows from the following general
statement~\citelist{\cite{bartels+rosenthal}*{Lemma~4.1}\cite{higson+roe}*{Lemma~3.3}}:  

Let $f:M\rightarrow N$ be a metrically coarse map between 
finite-dimensional, locally finite simplicial complexes $M$ and $N$ 
such that
$f\vert_{M^{(0)}}$ is proper, then $f$ is proper. 
Since every point of $M$ has distance at most~$1$ from $M^{(0)}$ and 
$M$ is locally finite, a subset $K\subset M$ is relatively compact if and only 
if $B_1(K)\cap M^{(0)}$ is finite. Here $B_1(K)$ denotes the
$1$-neighborhood of $K$. 
Now let $B(n,r)\subset N$ be the
ball of radius $r$ around some $n\in N$. Because $f$ is metrically coarse, 
there is an $S>0$ such that $B_1\bigl(f^{-1}(B(n,r))\bigr)$ is
contained in $f^{-1}\bigl(B(n,r+S)\bigr)$. 
Since $f$ is proper on 
the $0$-skeleton, 
\begin{equation*}
B_1\bigl(f^{-1}(B(n,r))\bigr)\cap M^{(0)}\subset f^{-1}\bigl(B(n,r+S)\bigr)\cap
M^{(0)}
\end{equation*}
is finite. Thus $f^{-1}(B(n,r))$ is relatively compact, and $f$ is proper. 
\end{proof}

\begin{lemma}\label{lem:homotopy retract}
Let $M$ be aspherical. Let $\scrU=\{A_i\times U_i\}_{i\in I}$ 
be a uniformly bounded $\scrR$-cover of $X\times\widetilde{M}$. 
Suppose $\scrN(\scrU)_x$ is finite-dimensional for a.e. $x\in X$. 
Then there is a metrically coarse, geometric
$\scrR$-map $\psi:\scrN(\scrU)\rightarrow X\times\widetilde{M}$. 
\end{lemma}

\begin{proof}
Let $I'\subset I$ a set of $\Gamma$-representatives. 
Pick for every $i\in I'$ a point $m_i\in U_i$, and extend the
definition of $m_i$ to $i\in I$ by $\gamma m_i=m_{\gamma i}$. 
Define $\psi:\scrN(\scrU)^{(0)}\rightarrow X\times\widetilde{M}$ by 
sending $(x,i)$ with $x\in A_i$ and $i\in I$ to $(x,m_i)$. For
a.e. $x\in X$ the map 
$\psi_x$ is proper since $\scrU_x$ is locally finite, and 
$\psi$ is metrically coarse since $\scrU$ is uniformly bounded. 
Now extend $\psi$ to $\scrN(\scrU)$ using Theorem~\ref{thm:extension
  theorem for coarse maps}. 
\end{proof}

\begin{lemma}\label{lem:geometric homotopy}
Assume that $M$ is aspherical. Let $\phi:X\times\widetilde{M}\rightarrow
X\times\widetilde{M}$ be a metrically coarse, geometric $\scrR$-map. 
Then there is a geometric $\scrR$-homotopy between
$\id_{X\times\widetilde{M}}$ and $\phi$.
\end{lemma}

\begin{proof}
Apply Theorem~\ref{thm:extension theorem for coarse maps} to 
extend the map $\phi\sqcup\id: X\times\widetilde{M}\times\{0,1\}\rightarrow
X\times\widetilde{M}$ to $X\times\widetilde{M}\times [0,1]$. 
\end{proof}

\section{Singular foliated homology}
\label{sec: homology theories of R-spaces}   

\subsection{The bundle of singular simplices of an
  $\scrR$-space}\label{subsec:bundle of singular simplices}

The \textit{bundle of singular $n$-simplices} of an $\scrR$-space $\Sigma$ is  
the disjoint union 
\[\sing_n(\Sigma)=\coprod_{x\in X}\sing_n(\Sigma_x), \]
where $\sing_n(\Sigma_x)$ is the set of singular $n$-simplices of 
$\Sigma_x$. The $\scrR$-action on $\Sigma$ induces one on
$\sing_n(\Sigma)$. Let $\pr:\sing_n(\Sigma)\rightarrow X$ denote the natural
projection. Since the $\scrR$-space $\Sigma$ has an
$\scrR$-fundamental domain (see Remark~\ref{rem:induced action}),  
$\sing_n(\Sigma)$ has one as
well: Take \textit{e.g.} the set of singular 
$n$-simplices whose first vertex lies in the fundamental domain of
$\Sigma$. 

Note that in the case $\Sigma=X\times\widetilde{M}$ we have 
$\sing_n(\Sigma)=X\times\sing_n(\widetilde{M})$. 

\begin{remark}\label{rem:combinatorial simplices in sing}
Let $\Phi$ be an $\scrR$-simplicial complex. Then there is a natural
inclusion $\Phi^{(n)}\subset\sing_n(\abs{\Phi})$ that maps 
$(v_0,\ldots,v_n)\in\Phi^{(n)}_x\subset
\Phi_x^{(0)}\ast\ldots\ast\Phi^{(n)}_x$ to the singular $n$-simplex 
$\Delta^n\rightarrow\abs{\Phi_x}, (t_0,\ldots,t_{n+1})\mapsto
t_0v_0+\ldots+t_nv_n$. 
\end{remark}

Of course, $\sing_n(\Sigma)$ for an $\scrR$-space $\Sigma$ 
is not an $X$-space in the sense of Definition 
since its fibers are in general uncountable. Although
Definition~\ref{def:transversal measure} thus cannot be applied, 
we will define a transversal measure for certain subsets 
of $\scrR\backslash\sing_n(\Sigma)$:  
From a standard embedding we obtain an injection 
\begin{equation}\label{eq:realization of standard embedding}
\Phi:\sing_n(\Sigma)\rightarrow X\times\sing_n(\bDelta(I\times
  \bbN))
\end{equation}
with respect to which we define the following notion. 
\begin{definition}\label{def:admissible subset}
A subset $W\subset\sing_n(\Sigma)$ is \textit{admissible}, if
there exists a countable subset 
$C\subset\sing_n(\bDelta(I\times\bbN))$ such that 
$\Phi(W)$ is a Borel subset of the $X$-space $X\times C$. 
A subset 
$W\subset\scrR\backslash\sing_n(\Sigma)$ is \textit{admissible} if 
its pullback to $\sing_n(\Sigma)$ is admissible. 
\end{definition}
By Lemma~\ref{lem:transversal measure well defined} below 
the property of being 
admissible and the Borel structure of an admissible set 
do not depend on the choice of~(\ref{eq:realization of standard
  embedding}). 
\begin{remark}\label{}
Obviously, intersections
and countable unions of admissible sets are admissible. An admissible
subset of $\sing_n(\Sigma)$ is an $X$-space with respect to
$\pr:\sing_n(\Sigma)\rightarrow X$
and carries, provided it
is $\scrR$-invariant, a discrete $\scrR$-action since
$\sing_n(\Sigma)$ possesses a fundamental domain. 
\end{remark}

\begin{definition}\label{def:natural measure on sing}
Let $\Sigma$ be an $\scrR$-space. 
For admissible $W\subset\sing_n(\Sigma)$ the function 
$x\mapsto\#\bigl(\pr^{-1}(x)\cap W\bigr)$ on $X$ is integrable, and 
\begin{equation*}
\nu(W)=\int_X\#\bigl(
\pr^{-1}(x)\cap W\bigr)d\mu(x)
\end{equation*}
is called the \textit{natural measure} of $W$. Of course, $\nu(W)$ 
equals the natural measure of $\Phi(W)\subset X\times C$ defined in 
Definition~\ref{def:measure on X spaces}. As in
Definition~\ref{def:transversal measure} one defines the
\textit{transversal measure} $\nu_t(W)$ of an admissible subset
$W\subset\scrR\backslash\sing_n(\Sigma)$. 
\end{definition}

\begin{lemma}\label{lem:transversal measure well defined}\hfill
\begin{enumerate}[a)]
\item Consider embeddings 
$\Phi_i:\sing_n(\Sigma)\rightarrow X\times
\sing_n(\bDelta(I_i\times\bbN))$, $i\in\{1,2\}$,
as in~(\ref{eq:realization of standard embedding}). Then
$\Phi_1\circ\Phi_2^{-1}:\sing_n(\Sigma)\rightarrow\sing_n(\Sigma)$ 
maps admissible to admissible sets, and is Borel on admissible sets. 
\item Let 
$\phi:\Sigma\rightarrow\Phi$ be a geometric $\scrR$-map. The 
map $\sing_n(\phi):\sing_n(\Sigma)\rightarrow\sing_n(\Phi)$ induced by
$\phi$ maps
admissible to admissible sets, and is Borel on admissible sets. 
\end{enumerate}
\end{lemma}

\begin{proof}
We leave it to the reader to verify the assertions using 
Lemma~\ref{lem:compacta and countable variance} and the fact that the 
image of a singular simplex is compact. 
\end{proof}

\subsection{Singular homology of $\scrR$-spaces}\label{subsec:singular
  homology}   

The goal of this section is to introduce 
the \textit{singular foliated homology} $\scrH_n(\Sigma)$ 
of an $\scrR$-space $\Sigma$. In spite of the notation one should 
think of $\scrH_n(\Sigma)$ as being a homology group of the foliated 
space $\scrR\backslash\Sigma$. 
In fact, $\scrH_n(\Sigma)$ will be a singular 
version of the sheaf-theoretic tangential homology of 
measured foliations~\cite{moore+schochet}. 
A homology class in
$\scrH_n(\Sigma)$ gives rise to a measurable family of classes in the 
locally finite homology of the leaves $\leaf_\Sigma(x)$ 
(see~Definition~\ref{def:leaf}). However, the functorial properties 
(and the actual definition) of $\scrH_\ast$ can be easier 
given for $\scrR$-spaces than for their quotients. \smallskip\\
In the following $\Sigma$ always denotes an $\scrR$-space. 
The restriction of the projection 
$\sing_n(\Sigma)\rightarrow\scrR\backslash\sing_n(\Sigma)$ 
to $\sing_n(\Sigma_x)$ is injective for $x\in X$. By the identification
$\leaf_\Sigma(x)\cong\Sigma_x$ we obtain an injection
$\sing_n(\leaf_\Sigma(x))\rightarrow\scrR\backslash\sing_n(\Sigma)$
that only depends on the $\scrR$-class of $x\in X$, that is, only 
on $\leaf_\Sigma(x)$ itself. Identifying $\sing_n(\leaf_\Sigma(x))$
with its image, we can view $\sing_n(\leaf_\Sigma(x))$ as a subset 
of $\scrR\backslash\sing_n(\Sigma)$. 

Notice that, if the $\scrR$-action comes from a $\Gamma$-action
(cf.~Remark~\ref{rem: R action from Gamma action}), then 
$\scrR\backslash\sing_n(\Sigma)=\Gamma\backslash\sing_n(\Sigma)$
and $\scrR\backslash\sing_n(\Sigma)=\Gamma\backslash\sing_n(\Sigma)$.  

\begin{definition}\label{def:foliated singular simplices}
A map $\sigma:A\rightarrow\scrR\backslash\sing_n(\Sigma)$, where 
$A\subset\bbR$ is a Borel subset of \textbf{finite} Lebesgue measure, 
is called a
\textit{foliated singular $n$-simplex of $\Sigma$}, if it has the 
following properties: 
\begin{enumerate}[a)]
\item The image $\im(\sigma)$ is admissible, 
and 
$\sigma:A\rightarrow\im(\sigma)$ is an $(\im(\sigma),\nu_t)$-space in the sense of 
Definition~\ref{def: X-space}.  
\item The Lebesgue measure coincides with the natural measure of 
the $(\im(\sigma),\nu_t)$-space $A$. 
\item The set of singular simplices 
$\im(\sigma)\cap\sing_n(\leaf_\Sigma(x))$ is locally finite in
$\leaf_\Sigma(x)$ for a.e. $x\in X$. 
\end{enumerate}
The set of foliated singular $n$-simplices is denoted by
$\scrS_n(\Sigma)$. 
\end{definition}

See Remark~\ref{rem:comments on definition of foliated simplex} below 
for comments. 

\begin{remark}\label{rem:finite fibers almost everywhere}
Since the natural measure of the $(\im(\sigma),\nu_t)$-space $A$ in the 
previous definition is finite, a.e. fiber of 
$\sigma:A\rightarrow\scrR\backslash\sing_n(\Sigma)$ is finite. 
\end{remark}

On each fiber we have the usual \textit{face operators}
$\partial_i:\sing_n(\Sigma_x)\rightarrow\sing_{n-1}(\Sigma_x)$
and \textit{degeneracy operators}
 $s_i:\sing_n(\Sigma_x)\rightarrow\sing_{n+1}(\Sigma_x)$
 for $i\in\{0,1,\ldots,n\}$. They induce $\scrR$-equivariant 
maps, denoted by the same symbols, 
 $\partial_i:\sing_n(\Sigma)\rightarrow\sing_{n-1}(\Sigma)$ and 
 $s_i:\sing_n(\Sigma)\rightarrow\sing_{n+1}(\Sigma)$.  
One immediately sees that $\partial_i,s_i$ map admissible sets to admissible
sets and are Borel on admissible sets. 

The set $\{\scrS_n(\Sigma)\}_{n\ge 0}$ becomes a \textit{simplicial
    set} by the face and degeneracy operators. Like for any simplicial
  set, there is an \textit{associated (unnormalized) chain complex 
$\Cn_\ast(\Sigma)$} with $\Cn_n(\Sigma)=\bbZ[\scrS_n(\Sigma)]$, and its 
differential $d:\Cn_n(\Sigma)\rightarrow \Cn_{n-1}(\Sigma)$ is the
alternating sum of face operators $d=\sum_{i=0}^n(-1)^i\partial_i$. 

\begin{definition}\label{def:function associated to a simplex}
For a singular foliated $n$-simplex $\sigma$ the function 
\[\omega_\sigma:\scrR\backslash\sing_n(\Sigma)\rightarrow\bbZ,~~
s\mapsto\#\sigma^{-1}(s),\] 
is called the \textit{multiplicity function}
of $\sigma$. We extend its definition linearly to elements in
$\bbZ[\scrS_n(\Sigma)]$. 
\end{definition}

\begin{remark}\label{rem:amibiguity in the definition of multiplicity function}
Note that $\omega_\sigma$ is supported in the admissible set 
$\im(\sigma)$. On $\im(\sigma)$, the multiplicity function
$\omega_\sigma$ is only well defined up to null-sets
(cf.~Remark~\ref{rem:finite fibers almost everywhere}). Subsequent
constructions that use multiplicity functions are insensitive to 
null-sets so that we can ignore this ambiguity. 
\end{remark}

Let $\rho\in\Cn_n(\Sigma)$. 
By condition $3$ in Definition~\ref{def:foliated singular simplices}, 
the formal sum 
\begin{equation}\label{eq:locally finite sum in the leaf}
\rho(x)=\sum_{s\in\sing_n(\leaf_\Sigma(x))}\omega_\rho(s)s
\end{equation} 
lies in $\Clf_n(\leaf_\Sigma(x))$, which is the chain group of locally finite
chains on $\leaf_\Sigma(x)$, for a.e. $x\in X$. The assignment 
$\Cn_n(\Sigma)\rightarrow\Clf_n(\leaf_\Sigma(x)), \rho\mapsto\rho(x)$,  
is compatible with the boundary operator meaning that the square 
\begin{equation}\label{eq:boundary operator and leaves}
\xymatrix{
\Cn_n(\Sigma)\ar[r]\ar[d]^d&\Clf_n\bigl(\leaf_\Sigma(x)\bigr)\ar[d]^d\\
\Cn_{n-1}(\Sigma)\ar[r]&\Clf_{n-1}\bigl(\leaf_\Sigma(x)\bigr)}
\end{equation}
is commutative for a.e. $x\in X$. 
Hence the chains whose multiplicity
function vanishes a.e. form a subcomplex. 

\begin{definition}\label{def:foliated chain group}
The \textit{singular foliated chain complex}
$\scrC_\ast(\Sigma)$ of $\Sigma$ 
is the quotient of $\Cn_\ast(\Sigma)$ by the subcomplex
of chains whose multiplicity function vanishes almost everywhere. 
The $n$-th homology of $\scrC_\ast(\Sigma)$ is denoted by 
$\scrH_n(\Sigma)$ and called the \textit{singular foliated
  homology} of $\Sigma$. 
\end{definition}

\begin{remark}\label{rem:comments on definition of foliated simplex}
In Definition~\ref{def:foliated singular simplices}, the condition 
that $\im(\sigma)$ is admissible and on the measures are needed for 
the approximation results of Lemma~\ref{lem:approximating chains} and
Theorem~\ref{thm:inclusion and support norm}. Without the condition of
locally finiteness, the formal chains in~(\ref{eq:locally finite sum
  in the leaf}) do not form a well-defined chain complex, which is
essential for defining $\scrC_\ast(\Sigma)$. 
\end{remark}

\subsection{Functoriality and homotopy invariance}\label{subsec:functoriality}

\begin{theorem}\label{thm:functoriality}
The singular foliated chain complex, thus the singular foliated
homology, 
of $\scrR$-spaces is functorial with
respect to geometric $\scrR$-maps.
\end{theorem}

\begin{proof}
Let $\phi:\Sigma\rightarrow\Phi$ be a geometric
$\scrR$-map. 
Since $\phi$ is fiberwise continuous, it induces a map 
$\sing_n(\phi):\sing_n(\Sigma)\rightarrow\sing_n(\Phi)$, which descends to 
a map on the quotients
$\scrR\backslash\sing_n(\phi):\scrR\backslash\sing_n(\Sigma)\rightarrow\scrR
\backslash\sing_n(\Phi)$.  
Let $\sigma:A\rightarrow\scrR\backslash\sing_n(\Sigma)$ 
be a singular foliated $n$-simplex. 
We have to show that 
$\bigl(\scrR\backslash\sing_n(\phi)\bigr)\circ\sigma$ is a singular
foliated $n$-simplex. 

Set $V=\im(\sigma)$, and let $\bar{V}$ denote
the pullback of $V$ to $\sing_n(\Sigma)$. The set $V$ is admissible. 
Moreover, $\bar{V}$ is an $(X,\mu)$-space and carries the natural
measure $\alpha=\nu^\mu$ (see Definition~\ref{def:measure on X spaces}). 
Let $W=\bigl(\scrR\backslash\sing_n(\phi)\bigr)(V)$, and let 
$\bar{W}$ be the pullback 
of $W$ to $\sing_n(\Phi)$. 
By Lemma~\ref{lem:transversal measure well defined}, $W$ and $\bar{W}$ are
admissible, and $\sing_n(\phi):\bar{V}\rightarrow\bar{W}$ is a map of 
$X$-spaces.
Let $\beta=\nu^\mu$ be the natural measure
of $\bar{W}$ as an $(X,\mu)$-space. Further, let $\alpha_t,\beta_t$
be the transversal measures of $\alpha, \beta$, respectively. 
Viewing $\bar{V}$ as an $(\bar{W},\beta)$-space via
$\sing_n(\phi)$ and, similarly, $V$ as an $(W,\beta_t)$-space, 
consider the natural measures $\nu^\beta,\nu^{\beta_t}$ on $\bar{V},V$
respectively. 
Since $\sing_n(\phi)\vert_V$ is a map of $X$-spaces, $\alpha=\nu^\beta$. 
The pullback of an $\scrR$-fundamental domain of $W$ to $V$ is one of
$V$, hence $\alpha_t=\nu^{\beta_t}$. Let $\lambda$ be the Lebesgue
measure of $A$. By definition of a singular foliated simplex, 
$\lambda=\nu^{\alpha_t}$. Natural measures are transitive in the sense 
that 
the natural measure $\lambda=\nu^{\alpha_t}=\nu^{(\nu^\beta_t)}$ 
on $A$ as an $V$-space via $\sigma$ 
equals the natural measure $\nu^{\beta_t}$ of $A$ as an $W$-space via
$\bigl(\scrR\backslash\sing_n(\phi)\bigr)\circ\sigma$. Thus,
$\lambda=\nu^{\beta_t}$ as
desired. 
Since $\phi$ is proper on
the fibers, the quotient map
$\scrR\backslash\Sigma\rightarrow\scrR\backslash\Phi$ 
is proper on leaves; hence 
$\bigl(\scrR\backslash\sing_n(\phi)\bigr)\circ\sigma$ 
satisfies the third property of
Definition~\ref{def:foliated singular simplices}. 

To sum up, $\phi$ induces a well-defined 
map $\scrS_n(\Sigma)\rightarrow\scrS_n(\Phi)$ and, by linear 
extension, $\Cn_n(\Sigma)\rightarrow\Cn_n(\Phi)$. The latter descends 
to a map $\scrC_n(\Sigma)\rightarrow\scrC_n(\Phi)$, which follows from 
the following observation, which we record in a Remark for later 
reference. 
\end{proof}

\begin{remark}\label{rem:functoriality on fibers}
The 
assignment 
$\Cn_n(\Sigma)\rightarrow\Clf_n(\leaf_\Sigma(x)), \rho\mapsto\rho(x)$ 
in~(\ref{eq:locally finite sum in the leaf}) 
descends to $\scrC_n(\Sigma)\rightarrow\Clf_n(\leaf_\Sigma(x))$. 
Let $\phi:\Sigma\rightarrow\Phi$ be a geometric $\scrR$-map.  
Then one easily verifies that 
the induced map $\scrC_n(\phi)$ is compatible with
$\rho\mapsto\rho(x)$ in the sense that the square (cf.~(\ref{eq:boundary
  operator and leaves})) 
\begin{equation}\label{eq:induced map on leaves}
\xymatrix{
\scrC_n(\Sigma)\ar[r]^-{x\mapsto\rho(x)}\ar[d]^{\scrC_n(\phi)}&
\Clf_n\bigl(\leaf_\Sigma(x)\bigr)\ar[d]^{\Clf_n(\phi)}\\
\scrC_n(\Phi)\ar[r]^-{x\mapsto\rho(x)}&\Clf_n\bigl(\leaf_\Phi(x)\bigr)}
\end{equation}
commutes for a.e. $x\in X$. 
\end{remark}

Next we present (in an informal way) 
a general principle that allows to transfer 
proofs from ordinary homology to singular foliated homology. 

\begin{remark}[Extension principle]\label{rem:extension principle}
Let $F_\ast$ be a functor 
\begin{equation*}
F_\ast:\{\text{top.~spaces}\}\rightarrow\{\text{chain complexes}\}
\end{equation*}
where $F_\ast(Y)$ is $\Cn_\ast(Y),\Cn_\ast(Y\times
[0,1])$ or $\Cn_{\ast+1}(Y)$. Instead of formalizing 
the extension principle in greater generality, 
we stick to these cases. 
Let $G_\ast$ be another such functor. Suppose there 
is natural transformation $\omega_\ast:F_\ast\rightarrow G_\ast$. Let 
$\Flf_\ast$ and $\Glf_\ast$ 
be the locally finite versions of $F_\ast$ and
$G_\ast$, respectively: If \textit{e.g.} $F_\ast(Y)=\Cn_\ast(Y)$, then
$\Flf_\ast(Y)=\Clf_\ast(Y)$.  By naturality we have 
a commutative 
diagram for a singular $p$-simplex $s:\Delta^p\rightarrow Y$: 
\begin{equation}\label{eq: acyclic model diagram}
\xymatrix{F_p(Y)\ar[r]^{\omega_p(Y)}&G_p(Y)\\
F_p(\Delta^p)\ar[u]^{F_p(s)}\ar[r]^{\omega_p(\Delta^p)}&G_p(\Delta^p)
\ar[u]^{G_p(s)}}
\end{equation}
If \textit{e.g.} $F_p(Y)=G_p(Y)=\Cn_p(Y)$, we have 
$\omega_p(Y)(s)=\Cn_p(s)\bigl(\omega_p(\Delta^p)(\id_{\Delta^p})\bigr)$. 
So \textquotedblleft supports of simplices are not enlarged\textquotedblright.
As a consequence (similar for all other examples of $F$ and $G$), 
$\omega_\ast$ naturally extends to 
$\Flf_\ast\rightarrow\Glf_\ast$. Now let 
\begin{equation*}
\scrF, \scrG:\{\text{$\scrR$-spaces}\}\rightarrow\{\text{chain
  complexes}\},~i\in\{0,1\},
\end{equation*}
be the corresponding versions 
of $F$ and $G$, respectively, for $\scrR$-spaces, that is,  
if \textit{e.g.} $F(Y)=\Cn_\ast(Y)$, then 
$\scrF(\Sigma)=\scrC_\ast(\Sigma)$ for an $\scrR$-space $\Sigma$. 

Then $\omega$ gives rise to 
a natural transformation
$\Omega:\scrF\rightarrow\scrG$ such that for every $\scrR$-space 
$\Sigma$ the diagram 
\begin{equation}\label{eq:abstract principle}
\xymatrix{
\scrF_p(\Sigma)\ar[r]^-{x\mapsto\rho(x)}\ar[d]^{\Omega(\Sigma)}&
\Flf_p\bigl(\leaf_\Sigma(x)\bigr)\ar[d]^{\omega(\leaf_\Sigma(x))}\\
\scrG_p(\Sigma)\ar[r]^-{x\mapsto\rho(x)}&\Glf_p\bigl(\leaf_\Sigma(x)\bigr)}
\end{equation}
commutes for a.e. $x\in X$ and all $p\ge 0$. 

We illustrate the 
idea for the case $F_p(Y)=G_p(Y)=\Cn_p(Y)$. 
Pick $a_1,\ldots, a_m\in\bbZ$ and 
$s_1,\ldots,s_m\in\sing_p(\Delta^p)$ such that 
$\omega_p(\Delta^p)=\sum_{i=1}^ma_is_i$. 
Let $\Sigma$ be an $\scrR$-space and 
$\sigma:A\rightarrow\scrR\backslash\sing_p(\Sigma)$ be a foliated
singular $p$-simplex. 
Then $\Omega$ is defined 
by 
$\Omega_p(\Sigma)(\sigma)=\sum_{i=1}^ma_i\rho_i$ where 
$\rho_i:A\rightarrow\scrR\backslash\sing_p(\Sigma)$ is given by 
the composition 
\begin{equation*}
A\xrightarrow{\sigma}\scrR\backslash\sing_p(\Sigma)\xrightarrow{s\mapsto
s\circ s_i}\scrR\backslash\sing_p(\Sigma).
\end{equation*}
\end{remark}

\begin{theorem}\label{thm: homotopy invariance}
Geometric $\scrR$-maps that are geometrically $\scrR$-homotopic
induce the same map in singular
foliated homology. 
\end{theorem}

\begin{proof}
It suffices to show that the two inclusion maps 
$i_0, i_1:\Sigma\rightarrow\Sigma\times [0,1]$ given by 
$i_0(x)=(x,0)$ and  $i_1(x)=(x,1)$  
induce chain homotopic maps on the singular foliated chain complexes. 

By~\cite{spanier}*{Theorem~3 on p.~174} there is natural transformation 
$D_\ast:F_\ast\rightarrow G_\ast$ of the functors 
$F_p(Y)=\Cn_p(Y)$ and $G_p(Y)=\Cn_{p+1}(Y\times [0,1])$ such that 
$D_\ast(Y)$ is a chain homotopy between the chain maps induced 
by $i_0$ and $i_1$. 
An application of the extension principle yields a corresponding chain homotopy 
for the singular foliated chain complexes. 
\end{proof}

\subsection{Support mass}\label{subsec:norms}

Throughout, let $\Sigma$ be an $\scrR$-space. 
By Definition~\ref{def:foliated singular simplices}, 
the support $\supp(\omega_\rho)\subset\scrR\backslash\sing_n(\Sigma)$ 
of the multiplicity function $\omega_\rho$ 
(Definition~\ref{def:function
  associated to a simplex}) of an element $\rho\in\scrC_n(\Sigma)$ 
is admissible, and $\omega_\rho$ is integrable 
with respect to the transversal measure $\nu_t$
(Definition~\ref{def:natural measure on sing}). 

\begin{definition}\label{def:support mass}
The \textit{support mass} of $\rho\in\scrC_n(\Sigma)$ is defined as 
\[\mass(\rho)=\nu_t\bigl(\supp(\omega_\rho)\bigr).\]
\end{definition}

Note that the support mass of a foliated singular simplex 
$\sigma:A\rightarrow\scrR\backslash\sing_n(\Sigma)$ is finite since 
$\mass(\sigma)\le\lambda(A)<\infty$ by Definition~\ref{def:foliated
  singular simplices}. 

\begin{lemma}\label{lem:support mass is a norm}
For all $\rho,\rho'\in\scrC_n(\Sigma)$ we have
$\mass(\rho+\rho')\le\mass(\rho)+\mass(\rho')$. 
\end{lemma}

\begin{proof}
This follows from
$\supp(\omega_{\rho+\rho'})=\supp(\omega_\rho+\omega_{\rho'})\subset
\supp(\omega_\rho)\cup\supp(\omega_{\rho'})$. 
\end{proof}

Note that the support mass is subadditive but not a norm since 
$\mass(d\rho)=\mass(\rho)$ for $\rho\in\scrC_n(\Sigma)$ and
$d\in\bbZ$. 

\begin{definition}\label{def:support mass on homology}
The \textit{support mass} of a homology class in $\scrH_n(\Sigma)$ is
the infimum of support masses of its representing cycles. 
We use the same notation $\mass(\kappa)$ for the support mass of 
a homology class $\kappa$. 
\end{definition}

\begin{remark}\label{rem:mass bounded by number of simplices}
Let $\Phi$ be an $\scrR$-simplicial complex, and let 
$\rho\in\scrC_n(\abs{\Phi})$. Assume that $\omega_\rho$ is supported 
in
$\scrR\backslash\Phi^{(n)}\subset\scrR\backslash\sing_n(\abs{\Phi})$
(for this inclusion see Remark~\ref{rem:combinatorial simplices in
  sing}). Then we have 
\begin{equation*}
\mass(\rho)\le
\nu_t\bigl(\scrR\backslash\Phi^{(n)}\bigr)=(n+1)!\cdot\sigma_n(\Phi).  
\end{equation*}
\end{remark}

\begin{lemma}\label{lem:norm-decreasing}
The homomorphisms $\scrC_\ast(\phi)$ and $\scrH_\ast(\phi)$ that are
induced by 
a geometric $\scrR$-map $\phi$ do not increase the support mass. 
\end{lemma}

\begin{proof} 
Let $\phi:\Sigma\rightarrow\Phi$ be a geometric $\scrR$-map and 
$\rho\in\scrC_n(\Sigma)$. Let $\rho'=\scrC_n(\phi)(\rho)$. 
We need only show that 
\begin{equation}\label{eq:to show}
\nu_t\bigl(\supp(\omega_{\rho'})\bigr)\le\nu_t\bigl(\supp(\omega_{\rho})\bigr).
\end{equation}\label{eq:subsets specific case}
First one verifies leafwise using 
diagram~(\ref{eq:induced map on leaves}) that 
\begin{equation}
\supp(\omega_{\rho'})\subset
\bigl(\scrR\backslash\sing_n(\phi)\bigr)\bigl(\supp(\omega_\rho)\bigr).  
\end{equation}
Then~(\ref{eq:to show}) is a 
consequence of the following general fact about transversal measures: 
By changing notation, 
let $\Sigma$ and 
$\Phi$ denote $(X,\mu)$-spaces with discrete $\scrR$-actions, and let 
$\phi:\Sigma\rightarrow\Phi$ be an $\scrR$-equivariant map of $X$-spaces. 
Let $A\subset\scrR\backslash\Sigma$ be a Borel subset. Then 
\begin{equation}\label{eq:subsets general case}
\nu_t\bigl((\scrR\backslash\phi)(A)\bigr)\le\nu_t(A). 
\end{equation}
The corresponding assertion $\nu(\phi(B))\le\nu(B)$ for
$B\subset\Sigma$ and the natural measure $\nu$ is clear from
Definition~\ref{def:measure on X spaces}. Now~(\ref{eq:subsets general
  case}) follows from the fact that if $\scrF$ is a
fundamental domain for $\Phi$, then 
$\phi^{-1}(\scrF)$ is one for $\Sigma$. 
\end{proof}

\begin{lemma}\label{lem:image of simplicial cycle}
Let $\Sigma,\Phi$ and $\phi$ as in Lemma~\ref{lem:realization of
  simplicial map is geometric}. 
Let $\kappa\in\scrC_n(\abs{\Sigma})$ be a cycle such that 
$\omega_\kappa$ is supported in $\scrR\backslash\Sigma^{(n)}$. 
Then $\scrH_n(\abs{\phi})([\kappa])\in\scrH_n(\abs{\Phi})$ 
has support mass at most $(n+1)!^2\cdot\sigma_n(\Phi)$. 
\end{lemma}

\begin{proof}
For every topological space $N$, 
there is chain map $Y_\ast(N):\Cn_\ast(N)\rightarrow\Cn_\ast(N)$ 
that realizes barycentric subdivision~\cite{bredon}*{p.~224}; 
$Y_\ast$ is natural in $N$ and thus a 
natural transformation in the sense of Remark~\ref{rem:extension
  principle}. 
The elementary properties of 
$Y_\ast$ are~\cite{bredon}*{p.~225}: 
\begin{enumerate}[a)]
\item Let $\sigma:\Delta^p\rightarrow N$ be degenerate in the sense
  that there is a surjective affine and vertex-preserving map 
$\alpha:\Delta^p\rightarrow\Delta^q$
with $q<p$ and a singular $q$-simplex 
$\sigma':\Delta^q\rightarrow N$ 
such that $\sigma=\sigma'\circ\alpha$. Then $\bigl(Y_p(N)\bigr)(\sigma)=0$. 
\item $\bigl(Y_p(\Delta^p)\bigr)(\id_{\Delta^p})$ 
is a linear combination of the affine 
simplices of the barycentric subdivision of $\Delta^p$. 
\item $Y_\ast$ is naturally chain homotopic to the identity. 
\end{enumerate}
The first property is not explicitly stated in~\cite{bredon} but 
follows from naturality and 
\begin{equation}\label{eq:signum of Y}
\bigl(Y_p(\Delta^q)\bigr)(\gamma\sigma)=\sign(\gamma)
\bigl(Y_p(\Delta^q)\bigr)(\sigma) 
\end{equation} 
for $\gamma\in\symm(q+1)$, $\sigma:\Delta^p\rightarrow\Delta^q$ and 
the natural $\symm(q+1)$-action on $\Delta^q$. As 
$Y_p(\Delta^q)$ is inductively defined~\cite{bredon}*{p.~224}, 
Equation~(\ref{eq:signum of Y}) can be easily proved by induction over $p$. 

By the extension principle (see Remark~\ref{rem:extension principle}), 
$Y$ gives rise to a natural transformation~$\scrY_\ast$ on the singular
foliated chain complex. 
Set $\rho=\scrC_n(\abs{\phi})(\kappa)$ for $\kappa$ as in the
hypothesis. Since $\phi$ is simplicial, we can decompose $\rho$ 
as $\rho=\rho_1+\rho_2$ such that $\omega_{\rho_1}$ is supported in 
the degenerate simplices and $\omega_{\rho_2}$ is supported 
in $\scrR\backslash\Phi^{(n)}$. The first and third property of~$Y_\ast$ and the 
analogous one of~$\scrY_\ast$ imply that the homology classes satisfy
\begin{equation*}
[\rho]=\bigl[(\scrY_n(\abs{\Phi}))(\rho)\bigr]=\bigl[(\scrY_n(\abs{\Phi})
)(\rho_2)\bigr]. 
\end{equation*}
By the second property the multiplicity function of 
$\bigl(\scrY_n(\abs{\Phi})\bigr)(\rho_2)$ is supported 
in~$\scrR\backslash\sd(\Phi)^{(n)}\subset\scrR\backslash\sing_n(\abs{\Phi})$. 
By Remark~\ref{rem:mass bounded by number of simplices} we finally get 
\begin{equation*}
\mass\bigl([\rho]\bigr)\le(n+1)!\cdot\sigma_n\bigl(\sd(\Phi)\bigr)=(n+1)!^2\cdot
\sigma_n(\Phi).  
\qedhere  
\end{equation*}
\end{proof}

For the rest of this Subsection retain the setting of 
Assumption~\ref{setup:framework for simplicial complexes}. 
The goal is to define a notion of support mass on 
$\Hn^\Gamma_n(\widetilde{M};\lz)$ and to compare it with the one 
on $\scrH_n(X\times\widetilde{M})$.  

\begin{definition}\label{def:norm on algebraic chain complex}
Let $p:\widetilde{M}\rightarrow M$ be the natural projection. 
The choice of 
a $\Gamma$-fundamental domain $\scrF\subset\sing_n(\widetilde{M})$ 
gives rise to an isomorphism
$\Cn_n(\widetilde{M})\xrightarrow{\cong}\bbZ\Gamma\otimes_\bbZ\Cn_n(M)$ 
of left $\bbZ\Gamma$-modules that maps $s\in\gamma\scrF$ to 
$\gamma\otimes p\circ s$. So we get the isomorphism 
\begin{equation}\label{eq:identification for support mass}
\lz\otimes_{\bbZ}\Cn_n(M)\cong \lz\otimes_{\bbZ\Gamma}\Cn_n(\widetilde{M}).
\end{equation}
For 
$s_1,\ldots,s_k\in\sing_n(M)$ with $s_i\ne s_j$ for
$i\ne j$ and $f_1,\ldots,f_k\in\lz$
we define 
\begin{equation*}
\mass\Bigl(\sum_{i=1}^kf_i\otimes_\bbZ s_i\Bigr)\defq
\sum_{i=1}^k\mu\bigl(\supp(f_i)\bigr) 
\end{equation*}
as the \textit{support mass} of the chain $\sum_if_i\otimes_\bbZ s_i$. 
The support mass of elements in 
$\lz\otimes_{\bbZ\Gamma}\Cn_n(\widetilde{M})$ is then defined
via~(\ref{eq:identification for support mass}). 
The support mass does not 
depend on the choice of $\scrF$. 
The support mass of a homology class of 
$\Hn^\Gamma_n(\widetilde{M};\lz)$ is defined as the infimum of support
masses of representing cycles. 
\end{definition}

In the sequel we use the term \textit{isometric} in the sense of 
\textit{preserving support masses}. 

\begin{lemma}\label{lem:inclusion-from-ordinary-to-foliated}
The chain homomorphism 
\begin{equation*}
\lambda_\ast:\lz\otimes_{\bbZ}\Cn_\ast(\widetilde{M})\rightarrow\scrC_\ast(X\times
\widetilde{M})
\end{equation*} 
that is uniquely determined by 
\begin{equation*}
\lambda_n(\chi_Y\otimes
s)=\Bigl(\sigma:j^{-1}(Y)\rightarrow\Gamma\backslash\bigl 
(X\times\sing_n(\widetilde{M})\bigr)\Bigr)\text{ with $\sigma(y)=(j(y),s)$.}
\end{equation*}
for $Y\subset X$ Borel, $s\in\sing_n(\widetilde{M})$ and a 
measure preserving Borel isomorphism $j:[0,1]\rightarrow
X$ does not depend on the choice of $j$. Further, 
$\lambda_\ast$ descends to an isometric map denoted by the same symbol  
\begin{equation*}
\lambda_\ast:\lz\otimes_{\bbZ\Gamma}\Cn_\ast(\widetilde{M})
\rightarrow\scrC_\ast(X\times\widetilde{M}). 
\end{equation*}
\end{lemma}

\begin{proof} This is a matter of straightforward verification. 
\end{proof}

\begin{lemma}\label{lem:approximating chains}
Let $\rho\in\scrC_n(X\times\widetilde{M})$ and $\epsilon>0$. Then 
there is $c\in\lz\otimes_{\bbZ\Gamma}\Cn_n(\widetilde{M})$ such that 
$\mass(\rho-\lambda_n(c))<\epsilon$. 
\end{lemma}

\begin{proof}
Without loss of generality, we can assume that
$\rho:A\rightarrow\scrR\backslash \bigl(
X\times\sing_n(\widetilde{M})\bigr)$ 
is a foliated singular simplex. 
Let $\scrF\subset\sing_n(\widetilde{M})$ be a $\Gamma$-fundamental domain. 
After identifying $X\times\scrF$ with
$\Gamma\backslash\bigl( X\times\sing_n(\widetilde{M})\bigr)$, 
there is a countable set $S=\{s_1,s_2,\ldots\}\subset\scrF$ such that 
$\im(\rho)\subset X\times S$. For $r,s>0$ set 
\begin{equation*}
W_{r,s}=\left\{(x,s_i);~x\in X,~\#\rho^{-1}(x,s_i)\le r,~1\le i\le
  s\right\}\subset X\times S.
\end{equation*}
Let $A_{r,s}=\rho^{-1}(W_{r,s})$. 
The sets $W_{r,s}$ form a directed system whose union is 
$\im(\rho)$. Thus for sufficiently large $r,s$ 
\begin{equation*}
\mass(\rho-\rho\vert_{A_{r,s}})=\mass(\rho\vert_{A-A_{r,s}})\le\lambda(A-A_{r,s})
\le\epsilon,
\end{equation*}
where $\lambda$ is the Lebesgue measure. 
It remains to show that $\rho\vert_{A_{r,s}}\in\im(\lambda_n)$. 
Fix $r,s$. There is a Borel partition $A_{r,s}=\bigcup_{k\le r,l\le
  s}B_{k,l}$ such that $\im(\sigma\vert_{B_{k,l}})\subset
X\times\{s_l\}$ and $\sigma\vert_{B_{k,l}}$ is
injective (theorem of
selection;~see~\textit{e.g.}~\cite{sauer(2005)}*{Lemma~3.1} 
for the version needed here). 
In particular, $\sigma\vert_{B_{k,l}}$ is a measure preserving Borel 
isomorphism $B_{k,l}\rightarrow X_{k,l}=X_{k,l}\times\{s_l\}$  
onto its image $X_{k,l}$. The chain 
\begin{equation*}
c_{r,s}=\sum_{k\le r,l\le s}\chi_{X_{k,l}}\otimes
s_l\in\lz\otimes_{\bbZ\Gamma}\Cn_n(\widetilde{M}).
\end{equation*}
satisfies $\lambda_n(c_{r,s})=\rho\vert_{A_{r,s}}$.  
\end{proof}

\begin{theorem}\label{thm:inclusion and support norm}
The homomorphism
$\Hn_n(\lambda_\ast):\Hn_n^\Gamma\bigl(\widetilde{M};\lz\bigr)\rightarrow
\scrH_n \bigl(X\times\widetilde{M}\bigr)$ 
is isometric. 
\end{theorem}

\begin{proof}
Let $\rho\in\scrC_n(X\times\widetilde{M})$ and
$c\in\lz\otimes_{\bbZ\Gamma}\Cn_n(\widetilde{M})$. Let $\epsilon>0$. 
Since $\lambda_n$ is isometric, $\scrH_n(\lambda_\ast)$ does not
increase support masses. It remains to show: If $\rho$ and 
$\lambda_n(c)$ are homologous, then there is a chain
$c_1\in\lz\otimes_{\bbZ\Gamma}\Cn_n(\widetilde{M})$ that is homologous
to $c$ and satisfies 
\begin{equation*}
\mass(c_1)<\mass(\rho)+\epsilon.
\end{equation*}
Let $\kappa\in\scrC_{n+1}(X\times\widetilde{M})$ 
be such that $\rho=\lambda_n(c)+d\kappa$. 
By Lemma~\ref{lem:approximating chains}, we find $c_0$ with 
$\mass(\lambda_{n+1}(c_0)-\kappa)<\epsilon/(n+1)$. 
Then $c_1=c+dc_0$ satisfies 
\begin{align*}
\mass(c_1)&=
\mass\bigl(\lambda_n(c_1)\bigr)\\
&\le\mass(\rho)+\mass\bigl(d(\lambda_{n+1}(c_0)-
\kappa)\bigr)\\
&\le
\mass(\rho)+(n+1)\mass\bigl(\lambda_{n+1}(c_0)-\kappa\bigr)<\mass(\rho)+
\epsilon. \qedhere
\end{align*}
\end{proof}

\subsection{Fundamental classes}
We retain the setting of Assumption~\ref{setup:framework for
  manifolds}. 

\begin{remark}[Fundamental classes]\label{rem:fundamental class in
    tensored complex} 
Let 
\begin{equation*}
j_\ast:\Cn_\ast(M)
=\bbZ\otimes_{\bbZ\Gamma}\Cn_\ast(\widetilde{M})\rightarrow\lz\otimes_{\bbZ\Gamma}
\Cn_\ast(\widetilde{M})
\end{equation*}
be the map coming from inclusion of constant functions. 
Let $\lambda_\ast$ be the map from 
Lemma~\ref{lem:inclusion-from-ordinary-to-foliated}. 
Let $[M]\in\Hn_n(M)$ be the fundamental class of $M$. 
By convention, the images 
$\Hn_n(j_\ast)([M])\in\Hn^\Gamma_n(\widetilde{M};\lz)$ and
$H_n(\lambda_\ast\circ j_\ast)([M])\in\scrH_n(X\times\widetilde{M})$
are also called \textit{fundamental classes}. 

Moreover, $\Hn_n(j_\ast)$ is an isomorphism if and only if 
the $\Gamma$-action on $(X,\mu)$ is ergodic: By equivariant 
Poincar\'{e} duality,
\[\Hn^\Gamma_n(\widetilde{M};\lz)\cong\Hn^0_\Gamma(\widetilde{M};
\lz)\cong\lz^\Gamma,\] 
and $\Hn_n(j_\ast)$ corresponds to the
inclusion $\bbZ\rightarrow\lz^\Gamma$ of constant functions under this 
isomorphism. 
\end{remark}

\begin{remark}[Foliated simplicial volume]\label{rem:connes foliated volume}
We point out the 
relation to Connes's simplicial volume for foliations. 
Suppose for the moment that the $\scrC_n(\Sigma)$ is defined by 
chains with \textit{real} coefficients, that is, as the quotient of
$\bbR[\scrS_n(\Sigma)]$ by the chains with a.e. vanishing 
multiplicity function. 
Consider the following norm on $\scrC_n(\Sigma)$ 
\[\abs{\rho}=\int_{\scrR\backslash\sing_n(\Sigma)}\abs{\omega_\rho}.\]
By Definition~\ref{def:foliated singular simplices}, the norm 
$\abs{\sigma}$ of a foliated singular simplex
$\sigma:A\rightarrow\scrR\backslash\sing_n(\Sigma)$ equals the
Lebesgue measure $\lambda(A)$ of $A$. The norm on $\scrC_n(\Sigma)$ 
induces a semi-norm on $\scrH_n(\Sigma)$ by taking the infimum of 
representing chains. Let $M$ be as in the previous remark. 
The norm of the fundamental class in $\scrH_n(X\times\widetilde{M})$ 
is the \textit{foliated simplicial volume} 
of the measured 
foliation $\Gamma\backslash (X\times\widetilde{M})$, which is 
attributed to Connes and described by Gromov
in~\cite{gromov(1991)}*{section~2.4.B}. 
\end{remark}

\section{Proof of Theorem~\ref{thm:folvol bound by packing
    inequalities}}\label{sec:proof of main inequality}  

\subsection{$\scrR$-covers with controlled multiplicity}
In this section we construct the $\scrR$-cover needed in the proof of
Theorem~\ref{thm:folvol bound by packing inequalities} following 
the strategy in Section~\ref{subsec:on the approach}. 
By convention, if $B(R)$ denotes a metric ball of radius $R$ in a metric space, 
then $B(R')$ denotes the concentric ball of radius $R'$. 

\begin{theorem}\label{thm:suitable-coverings}
Let $M$ be a Riemannian manifold and $(X,\mu)$ be a probability $\Gamma$-space 
as in 
Assumption~\ref{setup:framework for manifolds}. 
Assume further the following packing inequality: 
There is a constant 
$N_0\in\bbN$ such that each ball of radius~$1$ in 
$\widetilde{M}$ contains at most~$N_0$ disjoint balls of radius~$1/16$. 
Then there are 
countable families $\{A_i\}_{i\in I}$ of Borel subsets of $X$ and 
$\{B_i(3/16)\}_{i\in I}$ of balls of radius~$3/16$ in $\widetilde{M}$ 
such that 
\begin{enumerate}[a)]
\item $\scrU(3/16)\defq\{A_i\times B_i(3/16)\}_{i\in I}$ is an $\scrR$-covering of
$X\times\widetilde{M}$,
\item $\scrU(1/4)_{x}\defq\{B_i(1/4); x\in A_i\}_{i\in I}$ 
has multiplicity at most~$N_0$ for a.e.~$x\in X$. 
\end{enumerate}
\end{theorem} 

For the proof we need the following easy lemma. 
\begin{lemma}\label{lem:freeness lemma}
Let $(X,\mu)$ and $\Gamma$ be as in Assumption~\ref{setup:orbit
  equivalence relation}. Let $\gamma\in\Gamma\backslash\{1\}$. 
Let $A\subset X$ be a Borel set with
$\mu(A)>0$. Then there is a Borel subset $A'\subset A$ with
$\mu(A')>0$ such that $\mu(\gamma A'\cap A')=0$. 
\end{lemma}
\begin{proof}
By~\cite{vara}*{Theorem~3.2}, $X$ is equivariantly Borel isomorphic to 
a $\Gamma$-invariant Borel subset of 
a compact metric space $Y$ with a continuous $\Gamma$-action. Thus 
we may assume that $X\subset Y$, and the measure $\mu$ is extended to $Y$ 
by zero. By Ulam's theorem~\cite{cohn}*{Proposition~8.1.10 on p.~258},  
$\mu$ is regular on $Y$. 

Upon subtracting a null set from $A$, we can assume that 
$\gamma x\ne x$ for all $x\in A$. 
Next we show that there exists $x_0\in A$ such that $\mu(A\cap U)>0$ 
for every open neighborhood $U$ of $x_0$. Arguing by contradiction,
suppose that every $x\in A$ has an open neighborhood $U_x$ 
with $\mu(A\cap U_x)=0$. This yields $\mu(K)=0$ for every
compact $K\subset A$, thus $\mu(A)=0$ by regularity. 
Let $U, V$ be disjoint, open neighborhoods of $x_0$ and $\gamma x_0$,
respectively, such that $\gamma U\subset V$. Then $A'\defq A\cap U$
satisfies the conclusion. 
\end{proof}

\begin{proof}[Proof of Theorem~\ref{thm:suitable-coverings}]
An $\scrR$-packing by balls of radius $r$ is, by definition, an
$\scrR$-packing whose sets are of the type $A\times B$ where $B$ is a
ball of radius $r$. 
We say that an $\scrR$-packing $\scrU=\{A_i\times B_i(1/16)\}_{i\in
  I}$ of $X\times\widetilde{M}$ 
by balls of radius~$1/16$ is \textit{non-equivariantly maximal} 
if there does \textbf{not} exist a Borel subset
$A\subset X$ with $\mu(A)>0$ and a ball $B(1/16)\subset\widetilde{M}$ of
radius~$1/16$ such that 
\begin{equation}\label{eq:in-the-complement}
A\times B(1/16)\subset X\times\widetilde{M}\backslash\bigcup_{i\in I}
A_i\times B_i(1/16).
\end{equation}

\noindent\textit{First claim:} For every $A\subset X$ with $\mu(A)>0$ and 
every finite $F\subset\Gamma$ with $1\not\in F$ there exists $A'\subset A$ with 
$\mu(A')>0$ such that $\mu(A'\cap\gamma A')=0$ for 
all $\gamma\in F$. \smallskip\\
Let $F=\{\gamma_1,\ldots,\gamma_m\}$. 
Apply Lemma~\ref{lem:freeness lemma} repeatedly to obtain Borel sets
$A_1,\ldots, A_m$ of positive measure such that $A_{i+1}\subset A_i$
and $\mu(A_i\cap\gamma_i A_i)=0$. Set $A'\defq A_m$. \smallskip\\
\noindent\textit{Second claim:} If $\scrU$ is not non-equivariantly
maximal, then there exists 
a subset $A'\subset X$ of positive measure and a 
ball~$B(1/16)$ of radius~$1/16$ such that 
\begin{equation*}
A'\times B(1/16)\subset X\times\widetilde{M}\backslash\bigcup_{i\in I}
A_i\times B_i(1/16), 
\end{equation*} 
and 
$\scrU\cup\{\gamma A'\times\gamma B(1/16)\}_{\gamma\in\Gamma}$
is still an $\scrR$-packing.\smallskip\\
Pick $A\times B(1/16)$ as in~(\ref{eq:in-the-complement}). 
The set 
\begin{equation*}
\Gamma_B=\{\gamma\in\Gamma;~\gamma\ne 1,~\gamma B(1/16)\cap
B(1/16)\ne\emptyset\}
\end{equation*} 
is finite since $\Gamma$ acts properly on
$\widetilde{M}$. Now $A'$ is obtained from the first claim 
with $F=\Gamma_B$. \smallskip\\
The set of $\scrR$-packings on $X\times\widetilde{M}$ by balls of radius
$1/16$ is 
partially ordered  
as follows: $\scrU\le\scrV$ if and only if for every $A\times
B\in\scrU$ there exists $A'\subset X$ Borel such that $A=A'$ up 
to null-sets and $A'\times B\in\scrV$.  

Let $\{\scrU_k\}_{k\in K}$ be a totally ordered family of
$\scrR$-packings by balls of radius $1/16$. 
Let $\scrU_k=\{A_i\times B_i\}_{i\in I_k}$.
There exists an upper bound $\scrU$ 
of $\{\scrU_k\}_{k\in K}$: Set 
$I\defq\coprod_{k\in K}I_k$ (disjoint union). For $i,j\in I$,  
say $i\sim j$ if and only if $A_i=A_j$ up to null sets 
and $B_i=B_j$. Let $J$ be a $\Gamma$-invariant, complete set of 
$\sim$-representatives. Then 
$\scrU=\{A_j\times B_j\}_{j\in J}$ is an $\scrR$-packing and 
an upper bound of $\{\scrU_k\}_{k\in K}$. 

By Zorn's lemma there exists a maximal element $\scrU(1/16)=\{A_j\times
B_j(1/16)\}_{j\in J}$ (now a different~$J$). 
By the second claim, $\scrU(1/16)$ is non-equivariantly maximal. 
We may and will assume that $\mu(A_j)>0$ for every $j\in J$. 
By Lemma~\ref{lem:properties of R packings}, $J$ is countable. 

It remains to show that 
$\scrU(3/16)$ and $\scrU(1/4)$ have the stated properties. 
If there exist $m\in\widetilde{M}$ and $A\subset X$ with 
$\mu(A)>0$ such that 
\begin{equation*}
(x,m)\not\in\bigcup_{j\in J}A_j\times B_j(3/16)
\end{equation*}
for $x\in A$, then 
\begin{equation*}
A\times B(m,1/16)\cap\bigcup_{j\in J}A_j\times
B_j(1/16)=\emptyset
\end{equation*}
for the ball $B(m,1/16)$ of radius~$1/16$ around~$m$ 
contradicting non-equivariant maximality. 
It immediately follows from the packing inequality on $\widetilde{M}$
and the fact that $\scrU(1/16)_x$ is a packing for a.e.~$x\in X$ 
(Lemma~\ref{lem:properties of R packings}) that 
$\scrU(1/4)_x$ has multiplicity at most~$N_0$ for a.e.~$x\in X$. 
\end{proof}

\subsection{The map to the nerve}\label{subsec:map into nerve}

\begin{setup_subsection}
Retain the setting of Assumption~\ref{setup:framework
  for manifolds}. Further, we assume throughout this section:  
\begin{enumerate}[$\bullet$]
\item Let $M$ be equipped with a Riemannian metric 
such that the induced metric on the universal cover $\widetilde{M}$ 
satisfies the following packing inequality: 
There is a constant $N_0$ such that for $r<1$ 
each ball of radius $1$ in $\widetilde{M}$ 
contains at most $N_0r^{-n}$ balls of radius $r$. 
\item Let $\scrU=\scrU(1/4)=\{A_i\times B_i(1/4)\}_{i\in I}$ 
be the $\scrR$-cover constructed in 
Theorem~\ref{thm:suitable-coverings}, which has multiplicity 
$\le N_0$ (in fact, $\le 16^{-n}N_0$). 
\item We regard  
the nerve $\scrN(\scrU)$ as an 
$\scrR$-simplicial complex, and we write $\abs{\scrN(\scrU)}$ 
for the corresponding $\scrR$-space (realization), which 
we equip  with length metric 
(see Remark~\ref{rem:simplicial path on R spaces by convention}). 
Recall that there is a standard embedding $\scrN(\scrU)\subset X\times\bDelta(I)$ 
(see Remark~\ref{rem:standard embedding of the nerve}). 
\end{enumerate}
\end{setup_subsection}
We now define a geometric 
$\scrR$-map $\phi:X\times\widetilde{M}\rightarrow\abs{\scrN(\scrU)}$ 
such that $\phi_x$ is induced 
by a partition of unity of $\scrU_x$ for a.e.~$x\in X$. 
For $i\in I$ define $\phi_i:\widetilde{M}\rightarrow [0,1]$ by 
\[\phi_i(m)=\begin{cases}
                 1 & \text{if $m\in B_i(3/16)$}\\
              1-16d(m, B_i(3/16)) & \text{if $m\in
              B_i(1/4)\backslash B_i(3/16)$}\\
              0   & \text{if $m\not\in B_i(1/4)$},
            \end{cases}\]
and let

\begin{equation}\label{eq:map to the nerve}
\phi(x,m)=\left(x,\frac{1}{\sum_{i\in I}\chi_{A_i}(x)\phi_i(m)}\sum_{i\in
    I}\chi_{A_i}(x)\phi_i(m)i\right), 
\end{equation}
where $\chi_{A_i}$ denotes the characteristic function of $A_i$. 
Since $\scrU_x$ is locally finite for a.e.~$x\in X$, $\phi_x$ is
proper. Lemma~\ref{lem:properties of R-coverings} b) implies that 
$\phi$ is of countable variance, and equivariance of $\phi$ is
obvious. 

The goal of this subsection is to show the following theorem 
whose proof is given after a sequence of lemmas. 

\begin{theorem}\label{thm:homotopy further}
Let $\phi$ be the map in~(\ref{eq:map to the nerve}). 
There is a constant $\const_{n,N_0}$ that only depends on $n$ and
$N_0$ (but not on $M$) such that 
the image under $\scrH_n(\phi)$ 
of the fundamental class in $\scrH_n(X\times\widetilde{M})$ 
(cf.~Remark~\ref{rem:fundamental class in tensored complex}) 
has support mass at most $\const_{N_0,n}\vol(M)$. 
\end{theorem}

\begin{lemma}\label{lem:measure finiteness of nerve }
For every $k\ge 0$ we have $\sigma_k(\scrN(\scrU))<\infty$
(cf.~Definition~\ref{def:weighted number of simplices}). 
\end{lemma}

\begin{proof}
By compactness of $M$, we may 
choose a complete set $I'\subset I$ of $\Gamma$-representatives 
such that $\bigcup_{i\in I'}B_i(1/4)$ is relatively compact. 
Let $K$ be the compact closure of $\bigcup_{i\in I'}B_i(3/4)$. 
Let be $\scrF$ as in Equation~(\ref{eq:fundamental domain for nerve}) 
in the proof of Lemma~\ref{lem:number of cells in the nerve} 
with $U_i=B_i(1/4)$. If $(x,i_0,\ldots,i_k)\in\scrF_x$, then 
$B_{i_l}(1/4)\subset K$ for $l\in\{1,\ldots,k\}$. Since the 
sets $B_{i_l}(1/16)$ are disjoint (by the proof of
Theorem~\ref{thm:suitable-coverings}) and there is an upper bound for the
number of disjoint balls of radius $1/16$ that are contained in $K$,
there is a constant $0<C<\infty$ such $\#\scrF_x<C$ for a.e. $x\in X$. 
Then Equation~(\ref{eq:before fubini}) yields the conclusion. 
\end{proof}

\begin{lemma}\label{lem:differential-bounded-by-multiplicity}
There is a constant $\const_{N_0}$ that only depends on  
$N_0$ such that the map 
$\phi_x$ from~(\ref{eq:map to the nerve}) 
has Lipschitz constant $\le\const_{N_0}$ for a.e. $x\in X$. 
\end{lemma}

\begin{proof} 
For a.e. $x\in X$, 
$\phi_x$ is the map from $\widetilde{M}$ to the at most 
$N_0$-dimensional nerve of $\scrU_x$ induced by the partition of unity 
$\{\phi_i;~i\in I, x\in A_i\}$. A standard 
computation~\citelist{\cite{bartels+rosenthal}*{Lemma~4.6}
\cite{bell}*{Proposition~1}}    
shows that $\phi_x$ has a Lipschitz 
constant bounded in terms 
of~$N_0$ and the Lebesgue number of~$\scrU_x$ ($=1/16$ in our
case). 
\end{proof}

\begin{definition}\label{def:notation geometric simplex}
Let $\Sigma$ be an $\scrR$-simplicial complex and 
$s\in\Sigma^{(k)}$. We let $\abs{s}\subset\abs{\Sigma}$ 
denote the geometric $k$-simplex corresponding to the (combinatorial)
$k$-simplex $s\in\Sigma^{(k)}$. 
The interior of $\abs{s}$ is denoted by
$\abs{\mathring{s}}$. 
\end{definition}
For the proof of Theorem~\ref{thm:homotopy further}
we need the following general definition and lemma. 

\begin{definition}\label{def:family of projectors}
Let $\epsilon\ge 0$. 
Let $\omega:\abs{\Phi}\rightarrow\abs{\Psi}$ be an geometric
$\scrR$-map. 
Let $\Sigma\subset\Psi^{(k)}$ be an $\scrR$-invariant Borel subset. 
A family 
$\{P^s\}_{s\in\Sigma}$ of points $P^s\in\abs{\mathring{s}}$ is
called a \emph{family of $\epsilon$-projectors for $\omega$} if the
following conditions hold:  
\begin{enumerate}[a)]
\item $\{P^s\}_{s\in\Sigma}\subset\abs{\Psi}$ is admissible. 
\item distance $d(P^s,\im(\omega)\cap\abs{s})>\epsilon$
  for a.e. $s\in\Sigma$. 
\item $\{P^s\}_{s\in\Sigma}\subset\abs{\Psi}$ is an $\scrR$-invariant subset. 
\end{enumerate}
\end{definition}  
For a) notice that $\{P^s\}$ can be seen as a subset of
$\sing_0(\abs{\Psi})$ and \textit{admissible} is understood in the
sense of Definition~\ref{def:admissible subset}. 
In b) we set $d(P_x^s,\emptyset)=\infty$. 

\begin{lemma}\label{lem:family of projectors}
Let $\omega:\abs{\Phi}\rightarrow\abs{\Psi^{(k)}}$ be a geometric
$\scrR$-map that lands in the $k$-skeleton,
and let $\Sigma\subset\Psi^{(k)}$ be an $\scrR$-invariant Borel subset. 
Let $\{P^s\}_{s\in\Sigma}$ be a family of $\epsilon$-projectors for
$\omega$. Define  
$\psi:\abs{\Phi}\rightarrow\abs{\Psi}$ to be the
map obtained from $\omega$ by post-composition 
with the radial
projections (within $\abs{s}$)  
from $P^s$ to the boundary $\partial\abs{s}$ for every $s\in\Sigma$. 
We say that $\psi$ \emph{is obtained from $\omega$ via $\{P_x^s\}$}. 
Then the following holds or holds a.e., respectively.   
\begin{enumerate}[a)]
\item $\psi$ is a geometric
$\scrR$-map. 
\item Assume that $\epsilon>0$. 
If $\omega_x$ has a Lipschitz constant $C>0$, 
then $\psi_x$ has a Lipschitz constant $C/\epsilon$.
\item The $\psi_x$-preimage of the open star
of any vertex is contained in the corresponding
$\omega_x$-preimage. 
\item There is geometric $\scrR$-homotopy between 
$\omega$ and $\psi$. 
\end{enumerate}
\end{lemma}

\begin{proof}
The proof that $\psi$ is of countable variance is straightforward and
follows from the fact that $\{P^s\}$ is admissible. It is clear that 
$\psi$ is $\scrR$-equivariant and continuous on fibers. Since
$\psi_x$ is an approximation of $\omega_x$, 
Lemma~\ref{lem:homotopy for approximation} 
implies that $\psi$ is a geometric $\scrR$-map, which is geometrically
$\scrR$-homotopic to $\omega$. For the same reason 
assertion c) is true. 
The radial projections $\abs{s}\backslash
B_\epsilon(P^s)\rightarrow\partial\abs{s}\subset\abs{s}$ have
Lipschitz constant $\le\epsilon^{-1}$. Using the fact that
the metric on $\abs{\Psi}$ is a length metric, one easily sees that 
$\psi_x$ has Lipschitz constant $\le C/\epsilon$ for a.e. $x\in X$. 
\end{proof}

\begin{lemma}\label{lem:family of projectors exists}
Let $\omega:X\times\widetilde{M}\rightarrow\abs{\scrN(\scrU)^{(k)}}$ be 
a geometric $\scrR$-map such that $\omega_x^{-1}(\openstar(i))\subset
B_i(1/4)$ for a.e. $x\in X$ and every vertex $i\in I$ in $\scrN(\scrU)_x$. 
Let $\Sigma\subset\scrN(\scrU)^{(k)}$ be an $\scrR$-invariant subset
such that 
\begin{equation*}
\sup_{z\in\abs{\mathring{s}}}d(z,
\im(\omega)\cap\abs{s})\ge\epsilon\text{ for $s\in\Sigma$}. 
\end{equation*}
Then there exists
a family of $\epsilon/2$-projectors for $\omega$. 
\end{lemma}
\begin{proof}
By compactness of $M$, we may 
choose a complete set $I'\subset I$ of $\Gamma$-representatives 
of $I$ such that $\bigcup_{i\in I'}B_i(1/4)$ is relatively compact. 
Let $K$ be the compact closure of $\bigcup_{i\in I'}B_i(3/4)$. Let
$\scrF\subset\Sigma$ be an $\scrR$-fundamental domain such that 
every $s\in\scrF$ has a vertex in $I'$. 
Referring to the embedding $\scrN(\scrU)\subset X\times\bDelta(I)$ 
of Remark~\ref{rem:standard embedding of the nerve}, $\scrF$ is a countable, 
disjoint union of Borel sets $X_p\times\{s_p\}$, $p\ge 1$, 
where $X_p\subset X$ and $s_p\in\bDelta(I)^{(k)}$ has at least one
vertex in $I'$. We have 
$\omega_x^{-1}(\abs{s})\subset K$, thus 
\begin{equation}\label{eq:intersection - compact set suffices}
\omega_x(\widetilde{M})\cap\abs{s}=\omega_x(K)\cap\abs{s}, 
\end{equation}
for a.e. $x\in X$ and every $s\in\scrF_x$. 
Let $X=\bigcup_{j\in J}X_j$ be a countable Borel partition such that 
every restriction $\omega\vert_{X_j\times K}$ is a product
map. 
By~(\ref{eq:intersection - compact set suffices}) the set 
\begin{equation*}
M(j,p)\defq\omega_x(\widetilde{M})\cap\abs{s_p}
\end{equation*}
is constant for a.e. $x\in X_j\cap X_p$. 
By assumption we can pick a point $P_j^p\in\abs{\mathring{s_p}}$ for
every $j\in J$ and every $p\ge 1$ such that 
\begin{equation*}
d(P_j^p,M(j,s))>\epsilon/2. 
\end{equation*}
For $s=(x,s_p)\in\scrF$ with $x\in X_j\cap X_p$ define 
$P^s\defq P_j^p$. Extend the definition of $P^s$ to $s\in\Sigma$ by
equivariance. 
Then $\{P^s\}$ is admissible and a 
family of $\epsilon/2$-projectors for $\omega$. 
\end{proof}

The following is a version of~\cite{gromov(1982)}*{Lemma D, Section~3.4} in the 
context of $\scrR$-spaces. 

\begin{lemma}\label{lem:map-to-classifying-space} 
There exists a geometric $\scrR$-map $\psi:X\times\widetilde{M}\rightarrow
\abs{\scrN(\scrU)^{(n)}}$ into the $n$-skeleton such that 
for a.e. $x\in X$
\begin{enumerate}[a)]
\item $\psi$ and $\phi$ as defined in~(\ref{eq:map to the nerve}) 
are geometrically $\scrR$-homotopic, 
\item $\psi_x$ has a Lipschitz constant 
$C=C(n,N_0)$ that only depends on $n$ and $N_0$, and 
\item the $\psi_x$-preimage of the open star of the vertex $i\in I$ 
in $\abs{\scrN(\scrU)_x}$ is contained in $B_i(1/4)$. 
\end{enumerate}
\end{lemma}

\begin{proof}
Since $\scrU_x$ has multiplicity $\le N_0$, the map $\phi$ lands in the 
$N_0$-skeleton. If $N_0\le n$, setting $\psi=\phi$ will do. Let 
$N_0>n$. 
We inductively construct geometric $\scrR$-maps 
\begin{equation*}
\psi^{(k)}:X\times\widetilde{M}\rightarrow\abs{\scrN(\scrU)^{(N_0-k)}}
\end{equation*}
for $k=0,1,\ldots, N_0-n$ satisfying a), b) and c). 
Set $\psi^{(0)}=\phi$. Let $0\le k<N_0-n$. 
Assume $\psi^{(k)}$ is already defined; we will define 
$\psi^{(k+1)}$ as a map obtained from $\psi^{(k)}$ via a family of 
$\epsilon$-projectors $\{P^s\}_{s\in\scrN(\scrU)^{(N_0-k)}}$ for
$\psi^{(k)}$. 
If such a family $\{P^s\}$ exists with an $\epsilon>0$ that only
depends on $n$ and $N_0$, then the resulting map 
$\psi^{(k+1)}$ has the desired properties by Lemma~\ref{lem:family of
  projectors}. 
For $s\in\scrN(\scrU)_x^{(N_0-k)}$ 
let $M(x,s)=\psi_x^{(k)}(\widetilde{M})\cap\abs{s}$. To apply
Lemma~\ref{lem:family of projectors exists} and finish the proof, it
remains to show that for a.e. $x\in X$ and every
$s\in\scrN(\scrU)_x^{(N_0-k)}$ 
\begin{equation}\label{eq:epsilon}
\epsilon_{x,s}\defq\sup_{z\in\abs{\mathring{s}}} d\bigl(z,M(x,s)\bigr)
\end{equation}
is bounded from below in terms of $n$ and $N_0$. We set $\epsilon_{x,s}=\infty$
if $M(x,s)=\emptyset$. 
We may and will assume that $\epsilon_{x,s}<1/4$. 
Since the $n$-dimensional Hausdorff measure of 
$\psi_x^{(k)}(\widetilde{M})\cap\abs{\mathring{s}}$ is finite 
by the area formula for Lipschitz maps~\cite{morgan}*{Chapter~3} whilst the 
$n$-dimensional Hausdorff measure of 
$\abs{\mathring{s}}$ is infinite, $M(x,s)$ must miss a point 
in $\abs{\mathring{s}}$, thus $\epsilon_{x,s}>0$. 
There is a constant $D>0$ only depending on
$n$ and $N_0$ such that there are 
\begin{equation}\label{eq:estimate on k}
m\ge D\epsilon_{x,s}^{-(N_0-k)}
\end{equation}
disjoint $\epsilon_{x,s}$-balls in $\scrN(\scrU)^{(N_0-k)}$  
with centers in $M(x,s)$: To see this, 
pick a maximal packing by $\epsilon_{x,s}$-balls 
$B_1'(\epsilon_{x,s}),\ldots, B_m'(\epsilon_{x,s})$ whose centers lie in $M(x,s)$.  
By~(\ref{eq:epsilon}), $\abs{s}$ is covered by 
$B_1'(3\epsilon_{x,s}),\ldots, B_m'(3\epsilon_{x,s})$. A volume estimate 
yields~(\ref{eq:estimate on k}) for $D$ only depending on $n,N_0$. 

Each open ball $B_k'(\epsilon_{x,s})$ lies in the union of open stars 
of vertices of $s$. Hence 
the $\psi_x^{(k)}$-preimages of the open balls 
$B_1'(\epsilon_{x,s}),\ldots, B_m'(\epsilon_{x,s})$ 
are contained in a ball of radius~$3/4$. 
Each preimage contains a ball of
radius $r=\epsilon_{x,s} C^{-1}$, where 
$C$ is the Lipschitz constant of $\psi^{(k)}$. 
By the packing inequality, $m\le N_0r^{-n}=N_0C^n\epsilon_{x,s}^{-n}$. Combined 
with~(\ref{eq:estimate on k}) one obtains 
\begin{equation*}
\epsilon_{x,s}\ge\left(DN_0^{-1}C^{-n}\right)^{1/(-n+N_0-k)}>0.
\end{equation*}
\end{proof}

\begin{proof}[Proof of Theorem~\ref{thm:homotopy further}]
Let $\psi$ be the map from Lemma~\ref{lem:map-to-classifying-space}. 
Let $\epsilon>0$ be given. 
For $m\ge 1$ define 
\begin{align*}
\Sigma&\defq\{s\in\scrN(\scrU)^{(n)};~\abs{\mathring{s}}\not\subset\im(\psi)\}\\
\Sigma_m&\defq\{s\in\scrN(\scrU)^{(n)};~\sup_{z\in\abs{\mathring{s}}}d(z,\im(\psi)\cap
\abs{s})\ge 1/m\}.\notag  
\end{align*}
Since $\im(\psi)\cap\abs{s}$ is closed by properness of $\psi$, 
$s\in\Sigma$ implies that there is $m\ge 1$ with $s\in\Sigma_m$. 
Let $\nu_t$ denote the measure defined in
Definition~\ref{def:transversal measure}. We obtain that 
\begin{equation}\label{eq:approximation in measure}
\lim_{m\rightarrow\infty}\nu_t(\scrR\backslash\Sigma_m)=
\nu_t(\scrR\backslash\Sigma).  
\end{equation}
Since
$\nu_t(\scrR\backslash\Sigma)\le\nu_t(\scrR\backslash\scrN(\scrU))<\infty$
by Lemma~\ref{lem:measure finiteness of nerve },
there is $m\ge 1$ large enough such that 
$\nu_t(\scrR\backslash(\Sigma\backslash\Sigma_m))<\epsilon$. 
By Lemma~\ref{lem:family of projectors exists} there is a family 
$\{P^s\}_{s\in\Sigma_m}$ of $1/(2m)$-projectors for $\psi$. Let 
$\psi'$ be the map obtained from $\psi$ via
$\{P^s\}_{s\in\Sigma_m}$. 
Define 
\begin{align*}
\Sigma_{hit}&\defq\{s\in\scrN(\scrU)^{(n)};~\abs{\mathring{s}}\cap\im(\psi')\ne
\emptyset\}.\\
\Sigma_{full}&\defq\{s\in \scrN(\scrU)^{(n)};~\abs{\mathring{s}}\subset\im(\psi')\}
\end{align*}
For $s\in\scrN(\scrU)^{(n)}$ define   
\begin{equation*}
D(s)\defq\psi'^{-1}(\abs{\mathring{s}}). 
\end{equation*} 
If $s\in\Sigma_{full}$, we have 
$D(s)=\psi^{-1}(\abs{\mathring{s}})$ and 
$\psi'\vert_{D(s)}=\psi\vert_{D(s)}$. 
In particular, 
if $s\in\Sigma_{full}$, then $\psi'\vert_{D(s)}$ has the Lipschitz constant $C$ 
from Lemma~\ref{lem:map-to-classifying-space}, and by 
the area formula~\cite{morgan}*{Chapter~3}
\begin{equation}\label{eq: first volume estimate}
\vol(\Delta^n)=\vol(\psi'(D(s)))\le C^n\vol(D(s)). 
\end{equation}
By definition of $\psi'$ we have 
\begin{equation}\label{eq:set of hit simplices}
\Sigma_{hit}\subset\Sigma\backslash\Sigma_m\cup\Sigma_{full}. 
\end{equation}
We equip $X\times\widetilde{M}$ with the product of $\mu$ and 
the Riemannian measure on $\widetilde{M}$. 
Any $\scrR$-fundamental domain of $X\times\widetilde{M}$
has measure $\vol(M)$. Let $\scrF$ be a 
$\symm(n+1)\times\scrR$-fundamental domain of
$\Sigma_{hit}$. As usual, let
$\scrF_x=\scrF\cap\scrN(\scrU)^{(n)}_x$ for $x\in X$. 
Since the disjoint union $\bigcup_{s\in\scrF}D(s)$ is contained  
in an $\scrR$-fundamental domain, its measure is at most $\vol(M)$. 
Fubini's theorem yields
\begin{equation}\label{eq:second volume estimate}
\int_X\sum_{s\in\scrF_x}\vol(D(s))d\mu(x)
\le\vol(M). 
\end{equation}
By Lemma~\ref{lem:family of projectors} the maps $\psi'_x$, $x\in X$, have a
Lipschitz constant $Cm$, so $\psi'$ has a positive Lebesgue
number (see Definition~\ref{def:lebesgue number}). 
Upon subdividing the triangulation on $M$, 
the map $\psi'$
has an $\scrR$-simplicial approximation $\psi''$  
by Theorem~\ref{thm: simplicial approximation}. 
Note that 
$\im(\psi'')^{(n)}\subset\Sigma_{hit}$. Thus, 
\begin{align*}
\sigma_n(\im(\psi''))&=\nu_t\bigl(\symm(n+1)\times\scrR\backslash\im(\psi'')\bigr)\\
 &\le\nu_t(\symm(n+1)\times\scrR\backslash\Sigma_{hit})&\\
 &\le\nu_t\bigl(\symm(n+1)\times\scrR\backslash(\Sigma\backslash\Sigma_m)\bigr)+\nu_t\bigl(\symm(n+1)\times\scrR
 \backslash (\Sigma_{hit}\cap\Sigma_{full})\bigr)  &\text{by~(\ref{eq:set of hit simplices})}\\ 
 &\le\frac{\epsilon}{(n+1)!}+\int_X\#(\scrF_x\cap\Sigma_{full})d\mu(x)&\\
 &\le\frac{\epsilon}{(n+1)!}+\int_X\sum_{s\in\scrF_x}C^n 
 \frac{\vol(D(s))}{\vol(\Delta^n)}d\mu(x)&\text{by~(\ref{eq:  first volume
     estimate})}\\ 
 &\le\frac{\epsilon}{(n+1)!}+C^n\frac{\vol(M)}{\vol(\Delta^n)}.&\text{by~(\ref{eq:second
     volume estimate})}
\end{align*}
We apply Lemma~\ref{lem:image of simplicial cycle} 
to $\Sigma=X\times\widetilde{M}$, 
$\Phi=\im(\psi'')$ and the map $\psi''$. The triangulation of 
$M$ naturally gives rise to a cycle $\kappa'\in\Cn_n(M)$ representing 
the fundamental class. 
Let $\kappa$ be the image of $\kappa'$ under 
\begin{equation*}
\Cn_n(M)=\bbZ\otimes_{\bbZ}\Cn_n(\widetilde{M})\xrightarrow{j_n}\lz\otimes_{\bbZ\Gamma}
\Cn_n(\widetilde{M})\xrightarrow{\lambda_n}\scrC_n(X\times\widetilde{M}), 
\end{equation*}
where $j_n$ is the map from Remark~\ref{rem:fundamental class in
    tensored complex} and $\lambda_n$ is the map from
Lemma~\ref{lem:inclusion-from-ordinary-to-foliated}. 
Then $\kappa$ satisfies the hypothesis 
of Lemma~\ref{lem:image of simplicial cycle}. 
Since $\epsilon>0$ can be taken arbitrarily small in the estimate
above and $\scrH_n(\phi)=\scrH_n(\psi'')$, 
Lemma~\ref{lem:image of simplicial cycle} yields 
Theorem~\ref{thm:homotopy further}. 
\end{proof}

\subsection{Conclusion of proof}\label{subsec:homotopy retract}
\begin{proof}[Proof of Theorem~\ref{thm:folvol bound by packing
    inequalities}]
For $M$ as in the hypothesis and $\Gamma=\pi_1(M)$ pick a 
probability space $(X,\mu)$ as in Assumption~\ref{setup:orbit
  equivalence relation}. 
For example, we can always take $(X,\mu)=(\{0,1\}^\Gamma,\mu_{eq})$ where 
$\mu_{eq}$ is the infinite product of the measure $(1/2,1/2)$ 
on $\{0,1\}$. 

By Theorems~\ref{thm:suitable-coverings} and~\ref{thm:homotopy
  further} there 
are an $\scrR$-cover~$\scrU$ on~$X\times\widetilde{M}$ and 
a geometric $\scrR$-map
$\phi:X\times\widetilde{M}\rightarrow\scrN(\scrU)$ 
such that for the fundamental class $[M]$ in
$\scrH_n(X\times\widetilde{M})$ we have 
\begin{equation*}
\mass\bigl(\scrH_n(\phi)([M])\bigr)\le\const_{C,n}\vol(M), 
\end{equation*}
where $\const_{C,n}$ only depends on the dimension $n$ and 
the constant $C$ of the packing inequality. 
Because of Lemma~\ref{lem:differential-bounded-by-multiplicity} 
$\phi$ is metrically coarse, and, by construction, $\scrU$ is 
uniformly bounded, and $\scrN(\scrU)$ is finite-dimensional. 
By Lemmas~\ref{lem:homotopy retract} and~\ref{lem:geometric homotopy}
there is a metrically coarse, geometric $\scrR$-map 
$\psi:\scrN(\scrU)\rightarrow X\times\widetilde{M}$ such that 
the composition
$X\times\widetilde{M}\xrightarrow{\phi}\scrN(\scrU)\xrightarrow{\psi}
X\times\widetilde{M}$ is geometrically $\scrR$-homotopic to the identity. 
By Lemma~\ref{lem:norm-decreasing} we have 
\begin{equation*}
\mass\bigl([M]\bigr)=\mass\bigl(\scrH_n(\psi\circ\phi)([M])\bigr)\le 
\mass\bigl(\scrH_n(\phi)([M])\bigr)\le\const_{C,n}\vol(M).  
\end{equation*}
Now Theorem~\ref{thm:folvol bound by packing inequalities} follows 
from Theorem~\ref{thm: l2 betti and mass} in the
Appendix. 
\end{proof}

\section{Proofs of Theorems~\ref{thm:vanishing theorem}
  and~\ref{vanishing theorem for simplicial complexes}}
\label{sec:amenable covers}

\subsection{$\scrR$-covers from the Ornstein-Weiss-Rokhlin
  lemma}\label{subsec:rokhlin lemma} 
The crucial ingredient in the proofs of 
Theorems~\ref{thm:vanishing theorem} 
and~\ref{vanishing theorem for simplicial complexes} 
is the generalized Rokhlin-Lemma of
Ornstein and Weiss. 

\begin{theorem}\label{thm:bulletin version}
\textup{\citelist{\cite{ornstein+weiss(1980)}*{Proposition~4 and 
      Theorem~5}\cite{ornstein+weiss(1987)}*{Theorem~5 in~II.2}}}\hfill\\ 
Retain Assumption~\ref{setup:orbit equivalence relation}. Assume that 
$\Gamma$ is
amenable. Let $\epsilon>0$ and 
$\delta>0$, and let $K\subset\Gamma$ be a finite subset. Then there
exist an $N\in\bbN$, independent of $\delta$ and $K$, 
a sequence of $(K,\delta)$-invariant subsets $H_1,\ldots,H_N$ 
and Borel subsets $B_1,\ldots,B_N\subset X$ such that 
\begin{enumerate}[a)]
\item $\{\gamma B_i;~\gamma\in H_i\}$ are $\epsilon$-disjoint for
  every $i\in\{1,\ldots,N\}$, 
\item the sets $R_i\defq H_iB_i$ are pairwise disjoint, and 
\item $\mu\bigl(\bigcup_{i=1}^NR_i\bigr)>1-\epsilon$. 
\end{enumerate}
\end{theorem}

We say that Borel sets $A_1,\ldots,A_n\subset X$ are \textit{$\epsilon$-disjoint} if
there are pairwise disjoint Borel subsets $A_i'\subset A_i$ with
$\mu(A_i')>(1-\epsilon)\mu(A_i)$. 
For subsets $D,K\subset\Gamma$ we define 
\begin{equation*}
\partial_KD\defq\bigl\{\gamma\in D;~\exists\lambda\in (K\cup
K^{-1}):\lambda\gamma\not\in D\bigr\},
\end{equation*}
and we say that $D$ is 
\textit{$(K,\delta)$-invariant} if 
\begin{equation*}
\frac{\#\partial_KD}{\#D}<\delta. 
\end{equation*}

We need the following modified version of the previous theorem, 
which is nothing new. 
Here $\epsilon$-disjointness 
is replaced by disjointness but $N$ now 
depends on the whole setup. For convenience 
we include a proof. 

\begin{theorem}\label{thm:Rokhlin lemma -- weaker version}
Retain Assumption~\ref{subsec:conventions}. Assume that $\Gamma$ is
amenable. Let $\epsilon>0$ and 
$\delta>0$, and let $K\subset\Gamma$ be a finite subset. Then there
are an $N\in\bbN$ and 
a sequence of $(K,\delta)$-invariant subsets $H_1,\ldots,H_N$ 
and Borel subsets $A_1,\ldots,A_N\subset X$ such that 
\begin{enumerate}[a)]
\item $\{\gamma A_i;~\gamma\in H_i\}$ are disjoint for
  every $i\in\{1,\ldots,N\}$, 
\item the sets $R_i\defq H_iA_i$ are pairwise disjoint, and 
\item $\mu\bigl(\bigcup_{i=1}^NR_i\bigr)>1-\epsilon$. 
\end{enumerate}
\end{theorem}

\begin{proof}
Set $\delta_0\defq\delta/2$ and $\epsilon_0\defq\min\{1/2,\epsilon/2,
\delta_0/(2\#K+1)\}$. 
From applying Theorem~\ref{thm:bulletin version} for the constants 
$\epsilon_0^2$,
$\delta_0$ and $K$ we get 
$(K,\delta_0)$-invariant subsets
$H_1,\ldots,H_{N_0}\subset\Gamma$ and Borel subsets
$B_1,\ldots,B_{N_0}\subset X$. 
For every $i\in\{1,\ldots,N_0\}$
and $\gamma\in H_i$ 
let $B_{i,\gamma}\subset \gamma B_i$ be a Borel subset such that
$\mu(B_{i,\gamma})>(1-\epsilon_0^2)\mu(B_i)$ and the 
$(B_{i,\gamma})_{\gamma\in H_i}$ are pairwise 
disjoint. Set 
\begin{equation*}
S_i\defq\bigcup_{\gamma\in H_i}\{\gamma\}\times \gamma^{-1}B_{i,\gamma}\subset
H_i\times B_i.
\end{equation*}
Then the group action map $m:H_i\times B_i\rightarrow X$ is 
injective on $S_i$ and
\begin{equation}\label{eq:epsilon estimate}
(c\times\mu)(S_i)>(1-\epsilon_0^2)(c\times\mu)(H_i\times B_i),
\end{equation}
where $c$ is the counting measure. For any subset $S\subset\Gamma\times X$ let 
\begin{equation*}
c_S(x)\defq\#\bigl\{\gamma\in\Gamma;~(\gamma,x)\in S\bigr\}. 
\end{equation*}
Now define 
\begin{equation*}
A_i\defq\{x\in B_i;~c_{S_i}(x)>(1-\epsilon_0)\# H_i\}.
\end{equation*}
Next we show that 
\begin{equation}\label{eq:estimate for size of A_i}
\mu(A_i)>(1-\epsilon_0)\mu(B_i). 
\end{equation}
We have the obvious estimate 
\begin{align*}
(c\times\mu)(S_i)&\le\mu(A_i)\#
H_i+(1-\epsilon_0)\bigl(\mu(B_i)-\mu(A_i)\bigr)\# H_i\\
&=(1-\epsilon_0)\mu(B_i)\# H_i+\epsilon_0\mu(A_i)\#H_i, 
\end{align*}
which in combination with~(\ref{eq:epsilon estimate}) 
yields~(\ref{eq:estimate for size of A_i}). 
Let $S_i'=S_i\cap (H_i\times A_i)$. For $x\in A_i$ we 
have $c_{S_i'}(x)=c_{S_i}(x)>(1-\epsilon_0)\# H_i$. Together with 
the injectivity of $m\vert_{S_i'}$ and~(\ref{eq:estimate for size of
  A_i}) we get 
\begin{equation*}
\mu(H_iA_i)\ge (c\times\mu)(S_i')>(1-\epsilon_0)\mu(A_i)\#
H_i>(1-\epsilon_0)^2\mu(B_i)\# H_i\ge(1-\epsilon_0)^2\mu(H_iB_i).
\end{equation*} 
Since the sets $H_iB_i$, thus the sets $R_i\defq H_iA_i$, are pairwise disjoint, 
\begin{equation}
\mu\Bigl(\bigcup_{i=1}^{N_0}H_iA_i\Bigr)>
(1-\epsilon_0)^2\mu\Bigl(\bigcup_{i=1}^{N_0}H_iB_i\Bigr)>
(1-\epsilon_0)^2(1-\epsilon_0^2)>1-2\epsilon_0.  
\end{equation}
If $m(S_i')\ne H_iA_i$, we can enlarge $S_i'\subset H_i\times A_i$ 
while keeping $m\vert_{S_i'}$ injective. So 
we may and will assume that $m(S_i')=H_iA_i$. 
Partition each $A_i$ into finitely many subsets $A_{ij}$ such that 
$(H_i\times\{x\})\cap S'_i$ is a constant set $H_{ij}\subset H_i$ 
for $x\in A_{ij}$. We obtain that 
\begin{enumerate}[a')]
\item on each $A_{ij}\times H_{ij}$ the map $m$ is
  injective, 
\item the sets $\{A_{ij}H_{ij}\}_{i,j}$ are pairwise disjoint, and 
\item $\bigcup_jH_{ij}A_{ij}=H_iA_i$, thus
  $\mu(\bigcup_{i,j}H_{ij}A_{ij})>1-2\epsilon_0\ge 1-\epsilon$.  
\end{enumerate}
From $\#H_{ij}/\#H_i>1-\epsilon_0$ and the $(K,\delta_0)$-invariance 
of $H_i$ easily follows that 
$H_{ij}$ is 
$(K,\delta_0+(2\#K+1)\epsilon_0)$-invariant, thus
$(K,\delta)$-invariant. Now reindexing $(A_{ij})$ and $(H_{ij})$ 
as $A_1,\ldots,A_N$ and $H_1,\ldots,H_N$ gives the sets with the
stated properties. 
\end{proof}

In the sequel we refer to the following setup: 
\begin{setup}\label{setup:for theorems B and C}
Retain the setup of 
Assumption~\ref{setup:framework for simplicial complexes}, 
and assume that $M$ 
is covered by open amenable subsets 
$V(1),\dots, V(m)\subset M$ such that every point in $M$ is contained
in at most $n$ such subsets. We may and will assume that $V(j)$ is  
connected. 
\end{setup}

Write $\incl_j$ for the inclusion $V(j)\hookrightarrow M$. 
By hypothesis, $\Gamma(j)=\im\bigl(\pi_1(\incl_j)\bigr)$ 
is an amenable subgroup of
$\Gamma=\pi_1(M)$. Let $\bar{V}(j)$ be the regular covering of $V(j)$ 
associated to
$\ker\bigl(\pi_1(\incl_j)\bigr)\subset\pi_1\bigl(V(j)\bigr)$. 
It comes with
a free left action of $\Gamma(j)$. Let $p:\widetilde{M}\rightarrow M$ be the
universal covering of $M$. By covering theory, 
we have a pullback diagram:
\begin{equation}\label{eq:pullback}
\xymatrix{\Gamma\times_{\Gamma(j)}\bar{V }(j)\ar@{^(->}[r]\ar[d]&
  \widetilde{M}\ar[d]^p\\ 
V(j)\ar@{^(->}[r]^{\incl_j}&M}
\end{equation}
The upper map is $\Gamma$-equivariant. 
In other words, $p^{-1}(V(j))\subset\widetilde{M}$ decomposes into
connected components, indexed by $\Gamma/\Gamma(j)$, each of which is 
homeomorphic to $\bar{V}(j)$. 

For every Borel map $f:S_1\rightarrow S_2$ 
between standard Borel spaces with countable fibers there is 
a Borel subset $A\subset S_1$ with $f\vert_A$ injective and
$f(A)=f(S_1)$~\cite{kechris}*{corollary~15.2 on p.~89}.  

In particular, 
there is a Borel fundamental
domain $\scrF(j)$ for the $\Gamma(j)$-action on $\bar{V}(j)$, an open subset 
$U(j)\subset\bar{V}(j)$ and a symmetric subset $S(j)\subset\Gamma(j)$ 
such that 
\begin{equation}\label{eq:defining property of S(j)}
\scrF(j)\subset U(j)\subset S(j)\scrF(j).
\end{equation}
Moreover, we define 
\begin{equation}\label{eq:definition of bar(F)}
\bar{\scrF}(j)\defq S(j)\scrF(j). 
\end{equation}
Since $M$ is compact, we can 
take $\scrF(j)$ and $U(j)$ to be relatively compact in
$\widetilde{M}$. Since the
$\Gamma$-action on $\widetilde{M}$ is proper, 
$S(j)$ can be taken to be finite.  

For every $\delta>0$ we now construct a certain $\scrR$-covering 
$\scrU_\delta$ on $X\times\widetilde{M}$. Fix $\delta>0$. 
For every $j\in\{1,\ldots,m\}$ we apply 
Theorem~\ref{thm:Rokhlin lemma -- weaker version} 
to the $\Gamma(j)$-action on $X$ with the constants $\delta$ and 
\begin{equation}
K\defq S(j)^2,~~~\epsilon\defq\delta, 
\end{equation}
and thus obtain subsets $H_1(j),\ldots,H_N(j)$ of $\Gamma(j)$ and 
$A_1(j),\ldots, A_N(j)$ of $X$ such that 
the sets $\gamma A_i(j)$ for $\gamma\in H_i(j)$ and 
$R_i(j)\defq H_i(j)A_i(j)$ for $i\in\{1,\ldots,N\}$ 
are pairwise disjoint and 
each $H_i(j)$ is $(S(j)^2,\delta)$-invariant. 
By taking the maximum over $j\in\{1,\ldots,m\}$, 
we can pick an $N$ that is independent of $j$. This simplifies notation a bit. 
Define 
\begin{equation*}
A_{N+1}(j)\defq X\backslash\bigcup_{i=1}^NR_i(j)\text{ and }
H_{N+1}(j)\defq\{1\}. 
\end{equation*}
Then $\mu(A_{N+1}(j))<\epsilon=\delta$. For every $j\in\{1,\ldots,m\}$ 
the family 
\begin{equation*}
\scrU_\delta(j)\defq\bigl\{\gamma A_i(j)\times\gamma H_i(j)^{-1}U(j)\bigr\}_{1\le i\le
N+1,\gamma\in\Gamma}  
\end{equation*}
is an $\scrR$-cover of  
$X\times\bigl(\Gamma\times_{\Gamma(j)}\bar{V}(j)\bigr)=X\times
p^{-1}\bigl(V(j)\bigr)$. 
\begin{definition}\label{def:covering from ornstein-weiss}
Let $\scrU_\delta$ be the $\scrR$-cover 
of $X\times\widetilde{M}$ given by the union of $\scrR$-covers 
$\scrU_\delta(j)$ for $j\in\{1,\ldots,m\}$. Note here that a union of families is 
indexed by the disjoint union of their index sets. 
\end{definition}

\subsection{Estimating the number of $\scrR$-cells of the
  nerve $\scrN(\scrU_\delta)$}\label{subsec:estimating cells in the nerve}

\begin{theorem}\label{thm:weighted number of cells in nerve for amenabel
  covering}
Retain Assumption~\ref{setup:for theorems B and C}. 
If $k\ge n$, then 
\[\sigma_k\bigl(\scrN(\scrU_\delta)\bigr)\in O(\delta)\text{ for 
  $\delta\rightarrow 0$}.\] 
\end{theorem}

Here the notation means 
that there is a constant $C$ that depends on $k$, $M$ and the cover 
$\{V(j)\}_{1\le j\le m}$ but \textbf{not} on 
$\delta$ or $\scrU_\delta$ 
such that $\sigma_k\bigl(\scrN(\scrU_\delta)\bigr)\le C\delta$ 
for small $\delta>0$. 

\begin{proof}
Let $[m]$ denote the set $\{1,2,\ldots,m\}$. 
By Lemma~\ref{lem:number of cells in the nerve} we have 
\begin{equation*}
\sigma_k(\scrN(\scrU_\delta))=\frac{1}{(k+1)!}\sum_{\substack{(j_1,\ldots,j_{k+1})\in
    [m]^{k+1}\\(i_1,\ldots,i_{k+1})\in [N+1]^{k+1}}}  
\sum_{\substack{\gamma_2,\ldots, \gamma_{k+1}\\\text{as in~(\ref{eq:first subscript})}}}
\mu\bigr(A_{i_1}(j_1)\cap\gamma_2A_{i_2}(j_2)\cap\ldots\cap
\gamma_{k+1}A_{i_{k+1}}(j_{k+1})\bigr)  
\end{equation*}
where the inner sum runs over all
$\gamma_2,\ldots,\gamma_{k+1}\in\Gamma$ such that
\begin{equation}\label{eq:first subscript}
\begin{cases} 
H_{i_1}(j_1)^{-1}U(j_1)\cap\gamma_2H_{i_2}(j_2)^{-1}U(j_2)\cap\ldots\cap\gamma_k
H_{i_k}(j_k)^{-1} U(j_k)\ne\emptyset,\\
(\gamma_{i_r},i_{r},j_r)\ne (\gamma_{i_s},i_{s},j_s)\text{ for $1\le
  r\ne s\le k+1$ and $\gamma_1=1$}.
\end{cases}
\end{equation}
Define 
\begin{equation*}
\Sigma(j_1,\ldots,j_k)\defq\sum_{(i_1,\ldots,i_{k})\in [N+1]^{k}}
\sum_{\substack{\gamma_2,\ldots,\gamma_k,\\h_1,\ldots,h_k\\\text{as
      in~(\ref{eq:second subscript})}}}
\mu\bigr(A_{i_1}(j_1)\cap\gamma_2A_{i_2}(j_2)\ldots\cap
\gamma_{k}A_{i_{k}}(j_{k})\bigr) 
\end{equation*}
where the inner sum runs over all $\gamma_2,\ldots,\gamma_k\in\Gamma$
and $h_1,\ldots, h_k\in\Gamma$ such that $h_l\in H_{i_l}(j_l)$ and  
\begin{equation}\label{eq:second subscript}
\begin{cases}
h_1^{-1}\bar{\scrF}(j_1)\cap\gamma_{2}h_{2}^{-1}\bar{\scrF}(j_2)\cap\ldots\cap
\gamma_{k}h_{k}^{-1}\bar{\scrF}(j_k)\ne\emptyset, \\
(\gamma_{i_r},i_{r},j_r)\ne (\gamma_{i_s},i_{s},j_s)\text{ for $1\le
  r\ne s\le k$ and $\gamma_1=1$}.
\end{cases}
\end{equation}
The sum $\sum\Sigma(j_1,\ldots,j_{k+1})$ over all 
$(j_1,\ldots,j_{k+1})\in [m]^{k+1}$ clearly dominates 
$\sigma_k(\scrN(\scrU_\delta))$ because of~(\ref{eq:defining property
  of S(j)}) and~(\ref{eq:definition of bar(F)}). 
To prove the claim, it thus suffices to prove that 
\begin{equation}\label{eq:epsilon+delta bound}
\Sigma(j_1,\ldots,j_{k+1})\in O(\delta)\text{ for $\delta\rightarrow 0$.}
\end{equation}
for every $(j_1,\ldots,j_{k+1})\in [m]^{n+1}$ provided $k\ge n$. 
We need the following two lemmas to continue. 
\noqed
\end{proof}
\begin{lemma}\label{lem:decreasing k}
Let $k\ge 1$. There is a constant $C>0$ that 
does \textbf{not} depend on 
$\delta$ and $\scrU_\delta$ such that 
\begin{equation*}
\Sigma(j_1,\ldots,j_{k+1})\le C\cdot\Sigma(j_1,\ldots,j_k). 
\end{equation*}
\end{lemma}

\begin{proof}[Proof of lemma]
There is a finite set
$F\subset\Gamma$ such that for all $j,j'\in [m]$ 
\begin{equation*}
\gamma\bar{\scrF}(j)\cap\bar{\scrF}(j')\ne\emptyset\Rightarrow\gamma\in
F.
\end{equation*}
This is clear since $\Gamma$ acts properly on $\widetilde{M}$ and
each $\bar{\scrF}(j)$ is relatively compact. 
Now 
$\Sigma(j_1,\ldots,j_{k+1})$ is bounded by 
\begin{equation*}
\sum_{\substack{(i_1,\ldots,i_{k})\in [N+1]^{k}\\\gamma_l,h_l\text{as
      in~(\ref{eq:third subscript})}}}
\sum_{\substack{f\in F\\i_{k+1}\in [N+1]\\h\in H_{i_{k+1}}}}
\mu\left(A_{i_1}(j_1)\cap\ldots\cap
\gamma_{k}A_{i_{k}}(j_{k})\cap
\gamma_k h_k^{-1}fhA_{i_{k+1}}(j_{k+1})\right).  
\end{equation*}
In the first sum we sum over all $k$-tuples 
$(1,\gamma_2,\ldots,\gamma_k)\in\Gamma^{k+1}$ and $(h_1,\ldots, h_k)\in
H_{i_1}(j_1)\times\ldots\times H_{i_k}(j_k)$ that satisfy 
\begin{equation}\label{eq:third subscript}
\begin{cases}
h_1^{-1}\bar{\scrF}(j_1)\cap\gamma_{2}h_{2}^{-1}\bar{\scrF}(j_2)\cap\ldots
\cap\gamma_{k}h_{k}^{-1}\bar{\scrF}(j_k)
\ne\emptyset,\\
(\gamma_{i_r},i_{r},j_r)\ne (\gamma_{i_s},i_{s},j_s)\text{ for $1\le
  r\ne s\le k$ and $\gamma_1=1$}.
\end{cases}
\end{equation}
The family
$\bigl\{\gamma_k h_k^{-1}fhA_{i_{k+1}}(j_{k+1})\bigr\}_{i_{k+1}\in
  [N+1], h\in H_{i_{k+1}}}$ (for fixed $f,h_k,\gamma_k$) is a Borel partition of
$X$. Hence the claim is true for $C=\#F$. 
\end{proof}
\begin{lemma}\label{lem:k=2}
For $j\in [N+1]$, we have $\Sigma(j,j)\in O(\delta)$ for
$\delta\rightarrow 0$. 
\end{lemma}
\begin{proof}[Proof of lemma]
For $\bar{V}(j)\cap\gamma\bar{V}(j)\ne\emptyset$
it is necessary that $\gamma\in\Gamma(j)$ 
(see diagram~(\ref{eq:pullback})). 
If $h^{-1}\bar{\scrF}(j)\cap\gamma h'^{-1}\bar{\scrF}(j)\ne\emptyset$ for
$h,h'\in\Gamma(j)$ then $\gamma\in h^{-1}S(j)^2h'$ because 
$\scrF(j)$ is a $\Gamma(j)$-fundamental domain and
because of~(\ref{eq:definition of bar(F)}). Thus 
\begin{align*}
\Sigma(j,j)\le & \sum_{\substack{i_1,i_2\in[N+1]\\h_1\in
    H_{i_1}(j)}}\sum_{h_2\in H_{i_2}(j)}
\sum_{\substack{s\in S(j)^2\\h_1\ne sh_2\text{ if $i_1=i_2$}}} \mu\bigl(A_{i_1}(j)\cap
h_1^{-1}sh_2A_{i_2}(j)\bigr)\\ 
=&\sum_{\substack{i_1,i_2\in[N+1]\\h_1\in
    H_{i_1}(j)}}\sum_{h_2\in H_{i_2}(j)}
\sum_{\substack{s\in S(j)^2\\h_1\ne sh_2\text{ if $i_1=i_2$}}} \mu\bigl(h_1A_{i_1}(j)\cap sh_2A_{i_2}(j)\bigr)\\
\intertext{Notice that $\mu(h_1A_{i_1}(j)\cap sh_2A_{i_2}(j))=0$ if $sh_2\in
  H_{i_2}(j)$. Thus, } 
\Sigma(j,j)\le&\sum_{\substack{i_1,i_2\in[N+1]\\h_1\in
    H_{i_1}(j)}}\sum_{h_2\in \partial_{S(j)^2}H_{i_2}(j)}
\sum_{s\in S(j)^2} \mu\bigl(h_1A_{i_1}(j)\cap
sh_2A_{i_2}(j)\bigr)\\
=&\sum_{i_2\in[N+1]}\sum_{h_2\in \partial_{S(j)^2}H_{i_2}(j)}
\sum_{s\in S(j)^2}
\mu\bigl(sh_2A_{i_2}(j)\bigr)\\ 
=&\sum_{i_2\in[N+1]}\sum_{h_2\in \partial_{S(j)^2}H_{i_2}(j)}
\sum_{s\in S(j)^2}
\mu\bigl(H_{i_2}(j)A_{i_2}(j)\bigr)/\#H_{i_2}(j)\\ 
\le &
\# S(j)^2\delta\sum_{i_2\in[N+1]}\mu\bigl(H_{i_2}(j)A_{i_2}(j)\bigr)\\
=&\# S(j)^2\delta.\qedhere
\end{align*}
\end{proof}

\begin{proof}[Continuation of proof of 
Theorem~\ref{thm:weighted number of cells in nerve for amenabel
  covering}]
It remains to verify~(\ref{eq:epsilon+delta bound}). 
If $k\ge n$, then $\Sigma(j_1,\ldots,j_{k+1})=0$ unless 
two of the entries of $(j_1,\ldots, j_{k+1})$ are equal because 
of~(\ref{eq:second subscript}) and the fact that 
$\{V(j)\}_{1\le j\le m}$ has multiplicity at most $n$. 

We can assume that $j=j_1=j_2$ without loss of generality. 
Then by repeated application of Lemma~\ref{lem:decreasing k}, we
obtain that 
$\Sigma(j_1,\ldots,j_{k+1})\le C^{k-1}\Sigma(j,j)$ where $C$ is the
constant from Lemma~\ref{lem:decreasing k},  
and Lemma~\ref{lem:k=2} finally yields~(\ref{eq:epsilon+delta bound}). 
\end{proof}

\subsection{Conclusion of proofs}\label{subsec:conclusion of proofs}
Let $M$ be as in the hypothesis of Theorem~\ref{thm:vanishing
  theorem} or Theorem~\ref{vanishing theorem for simplicial complexes}. 
For $\Gamma=\pi_1(M)$ let $(X,\mu)$ 
be as in Assumption~\ref{setup:orbit equivalence
  relation}. For example, take 
the Bernoulli space $(X,\mu)=(\{0,1\}^\Gamma,\mu_{eq})$ 
where $\mu_{eq}$ is the infinite product of 
the measure $(1/2,1/2)$ on $\{0,1\}$. 
Let $\scrU_\delta$ be the $\scrR$-cover on 
$X\times\widetilde{M}$ constructed in Subsection~\ref{subsec:rokhlin
  lemma}. It is immediate from the construction that $\scrU_\delta$ is
uniformly bounded and $\scrN(\scrU_\delta)$ is finite-dimensional. 

By Lemmas~\ref{lem:map to nerve;finite index set} and~\ref{lem:homotopy
  retract} 
there exists an $\scrR$-simplicial (thus, metrically coarse) map
\begin{equation}\label{eq:R simplicial map}
\phi:X\times\widetilde{M}\rightarrow\scrN(\scrU_\delta)
\end{equation}
(after subdividing the domain).  
By Lemma~\ref{lem:geometric homotopy} there is 
a metrically coarse geometric $\scrR$-map 
\begin{equation}\label{eq:second R simplicial map}
\psi:\scrN(\scrU_\delta)\rightarrow
X\times\widetilde{M}
\end{equation}
such that there is a geometric $\scrR$-homotopy between
$\id_{X\times\widetilde{M}}$ and 
$\psi\circ\phi$. 

\begin{proof}[End of proof of Theorem~\ref{thm:vanishing theorem}]
The triangulation of 
$M$ naturally gives rise to a cycle $\kappa'\in\Cn_n(M)$. 
Let $\kappa$ be the image of $\kappa'$ under 
\begin{equation*}
\Cn(M)=\bbZ\otimes_{\bbZ\Gamma}\Cn(\widetilde{M})\xrightarrow{j_n}\lz\otimes_{\bbZ\Gamma}
\Cn(\widetilde{M})\xrightarrow{\lambda_n}\scrC_n(X\times\widetilde{M}), 
\end{equation*}
where $j_n$ and $\lambda_n$ are the maps from
Remark~\ref{rem:fundamental class in
    tensored complex} and 
Lemma~\ref{lem:inclusion-from-ordinary-to-foliated}, respectively. 
Then $\kappa$ represents the fundamental class $[M]$ in
$\scrH_n(X\times\widetilde{M})$, and $\kappa$ and $\phi$  
satisfy the hypothesis 
of Lemma~\ref{lem:image of simplicial cycle}. 

By Lemma~\ref{lem:image of simplicial cycle} and 
Theorem~\ref{thm:weighted number of cells in nerve for amenabel
  covering} we have 
\begin{equation*}
\mass\bigl(\scrH_n(\phi)([M])\bigr)\le
(n+1)!^2\cdot\sigma_n\bigl(\scrN(\scrU_\delta)\bigr)\in O(\delta). 
\end{equation*}
Lemma~\ref{lem:norm-decreasing} yields 
\begin{equation*}
\mass\bigl([M]\bigr)=\mass\bigl(\scrH(\psi\circ\phi)([M])\bigr)\le 
\mass\bigl(\scrH_n(\phi)([M])\bigr)\in O(\delta).
\end{equation*}
Now Theorem~\ref{thm:vanishing theorem} follows from Theorem~\ref{thm:
  l2 betti and mass} and by letting $\delta\rightarrow 0$. 
\end{proof}

Before we come to the proof of Theorem~\ref{vanishing theorem for
  simplicial complexes}, we need the following lemma that we extract
from Gaboriau's theory of $L^2$-Betti numbers $\betti_k(\Sigma)$ 
for arbitrary $\scrR$-simplicial
complexes~\cite{gaboriau(2002b)}*{Section~4.3}.  

\begin{lemma}\label{lem:abstract lemma from gaboriau theory}
Let $\Sigma,\Psi$ be $\scrR$-simplicial complexes such that 
$\Sigma_x$ is contractible for a.e. $x\in X$. Let 
$\phi:\Sigma\rightarrow\Psi$ be an $\scrR$-simplicial map. Then 
\begin{equation*}
\betti_k(\Sigma)\le\betti_k(\Psi)\text{ for all $k\ge 0$.}
\end{equation*}
\end{lemma}

\begin{proof}
The proof is essentially the same as the one 
of~\cite{gaboriau(2002b)}*{Th\'{e}or\`{e}me~3.13}; there the
hypothesis of the existence of $\phi$ is replaced by 
an assumption on the connectivity of $\Psi$. 
We only indicate the necessary modifications: The map $\phi$ induces 
a \textit{champ bor\'{e}lien \'{e}quivariant}
$t^\ast_x:\Cn_\ast(\Sigma_x)\rightarrow\Cn_\ast(\Psi_x)$ in the sense 
of~\cite{gaboriau(2002b)}*{D\'{e}finition~3.4}. For the 
$t^\ast_x$ we can find 
$s_x^\ast, r_x^\ast$ as in the diagram of~\cite{gaboriau(2002b)}*{Lemme~4.6}. 
The proof~\cite{gaboriau(2002b)}*{Lemme~4.6} almost literally translates to 
our situation. Now the assertion follows 
from~\cite{gaboriau(2002b)}*{Th\'{e}or\`{e}me~4.8}
\end{proof}

\begin{proof}[End of proof of Theorem~\ref{vanishing theorem for
    simplicial complexes}]
Applying Lemma~\ref{lem:abstract lemma from gaboriau theory} 
to~$\Sigma=X\times\widetilde{M}$, $\Psi=\scrN(\scrU_\delta)$ and 
the map in~(\ref{eq:R simplicial map}), we get the estimate 
\begin{equation}\label{eq:l2 betti estimate}
\betti_k\bigl(\widetilde{M}\bigr)=\betti_k\bigl(X\times\widetilde{M}\bigr)\le
\betti_k\bigl(\scrN(\scrU_\delta)\bigr)
\end{equation}
for all~$k\ge 0$. We have
$\betti_k\bigl(\scrN(\scrU_\delta)\bigr)\le\sigma_k\bigl(\scrN(\scrU_\delta)
\bigr)$~\cite{gaboriau(2002b)}*{Proposition~3.2 (3)}. Since 
$\sigma_k\bigl(\scrN(\scrU_\delta)\bigr)\in O(\delta)$ for $k\ge n$ by 
Theorem~\ref{thm:weighted number of cells in nerve for amenabel covering}, 
the assertion follows for~$\delta\rightarrow 0$. 
\end{proof}

\appendix
\section{$\ltwo$-Betti numbers and the support mass of the fundamental class}       
\label{app:l2 betti numbers and foliated volume}

Throughout, 
we retain the setting of Assumption~\ref{setup:framework for
  manifolds}. This appendix is devoted to 
Theorem~\ref{thm: l2 betti and mass} below. 
A statement of this kind was posed in Gromov's
book~\cite{gromov(1999)}*{5.38 e) on p.~304} as an exercise. A proof of 
Theorem~\ref{thm: l2 betti and mass} is given in the doctoral 
thesis~\cite{schmidt}. The proof we present at the end of this
Appendix 
is basically 
that of~\cite{schmidt} and only differs in the 
discussion of equivariant Poincar\'{e} duality. 

\begin{theorem}\label{thm: l2 betti and mass}
Let $[M]$ be the fundamental class in
$\scrH_n(X\times\widetilde{M})$ (cf.~Remark~\ref{rem:fundamental class in
    tensored complex}). Then 
\begin{equation*}
\betti_p(\widetilde{M})\le\mass\bigl([M]\bigr)\text{ for all $p\ge 0$.}
\end{equation*}
\end{theorem}

The proof of Theorem~\ref{thm: l2 betti and mass} 
is essentially an application of equivariant Poincar\'{e} duality. 
The presentation of the latter in~\cite{lueck(2002b)} very much 
fits for our purposes. First we discuss the necessary algebraic 
objects.

\begin{definition_o}
The \textit{twisted group ring} $\lc\rtimes\Gamma$ is, as an $\lc$-module,
free with basis $\Gamma$. Its ring multiplication extends uniquely 
that of $\lc$ and $\bbC\Gamma$ such that 
the commutation rule $f\gamma=\gamma f^\gamma$ holds. 
Recall that $\Gamma$ acts on
$\lz$ from the right by $f^\gamma(x)=(f\gamma)(x)=f(\gamma x)$. 
There is a 
ring involution on $\lc\rtimes\Gamma$ given by 
$\bar{\gamma}=\gamma^{-1}$ where $\bar{f}(x)=\bar{f(x)}$ is complex
conjugation. 

The twisted group ring $\lc\rtimes\Gamma$ is equipped with the trace 
\begin{equation*}
\tr:\lc\rtimes\Gamma\rightarrow\bbC,~ 
\tr\Bigl(\sum_{i=1}^kf_k\gamma_k\Bigr)=\int_X f_1(x)d\mu(x)
\end{equation*}
where
$\gamma_1=1\in\Gamma$, $\gamma_2,\ldots,\gamma_k\in\Gamma\backslash\{1\}$
and $f_1,\ldots,f_k\in\lc$. The corresponding GNS-construction 
defines a von Neumann algebra $L(X,\Gamma)$ with a finite trace,
which is commonly referred to as the \textit{group measure space
  construction}. The \textit{group von Neumann algebra} $L(\Gamma)$ 
is contained in $L(X,\Gamma)$. 
\end{definition_o}

Let  $\Cn_\ast,\Dn_\ast$ denote $\bbZ\Gamma$-chain complexes. Let 
$\Cn^{-\ast}$ be the chain complex whose $p$-th chain module is 
$\hom_{\bbZ\Gamma}(C_{-p},\bbZ\Gamma)$ with the induced
differential. The minus sign causes $\Cn^{-\ast}$ to be chain complex rather
than a cochain complex. Naturally, $\Cn^{-\ast}$ is a \textit{right} 
$\bbZ\Gamma$-module but every right module over a ring with involution
can be viewed as a \textit{left} module. Furthermore, 
$\hom_{\zg}\bigl(\Cn^{-\ast},\Dn_\ast\bigr)$ is the
\textit{hom-complex}; it is again a chain complex whose 
$p$-th chain group consists of degree $p$ chain maps
$\Cn^{-\ast}\rightarrow\Dn_\ast$. Its $p$-th homology consists of
homotopy classes of such. 

For brevity, let us write $\lzo$ instead of $\lz$ in the sequel. 
We have the following commutative square of chain homomorphisms 
\begin{equation*}
\xymatrix{
\bbZ\otimes_{\zg} \bigl(\Cn_\ast\otimes
\Dn_\ast\bigr)\ar[r]^{\phi_u}\ar[d]&\hom_{\zg}\bigl(\Cn^{-\ast},\Dn_\ast
\bigr)\ar[d]\\
\lzo\otimes_{\zg} \bigl(\Cn_\ast\otimes
\Dn_\ast\bigr)\ar[r]^-{\phi_d}&\hom_{\zg}\bigl(\Cn^{-\ast},
\lzog\otimes_{\zg}\Dn_\ast\bigr)}
\end{equation*}
The verticals are the obvious inclusion and induction
$\lzog\otimes_{\zg}\_$ respectively; note here that 
the lower right corner is
canonically isomorphic to
\begin{equation*}
\hom_{\lzo\rtimes\Gamma}(\lzo\rtimes\Gamma\otimes_{\zg}\Cn^{-\ast}, 
\lzog\otimes_{\zg}\Dn_\ast).
\end{equation*} 
The upper homomorphism 
is given by $\phi_u(1\otimes x\otimes y)(g)=\overline{g(x)}y$ for 
$g\in \Cn^{-\ast}$; the lower homomorphism is defined by 
$\phi_d(f\otimes x\otimes y)(g)=\overline{f g(x)}\otimes
y$. \smallskip\\
Let $AW_\ast:\Cn_\ast(\widetilde{M}\times\widetilde{M})\rightarrow
\Cn_\ast(\widetilde{M})\otimes_\bbZ 
\Cn_\ast(\widetilde{M})$ be the \textit{Alexander-Whitney map}, and let 
$\Delta_\ast: \Cn_\ast(\widetilde{M})\rightarrow
\Cn_\ast(\widetilde{M}\times\widetilde{M})$ be the chain map coming from the
diagonal inclusion. Both maps are $\bbZ\Gamma$-equivariant with the 
appropriate diagonal $\Gamma$-action on the target. \smallskip\\

We obtain the following commutative square: 
\begin{equation*}
\xymatrix{\bbZ\otimes_{\bbZ\Gamma}\Cn_\ast(\widetilde{M})\ar[rr]^-{AW_\ast\circ\Delta_\ast}\ar[d]& &
  \bbZ\otimes_{\bbZ\Gamma}\bigl(\Cn_\ast(\widetilde{M})\otimes\Cn_\ast(\widetilde{M})\bigr)\ar[d]\\
\lzo\otimes_{\bbZ\Gamma}\Cn_\ast(\widetilde{M})\ar[rr]^-{\id\otimes
  AW_\ast\circ\Delta_\ast}& &
  \lzo\otimes_{\bbZ\Gamma}\bigl(\Cn_\ast(\widetilde{M})\otimes \Cn_\ast(\widetilde{M})\bigr)}
\end{equation*}

If we set $\Cn_\ast=\Cn_\ast(\widetilde{M})$, concatenate both squares above
and take the $n$-th homology, we obtain the following diagram. 
\begin{equation*}
\xymatrix{
\Hn_n(M)\ar[r]^-{\cap\_}\ar[d]&\left[\Cn^{-\ast}(\widetilde{M}),
  \Cn_{n+\ast}(\widetilde{M})\right]\ar[d]\\ 
\Hn^\Gamma_n(\widetilde{M};\lzo)\ar[r]^-{\cap\_}&\left[\lzo\rtimes\bbZ\Gamma
  \otimes_{\zg}\Cn^{-\ast}(\widetilde{M}),
  \lzo\rtimes\Gamma\otimes_{\zg} \Cn_{n+\ast}(\widetilde{M})\right]}
\end{equation*}
The square brackets denote homotopy classes of chain maps. A
representative of the image of $z$ under one of the horizontal maps 
is a chain map that 
is called the \textit{cap product with $z$}. \smallskip\\
Now consider an element $z\in H_n(M)$ and its image $z'$ in 
$\Hn^\Gamma_n(\widetilde{M};\lzo)$. 
The cap product with $z'$ defines, after induction 
with $L(X,\Gamma)\otimes_{\lzog}\_$, an 
$L(X,\Gamma)$-homomorphism (up to homotopy) 
\begin{equation}\label{eq:cap product}
\_\cap z':~L(X,\Gamma)\otimes_{\zg}\Cn^{-\ast}(\widetilde{M})\longrightarrow
L(X,\Gamma)\otimes_{\zg}\Cn_{n+\ast}(\widetilde{M}).
\end{equation} 
According to equivariant Poincar\'{e} duality,  
the map~(\ref{eq:cap product}) is a $\zg$-homotopy equivalence 
if $z\in H_n(M)$ is the \textit{fundamental class} of $M$ (already 
before induction with $L(X,\Gamma)\otimes_{\lzog}\_$). 
Let us take a closer look at this map. Suppose $z'$ is represented by 
the cycle $\sum_{k=1}^m f_k\otimes \sigma_k$ with $f_k\in\lzo,
\sigma_k:\Delta^n\rightarrow\widetilde{M}$. 
With the standard formula for the 
Alexander-Whitney map, we can unravel the definition of~(\ref{eq:cap
  product}) and see that it 
sends the element $1\otimes g\in
L(X,\Gamma)\otimes_{\zg}\Cn^{j}(\widetilde{M})$ to 
\begin{equation*}
\sum_{k=1}^m\overline{f_kg(\sigma_k\rfloor_j)}\otimes\sigma_k\lfloor_{n-j}\in
L(X,\Gamma)\otimes_{\zg}\Cn_{n-j}(\widetilde{M}). 
\end{equation*}
Here $\sigma\rfloor_j$ and $\sigma\lfloor_{n-j}$ denote 
the front $j$-face and the back $(n-j)$-face of $\sigma$
respectively. 
It follows that the kernel $\ker(ev_j)$ of  
\begin{gather*}
ev_j: L(X,\Gamma)\otimes_{\zg}
\Cn^j(\widetilde{M})\rightarrow\bigoplus_{k=1}^mL(X,\Gamma)\chi_{\supp
  f_k},\\
x\otimes g\mapsto
\bigl(x\overline{f_kg(\sigma_k\rfloor_j)}=x\overline{g(\sigma_k\rfloor_j)}f_k
\bigr)_k 
\end{gather*}
is contained in $\ker(\_\cap z')$; here $\chi_{\supp f_k}$ is the
characteristic function of the support of $f_k$. \smallskip\\
To follow the proof of Theorem~\ref{thm: l2 betti and mass}, 
the reader must be aware of the following facts. 
\begin{enumerate}[a)]
\item There is a \textit{dimension} for every (algebraic)
  module $M$ over a von Neumann algebra $\cala$ with a 
finite trace $\tr_\cala:\cala\rightarrow\bbC$, denoted by
$\dim_\cala(M)\in [0,\infty]$.  
\item The dimension $\dim_\cala$ is additive for short exact sequences
  of $\cala$-modules. 
\item For an idempotent $p\in\cala$, $\dim_\cala(\cala p)=\tr_\cala(p)$. 
\item $\betti_i(\widetilde{M})=\dim_{L(\Gamma)}
\bigl(\Hn_i^\Gamma(\widetilde{M}; L(\Gamma))\bigr)$
\item $\dim_{L(\Gamma)}\bigl(\Hn_i^\Gamma(\widetilde{M};
  L(\Gamma))\bigr)=\dim_{L(X,\Gamma)}\bigl(\Hn_i^\Gamma(\widetilde{M};
  L(X,\Gamma))\bigr)$. 
\end{enumerate}
The items a)-c) belong to the core of L\"uck's dimension
theory for a von Neumann algebra with a finite
trace~\cite{lueck(2002)}*{chapter~6}; d) is a consequence 
of~\cite{sauer(2005)}*{Theorem~2.6 and Theorem~4.3}. 
See also~\cite{sauer(2006)}*{Theorem~6.8}. 

\begin{proof}[Proof of Theorem~\ref{thm: l2 betti and mass}]
Let $\epsilon>0$. By Theorem~\ref{thm:inclusion and support norm} 
the fundamental class $z'$ in $\Hn_n^\Gamma(\widetilde{M};\lzo)$ can be 
represented by a cycle 
\begin{equation*}
\sum_{k=1}^mf_k\otimes\sigma_k\in\lzo\otimes_{\bbZ\Gamma}
\Cn_n(\widetilde{M})
\end{equation*}
such that 
\begin{equation}\label{eq:explicit mass}
\sum_{k=1}^m\mu\bigl(\supp(f_k)\bigr)<\mass([M])+\epsilon
\end{equation}
where $[M]$ denotes the fundamental class in
$\scrH_n(X\times\widetilde{M})$.

Let $\delta^j$ denote the $j$-th differential of 
$L(X,\Gamma)\otimes_{\zg}\Cn^\ast(\widetilde{M})$. There is always 
a projection  
$\ker(\delta^j)\twoheadrightarrow\Hn_\Gamma^j(\widetilde{M};L(X,\Gamma))$. 
If we compose it with the 
homology homomorphism induced by~(\ref{eq:cap product}), we obtain 
a map 
\begin{equation*}
\ker(\delta^j)\twoheadrightarrow\Hn_\Gamma^j\bigl(\widetilde{M};L(X,\Gamma)
\bigr) 
\xrightarrow{\cong} 
\Hn^\Gamma_{n-j}\bigl(\widetilde{M};L(X,\Gamma)\bigr) 
\end{equation*}
that factors over $\ker(ev_j)$ since $\ker(ev_j)$ is contained in
$\ker(\_\cap z')$. Thus 
\begin{equation*}
\betti_{n-j}(\widetilde{M})=
\dim_{L(X,\Gamma)}\bigl(\Hn_{n-j}^\Gamma(\widetilde{M},L(X,\Gamma))\bigr)\le
\dim_{L(X,\Gamma)}\bigl(\ker(\delta^j)/\ker(ev_j)\bigr).
\end{equation*}
Since $\ker(\delta^j)/\ker(ev_j)$ injects into 
\begin{equation*}
\bigoplus_{k=1}^mL(X,\Gamma)\chi_{\supp f_k},
\end{equation*}
we conclude with~(\ref{eq:explicit mass}) that 
\begin{equation*}
\betti_{n-j}(\widetilde{M})\le\sum_{k=1}^m\mu\bigl(\supp(f_k)\bigr)<\mass
\bigl([M]\bigr)+\epsilon.\qedhere
\end{equation*}
\end{proof}

\end{document}